\documentclass[10pt]{article}
\usepackage[latin1]{inputenc}
\usepackage{epsfig}
\usepackage{amsmath}
\usepackage{amssymb}
\usepackage{graphicx}

{\catcode `\@=11 \global\let\AddToReset=\@addtoreset}
\AddToReset{equation}{section}

\newtheorem{@definition}{\sc Definition}[section]

\newtheorem{@remark}{\sc Remark}[section]

\newtheorem{@example}{\sc Example}[section]

\newcommand{\beqn}{\begin{displaymath}}
\newcommand{\eeqn}{\end{displaymath}}
\newcommand{\beq}{\begin{equation}}  
\newcommand{\eeq}{\end{equation}}

\newcommand{\noi}{\noindent}

\def\mathsf{\bf}
\def\N{\mathbb{N}}

\def\R{\mathbb{R}}
\def\Z{\mathbb{Z}}

\def\i{\mathrm i}
\def\d{\mathrm d}
\def\e{\mathrm e}

\def\E{\mathrm E}
\def\P{\mathrm P}

\def\text{\mbox}

\def\1{{\bf 1}}

\newcommand{\mbf}[1]{\mbox{\boldmath $#1$}}

\newcommand{\Var}{\mathrm{Var}\,}
\newcommand{\Cov}{\mbox{Cov}\,}

\def\limiteloiN{\renewcommand{\arraystretch}{0.5}
\begin{array}[t]{c}
\stackrel{{\cal D}}{\longrightarrow} \\
{\scriptstyle N\rightarrow\infty}
\end{array}\renewcommand{\arraystretch}{1}}

\def\limiteloin{\renewcommand{\arraystretch}{0.5}
\begin{array}[t]{c}
\stackrel{{\cal D}}{\longrightarrow} \\
{\scriptstyle n \rightarrow\infty}
\end{array}\renewcommand{\arraystretch}{1}}

\def\limiteproban{\renewcommand{\arraystretch}{0.5}
\begin{array}[t]{c}
\stackrel{{\cal P}}{\longrightarrow} \\
{\scriptstyle n \rightarrow\infty}
\end{array}\renewcommand{\arraystretch}{1}}

\def\limiten{\renewcommand{\arraystretch}{0.5}
\begin{array}[t]{c}
\stackrel{}{\longrightarrow} \\
{\scriptstyle n\rightarrow\infty}
\end{array}\renewcommand{\arraystretch}{1}}

\def\limitepsn{\renewcommand{\arraystretch}{0.5}
\begin{array}[t]{c}
\stackrel{a.s.}{\longrightarrow} \\
{\scriptstyle n\rightarrow\infty}
\end{array}\renewcommand{\arraystretch}{1}}

\def\limitefddl{\renewcommand{\arraystretch}{0.5}
\begin{array}[t]{c}
\stackrel{f.d.d.}{\longrightarrow} \\
{\scriptstyle \lambda \rightarrow 0}
\end{array}\renewcommand{\arraystretch}{1}}

\def\egalfdd{\renewcommand{\arraystretch}{0.5}
\begin{array}[t]{c}
\stackrel{f.d.d.}{=}
\end{array}\renewcommand{\arraystretch}{1}}

\def\egallaw{\renewcommand{\arraystretch}{0.5}
\begin{array}[t]{c}
\stackrel{law}{=}
\end{array}\renewcommand{\arraystretch}{1}}

\def\limitefddn{\renewcommand{\arraystretch}{0.5}
\begin{array}[t]{c}
\stackrel{f.d.d.}{\longrightarrow} \\
{\scriptstyle n\rightarrow\infty}
\end{array}\renewcommand{\arraystretch}{1}}

\def\limitet{\renewcommand{\arraystretch}{0.5}
\begin{array}[t]{c}
\stackrel{}{\longrightarrow} \\
{\scriptstyle t\rightarrow\infty}
\end{array}\renewcommand{\arraystretch}{1}}

\def\limitet0{\renewcommand{\arraystretch}{0.5}
\begin{array}[t]{c}
\stackrel{}{\longrightarrow} \\
{\scriptstyle t\rightarrow 0}
\end{array}\renewcommand{\arraystretch}{1}}

\newtheorem{thm}{Theorem}
\newtheorem{rem}{Remark}
\newtheorem{cor}{Corollary}
\newtheorem{lem}{Lemma}
\newtheorem{prop}{Proposition}

\newtheorem{popy}{Property}

\def\Cov{\mathrm{Cov}\,}
\def\Cor{\mathrm{Cor}\,}

\begin{document}

\title{Nonparametric estimation of the local Hurst function of multifractional Gaussian processes}
\author{Jean-Marc Bardet and Donatas~Surgailis
\footnote{Research of this author was supported in part by the Grant MIP-11155 from the Research Council of Lithuania.}}
\maketitle

\begin{abstract}
A new nonparametric estimator of the local Hurst function of a multifractional Gaussian process based on the increment ratio (IR) statistic is defined. In a general frame, the point-wise and uniform weak and strong
consistency and a multidimensional central limit theorem for this estimator are established. Similar results are obtained
for a  refinement of the generalized quadratic variations (QV)  estimator. The example of the multifractional Brownian motion  is studied
in detail. A simulation study is included  showing that the IR-estimator is more accurate than the QV-estimator.
\end{abstract}
\noindent Keywords \  Nonparametric estimators \ Hurst function \ tangent process \
multifractional Brownian motion \ Gaussian process \ central limit theorem.
\noindent MSC Primary: 62G05 \ Secondary 62G20 \ 60F05 \ 60G22

\section{Introduction}

The present paper discusses nonparametric estimation of the local Hurst function $H(t), t \in (0,1)$ of a multifractional
Gaussian process $Z =   \big( Z(t)\big)_{t \in (0,1)} $ observed at discrete times $k/n, k=1,\cdots, n-1$. The process $Z$ is characterized by
the fact that its infinitesimal increments at each point $t \in (0,1)$ resemble a multiple of a fractional Brownian motion  $B_{H(t)}$
with Hurst parameter $H(t) \in (0,1)$:
\begin{equation}\label{ass}
\big(\lambda^{-H(t)} (Z(t + \lambda \tau) - Z(t))\big)_{\tau \ge 0} \, \limitefddl \, \big(c(t) B_{H(t)}(\tau)\big)_{\tau \ge 0},
\end{equation}
in the sense of the weak convergence of finite-dimensional distributions, where $c(t)>0, t \in (0,1)$ is a deterministic function.
Recall that a fractional Brownian motion  (FBM) with Hurst parameter $H \in (0,1)$ is a Gaussian process $\big(B_{H}(\tau)\big)_{\tau \ge 0}$  having stationary increments and satisfying $\E B_H(\tau)=0$ and $\E B^2_H(\tau)= \tau^{2H}, \tau \ge 0$. 
Property (\ref{ass}) is called asymptotic self-similarity
and the limiting process a tangent process at $t$. For a Gaussian zero-mean process $Z$, relation (\ref{ass}) is equivalent to
\begin{equation}\label{ass1}
\E (Z(t + \tau) - Z(t))^2  \, \sim \, c^2(t) \tau^{2H(t)}, \qquad \tau \to 0.
\end{equation}
A  typical example of such Gaussian process $Z$ is a multifractional Brownian motion (MBM) defined through
the following harmonizable representation
\begin{equation}
X(t)  = \ B_{H(t)}(t) \ =  \ K(H(t)) \,  \int_{\R} \frac{{\e}^{{\i}t x} -
1}{|x|^{H(t)+1/2}} W(\d x), \quad\mbox{with}\quad K(H) :=\big( H\,\Gamma (2H) \sin(\pi \, H)/\pi\big )^{1/2},
   \label{MBMintro}
\end{equation}
where $W(\d x)$ is a complex-valued Gaussian noise satisfying
$ \E \big[\int_{\R} g(x)\,W(\d x) \overline{\int_{\R} h(x)\, W(\d x)} \big] =
\int_{\R} g(x) \overline{h(x)}\, \d x, \,  h, g \in  L^2 (\R)$
 (see more details in section 7.2.2 of \cite{ST})
A first version of  MBM was first introduced by Peltier and L\'evy-V\'ehel \cite{lvmh1995} using a slightly different time-domain representation, while the MBM in (\ref{MBMintro}) was defined, studied and generalized in several papers. In particularly, Stoev and Taqqu
\cite{StoevTaqqu} defined a generalization   $Y_{(a^+,a^-)}$ of (\ref{MBMintro}) depending on parameters $(a^+, a^-) \in \R^2 \backslash \{(0,0)\} $
such that $Y_{(1,1)}$ agrees with $X$ of (\ref{MBMintro}) up to a multiplicative deterministic function,
see Section \ref{Ya} below for details. In the present paper, the
functions $H(\cdot)$ and $c(\cdot)$ are assumed to belong to a H\"older class ${\cal C}^\eta (0,1), \, \eta >0.$
The smoothness
parameter $\eta $ plays an important role in the asymptotic results discussed below. In particular, for $0<\eta < 1$ the local Hurst function of MBM  can be very irregular making the problem of its estimation
more difficult.

\vskip.2cm

Another important parameter for estimation of $H(t)$ is the bandwidth parameter $\alpha \in (0,1)$. Following
a general approach in non-parametric estimation,  discrete observations $Z(k/n), 1\le k \le n-1$
are localized in a neighborhood ${\cal V}_{n,\alpha}(t) = \{ k/n: |k/n -t| \le n^{-\alpha}\} $ of a given point $t\in (0,1)$  containing at most
$[2n^{1-\alpha}]+1 $ points, and an estimate
of $H(t)$ from the localized sample
$Z(k/n),\, k \in  {\cal V}_{n,\alpha}(t)$ is constructed.
Choosing a smaller $\alpha $ means that the estimate of  $H(t)$
uses more
observations of $Z$ and therefore  the variance of it can be expected to be smaller. On the other hand,
increasing the ``width'' of ${\cal V}_{n,\alpha}(t)$ (by taking a smaller $\alpha $) usually  increases the bias of the estimator
particularly
when $\eta $ is small, since the function $H(\cdot)$ is more likely to change in a larger interval ${\cal V}_{n,\alpha}(t)$.
Hence the accuracy of estimation of $H(t)$ depends on $\alpha $ and $\eta$ in a crucial way.

\vskip.2cm
The present paper discusses two discusses two estimators of $H(t)$.  The first estimator $\widehat H^{(I\!R)}_{n,\alpha}(t)$ is a localized version of the increment ratio (IR) estimator introduced in \cite{bs2010}. It is written as a localized sum of ratios of
second-order variations $V_n^{a^*} Z(k/n) = Z((k+2)/n) - 2Z((k+1)/n) + Z(k/n), \, k \in  {\cal V}_{n,\alpha}(t)$ and is related to zero intersections of  the  variations sequence, see (\ref{defHIR}), also  (\cite{bs2010}, Remark 1.1).  We also
discuss a new estimator  $\widehat H^{(I\!R2)}_{n,\alpha}(t)$ which is a pseudo-generalized
least squares (PGLS) version of $\widehat H^{(I\!R)}_{n,\alpha}(t)$ constructed from several dilatations of the second-order variations. At the same time, we study the generalized quadratic variations (QV) estimator $\widehat H^{(QV)}_{n,\alpha}(t)$  of \cite{coeur2005} and its refinement
(a PGLS version) $\widehat H^{(QV2)}_{n,\alpha}(t).$  The main reasons for including the QV estimators in this paper is
the fact that the study and the asymptotic results for IR and QV estimators are very similar (both estimators
are expressed in terms of sums of nonlinear functions of Gaussian vectors of the Hermite rank $2$), and
our desire to compare finite-sample performance of these two classes of estimators.

\vskip.2cm

Let us briefly describe the main results of this paper. Section \ref{sec1} introduces the two main Assumptions {\bf (A)}$_{\kappa}$
and {\bf (B)}$_{\alpha}$ on a Gaussian process $Z$. The first assumption specifies the rate of convergence in a similar asymptotic relation to (\ref{ass1}) and controls the bias of our estimators. The second assumption concerns the decay rate of the covariance
of generalized variations of $Z$ and thus helps to control the variance of the estimators. Under these assumptions,
Theorem \ref{Limitgeneral} (i)-(iv) obtains the weak and strong consistency of the estimators  $\widehat H^{(I\!R)}_{n,\alpha}(t)$ and $\widehat H^{(QV)}_{n,\alpha}(t)$
at a given point $t\in (0,1)$ and some  uniform weak and strong consistency rates of these estimators over arbitrary interval $t \in (\epsilon, 1-\epsilon), \, \epsilon >0$
under increasingly restrictive conditions on $\alpha, \eta, \kappa $ and other parameters. A multidimensional central limit theorem (CLT)  satisfied by the above estimators is proved in
Theorem \ref{Limitgeneral2}. Corollaries \ref{cor1} and \ref{cor2} extend Theorems \ref{Limitgeneral} and \ref{Limitgeneral2} for the PGLS versions of the IR and QV estimators which are studied in Section \ref{Impro}. Proposition \ref{limitMBM} provides limit theorems to the IR and QV estimators for the MBM (including its generalization in \cite{StoevTaqqu}).  Section \ref{simu} is devoted to simulations while Section \ref{proofs} and Section \ref{proofsMBM} contain the proofs. 
Let us note that the simulations in Section \ref{simu} performed on the MBM in (\ref{MBMintro}) for some typical Hurst functions 
and three different values of the smoothness parameter $\eta$ show a very good performance of the estimator
$\widehat H^{(I\!R2)}_{n,\alpha}(t)$. The last estimator provides in almost all cases a smaller MISE criterion than $\widehat H^{(QV2)}_{n,\alpha}(t)$. We also note that the numerical results for both classes of estimators are clearly (in the case of IR) or slightly (in the case of QV) in favor
of their PGLS versions.

\vskip.2cm

Nonparametric estimation of the local Hurst function was the subject of several papers
(\cite{ayache2005}, \cite{ayache2004}, \cite{Begyn2007}, \cite{Begyn2007bis}, \cite{ben2000}, \cite{ben1998a}, \cite{coeur2005}, \cite{istas1997},  \cite{Lacaux2004}). However most of these results are less comprehensive and sharp when applied to the MBM while some contain erroneous statements.
See Remarks \ref{mis}, \ref{remlacaux}, and \ref{rembegyn}, below, for more detailed  discussion and some controversies.
The MBM is probably the most important example of a multifractional Gaussian process at the present time.
It is considered as a benchmark for the estimation problem  of $H(t)$ in the present paper and is  discussed in detail in Section \ref{Ya}.
We believe that the conditions of Proposition \ref{limitMBM} for the IR and QV estimators of the Hurst function of the
MBM are optimal or close to optimal since they are derived from the asymptotic expansion of
the covariance of generalized variations of the MBM in Lemma \ref{lem0}. See also Remark \ref{mis}.

\vskip.2cm

{\it Notation.} In what follows, $C$ stands for a constant whose precise value is unimportant and which may change from line to line. Also, we write
$\limiteproban, \limitepsn, \limiteloin $ for convergence in probability, almost sure (a.s.) convergence and the (weak) convergence
of probability distributions, respectively. $[\eta ]$ denotes the integer part of a real number $\eta$.
${\cal C}^\eta(a,b)$ is the space of all  $[\eta]\ge 0$ times continuously differentiable functions $f: (a,b) \mapsto \R$
on interval $(a,b) \subset \R$ such that $f^{([\eta])}$, the $[\eta]$th derivative of $f$, is an $(\eta-[\eta])$-Hölderian function on $(a,b)$, {\it i.e.} there exists $K\geq 0$ such that  $|f^{([\eta])}(x)-f^{([\eta])}(y)|\leq K \, |x-y|^{\eta-[\eta]}$ for all $x,y \in (a,b)$.


\section{General results}\label{sec1}

\subsection{Definitions of estimators}

Consider a filter $a:=(a_0,\cdots,a_q) \in \R^{q+1}$ such that there exists $m \in \N^*$ satisfying
\begin{eqnarray}\label{moment}
\sum_{\ell=0}^q \ell^p a_\ell=0~~\mbox{for}~ p=0,\cdots,(m-1)~~~\mbox{and}~~~\sum_{\ell=0}^q \ell^{m} a_\ell \neq 0.
\end{eqnarray}
Denote ${\cal A}(m,q)$ the above class of filters. For $a \in {\cal A}(m,q)$ and
$n \in \N^*$ define
the corresponding generalized variations of $Z$ by
\begin{equation} \label{vardef}
V_n^a Z(t):=\sum_{\ell=0}^q a_\ell \, Z(t+\ell/n), \qquad V^a Z(t) \equiv V_1^a Z(t).
\end{equation}
The main examples of generalized variations are the usual simple variations corresponding to  $a=(1,-1)\in {\cal A}(1,1), \, m=q=1$,  and the second order variations corresponding to $a=a^*:=(1,-2,1)\in {\cal A}(2,2), \, m=q=2$. For a filter $a=(a_0,\cdots,a_q)\in {\cal A}(m,q)$ and $j \in \N $
define
its $j$th dilatation $a^{(j)}=(a^{(j)}_0,\cdots,a^{(j)}_{jq})  \in {\cal A}(m,jq) $ by
$$
a^{(j)}_{ij}:=a_i~~\mbox{and}~~a^{(j)}_{k}:=0~~\mbox{if}~k \notin j \N.
$$
For $\alpha\in (0,1)$, define a neighborhood of $t\in (0,1)$ and its cardinal by:
$${\cal V}_{n,\alpha}(t):=\Big \{k \in \{1,2,\cdots,n-q-1\},\, |\frac kn-t|\leq n^{-\alpha} \Big \}~~\mbox{and}~~v_{n,\alpha}(t):=\# \,{\cal V}_{n,\alpha}(t). $$
Note that for any $t\in (0,1)$ and $n$ large enough, $v_{n,\alpha}(t)=[2n^{1-\alpha}]$ or $[2n^{1-\alpha}]+1$ depending on the parity of $ [2n^{1-\alpha}]$ = the integer part of
$2n^{1-\alpha}$.
For $H \in (0,1)$, define
\begin{eqnarray}
\label{Lambda2} \Lambda_2(H)&:=&\E \psi\big(V_n^{a^*} B_H(0),V_n^{a^*} B_H(\frac 1 n)\big), 
\end{eqnarray}
where
\begin{equation}\label{psi}
\psi(x^{(1)}, x^{(2)}) := \frac{|x^{(1)}  + x^{(2)}|}{|x^{(1)}|  + |x^{(2)}|}, \qquad (x^{(1)}, x^{(2)}) \in \R^2
\end{equation}
is a function taking values in $[0,1]$, more explicitly,
\begin{eqnarray}
\Lambda_2(H)&=&
\frac{1}{\pi}\arccos (-\rho_2 (H)) +
\frac{1}{\pi}
\sqrt{\frac{1+\rho_2 (H)}{1-\rho_2 (H)}} \log \left(\frac{2}{1+\rho_2 (H)}\right), \label{Lambda2}\\
\label{rho2}
\rho_2(H)&:=&{\rm Cor}(V_n^{a^*} B_H(0), V_n^{a^*} B_H(\frac 1 n)) \ = \  \frac{-3^{2H}+2^{2H+2}-7}{8-2^{2H+1}}.
\end{eqnarray}
The function $\Lambda_2$ does not depend on $n$ by self-similarity of FBM and is monotone increasing on $(0,1)$, see \cite{bs2010}.
The estimators $\widehat H^{(QV)}_{n,\alpha}(t)$ and $\widehat H^{(I\!R)}_{n,\alpha}(t)$ are defined for $t\in (0,1)$ from
observed sample $\big(Z(1/n), Z(2/n),\cdots,Z((n-1)/n)\big)$ as follows:
\begin{eqnarray}\label{defHQV}
\widehat H^{(QV)}_{n,\alpha}(t)&\hspace{-2mm}:=&\hspace{-2mm}\frac 1 2 \frac {A}{A A^\intercal} \, \Big ( \log \Big (\frac 1 {v_{n,\alpha}(t)}
\sum_{k\in {\cal V}_{n,\alpha}(t)} \big( V_n^{a^{(i)}} Z(\frac k n)\big )^2\Big )  \Big )^\intercal_{1\leq i\le  p}, \\
\label{defHIR}
\widehat H^{(I\!R)}_{n,\alpha}(t)&\hspace{-2mm}:= &\hspace{-2mm}\Lambda_2^{-1} \Big (\frac 1 {v_{n,\alpha}(t)}
\sum_{k\in {\cal V}_{n,\alpha}(t)}
\psi \big(V_n^{a^*} Z(\frac k n), V_n^{a^*} Z(\frac {k+1} n)\big) \Big),
\end{eqnarray}
where $A :=\big (\log i -\frac 1 p \sum_{j=1}^p \log j \big )_{1\leq i\leq p} \in \R^p$ is a row vector, $A^{\intercal}$ its transposed vector
(vector-column),
and $\Lambda_2^{-1}$ is the inverse function of $\Lambda_2 $.
Thus, $\widehat H^{(QV)}_{n,\alpha}(t)$ is defined from a log-regression of generalized quadratic variations, see also \cite{coeur2005}, whereas $\widehat H^{(I\!R)}_{n,\alpha}(t)$ is defined from a sample mean of ratios of the second-order variations in a neighborhood of a given point $t \in (0,1)$. A discrete time version of the last estimator was first studied in \cite{surt2008}.

\subsection{Assumptions}

First, recall (see e.g. \cite{surt2008}) that  for a filter $a=(a_0,\cdots,a_q)\in {\cal A}(m,q)$ with $m\geq 1$ and for any fixed $H\in (0,1)$ and $n \in \N$,
$\left(V^{a}_n  B_{H}(j/n) \right)_{j \in \Z}$ is a stationary
Gaussian process,  $\big(V^{a}_n  B_{H}(\frac j n) \big)_{j \in \Z}
\egalfdd n^{-H}\big(V^{a} B_{H}(j) \big)_{j \in \Z} $ by self-similarity, furthermore,
\begin{eqnarray}
\nonumber {\rm Var}\big(V_n^{a^{(i)}}  B_{H} (j)\big)&\hspace{-3mm} =&
\hspace{-3mm} \frac{i^{2H}}{n^{2H}} C(H,a), \quad
C(H,a):= {\rm Var}\big(V^{a}  B_{H} (0)\big) \, = \
-\frac 1 2 \sum_{k,k'= 0}^q a_k a_{k'} |k-k|^{2H}, \quad i \in \N^*, \\
{\rm Cov}\big(V_n^{a}  B_{H}(\frac k n), V^{a}  B_{H}(\frac {k'} n )\big)
\nonumber &\hspace{-3mm}=&\hspace{-3mm}-\frac 1 {2\, n^{2H}} \,  \sum_{i=0}^q \sum_{j=0}^q a_ia_j \big |k'-k+j-i \big |^{2H}  \\
\label{equivar} &\hspace{-3mm}\sim&\hspace{-3mm} \Big ( (-1)^{m+1}\frac {\big (\sum_{i=0}^q a_i i^m\big )^2}{2 (m!)^2} \,  \prod_{\ell=0}^{2m-1} (2H-\ell) \Big )\, \frac {|k'-k|^{2H-2m}}{n^{2H}}, \qquad |k'-k| \to \infty.
\end{eqnarray}
Note $C(H,a) >0$ for any $H \in (0,1), a \in {\cal A}(m,q) $ which fact follows from the harmonizable
representation of FBM analogous to (\ref{MBMintro}) yielding $C(H,a) =
K^2(H) \int_{\R} |\sum_{k=0}^q a_k (\e^{\i k x} -1)|^2 |x|^{-2H-1} \d x >0$ since the integrand does not
identically vanish on $\R$.
In the particular case $a=a^*=(1,-2,1)$ $(m=2)$, then
\begin{eqnarray*}
{\rm Var}\big(V^{a^{*}}  B_{H} (j)\big)&=&4- 4^H, \\
 {\rm Cor}\big(V^{a^{*}}  B_{H}(j), V^{a^{*}}  B_{H}(j+1)\big) &=& \rho_2(H),
\end{eqnarray*}
where $\rho_2(H)$ is defined in (\ref{rho2}).
~\\

For a given $m \ge 1$, a general Gaussian process  $ Z= (Z(t))_{t \in (0,1)}$
and real numbers $\kappa\geq 0$ and $\alpha \in (0,1)$, we introduce the following assumptions:

\begin{enumerate}
\item [{\bf (A)}$_{\kappa}$] There exist continuous functions $0<H(t) <1$ and
$c(t)>0$ for $t \in (0,1)$ satisfying $H (\cdot) \in {\cal C}^\eta(0,1)$ and $c (\cdot) \in {\cal C}^{\eta}(0,1), \ \eta >0$ such that for any $0<\varepsilon <1/2,\, j\ge 1$ and any filter $a \in {\cal A}(m,q), q\geq 1$,
\begin{eqnarray}
\label{prop12}
\max_{[n\varepsilon] \le k \leq [(1-\varepsilon)n]} n^\kappa \,
\Big| \frac{\Var \big (V_n^{a} Z(k/n)  \big ) }
{\Var \big (V_n^a B_{H(k/n)}(0)\big ) } - c(\frac k n )\Big|\limiten 0.
\end{eqnarray}

\smallskip

\item [{\bf (B)}$_{\alpha}$]
There exist  $C>0, \  \gamma > 1/2 $ and $0\le \theta < \gamma/2$
such that for any  $n\in {\N}^*$, $0<\varepsilon<1/2, \, [n\varepsilon]\le k <k' \leq [(1-\varepsilon)n]$
and any filter $a \in {\cal A}(m,q), q\geq 1$,
\begin{equation*}
\Big | {\rm Cor} \big (V_n^{a} Z(\frac k n),  V_n^{a} Z(\frac {k'} n)\big )
\Big | \  \le  \ C\, n^{(1-\alpha)\theta} (|k'-k|\wedge n^{1-\alpha})^{-\gamma}.
\end{equation*}
\end{enumerate}

\medskip

\begin{rem} \label{tangent} {\rm
If condition {\bf (A)$_0$} is satisfied for $m=1$
then $Z$ admits a FBM as a tangent process at each point $t \in (0,1)$, in the sense that
for any $k \in \Z$,
\begin{eqnarray} \label{tangentB}
n^{H(t)}\Big(Z(\frac {[nt] + k+ j } n) - Z(\frac {[nt] +k} n)\Big)_{j \in \Z} &\limitefddn&\left(\sqrt{c(t)}
B_{H(t)}(j) \right)_{j \in \Z}.
\end{eqnarray}
Assumption  {\bf (A)}$_{\kappa}$ with $\kappa >0$  allows to control the uniform convergence rate of generalized variations of $Z$ to the corresponding
variations of the tangent process. Condition {\bf (B)$_{\alpha}$} provides a  bound for the correlation of the process
$\big(V_n^a  Z(k/n)\big)_{[n\varepsilon]\le k, k' \leq [(1-\varepsilon)n]} $ (the choice of this bound can be deduced from the expansion (\ref{equivar})). Under such conditions, the strong consistency and a CLT for the estimators
$\widehat H^{(QV)}_{n,\alpha},
\widehat H^{(I\!R)}_{n,\alpha}$ respectively defined
in (\ref{defHQV}), (\ref{defHIR}) can be established (see Sec. 2.3 below).

}

\end{rem}

\noindent  The following statement deduced from Assumption {\bf (A)$_\kappa$} will be often used
in the remaining sections.

\begin{popy} \label{popy1} Assumption {\bf (A)$_\kappa$} implies that for any  $a \in {\cal A}(m,q)$ and any $\varepsilon >0, \, j \in \N$
\begin{eqnarray}\label{prop11}
&&\max_{[n\varepsilon] \le k \leq [(1-\varepsilon)n]}
n^{\kappa} \Big| {\rm Cor}\big (V_n^{a} Z(\frac k n ), V_n^{a} Z(\frac {k+j} n) \big )
- {\rm Cor}\big (V^{a}  B_{H(\frac k n)}(0), V^{a} B_{H(\frac k n)}(u)\big ) \Big| \limiten 0. 
\end{eqnarray}
\end{popy}

\noi {\it Proof.} We have
\begin{eqnarray*}
\E \big[V_n^{a} Z(\frac k n )V_n^{a} Z(\frac {k+j} n)\big]
&=&\frac{1}{2}\big\{\Var \big(V_n^{a} Z(\frac k n ) +  V_n^{a} Z(\frac {k+j} n )\big)
- \Var \big(V_n^{a} Z(\frac k n ) -  V_n^{a} Z(\frac {k+j} n )\big)\big\}\\
&=& \frac{1}{2}\big\{\Var \big(V_n^{a^+} Z(\frac k n )\big)
- \Var \big(V_n^{a^-} Z(\frac k n )\big)\big\},
\end{eqnarray*}
where generalized variations $V_n^{a^\pm} Z(\frac k n )
= \sum_{i=0}^{j+q} a^\pm_i Z(\frac {k+i} n)$ correspond to filters  $a^\pm  = (a_0^\pm, \cdots, a^\pm_{q + j}) \in \R^{ q + j+1}$ with
\begin{equation}
a^+_i := a_i \1(0 \le i \le q) + a_{i-j} \1(j \le i \le j+q), \qquad  a^-_i := a_i \1(0 \le i \le q) - a_{i-j} \1(j \le i \le j+q).
\end{equation}
Note that $a \in {\cal A}(m,q)$
implies $a^\pm \in  {\cal A}(m,j+ q)$ for any $j =0,1,\cdots.$ From
 {\bf (A)$_\kappa$} with $a = a^\pm$ we have that
uniformly in $[n\varepsilon] \le k \leq [(1-\varepsilon)n],$
$
V_n^{a^\pm} Z(\frac k n ) = c(\frac k n) \Var V_n^{a^\pm} B_{H(k/n)}(0) \big(1 + o(n^{-\kappa})\big),$
implying
\begin{eqnarray}
\Cov \big(V_n^{a} Z(\frac k n ),V_n^{a} Z(\frac {k+j} n)\big)&=&\frac 1 2 c(\frac k n)
\big\{\Var (V_n^{a^+} B_{H(k/n)}(0)) \big(1 + o(n^{-\kappa})\big) - \Var (V_n^{a^-} B_{H(k/n)}(0)) \big(1 + o(n^{-\kappa})\big)\big\}
\nonumber \\
&=&\frac 1 2 c(\frac k n)
\big\{\Var (V_n^{a^+} B_{H(k/n)}(0))  - \Var (V_n^{a^-} B_{H(k/n)}(0))\big\} +  o(n^{-2H(k/n) -\kappa}) \nonumber \\
&=&c(\frac k n){\rm Cov}\big(V_n^{a}  B_{H(k/n)}(0), V^{a}_n  B_{H(k/n)}(\frac {j} n )\big) + o(n^{-2H(k/n) -\kappa}).\label{CovVa}
\end{eqnarray}
We also have 
\begin{equation}\label{VarVa}
V_n^{a} Z(\frac k n ) = c(\frac k n) \Var (V_n^{a} B_{H(k/n)}(0)) \big(1 + o(n^{-\kappa})\big), \qquad  
V_n^{a} Z(\frac {k+j} n ) = c(\frac k n) \Var (V_n^{a} B_{H(k/n)}(0)) \big(1 + o(n^{-\kappa})\big), 
\end{equation}
where the last relation follows by writing  
$\Var \big(V_n^{a} Z(\frac {k+j} n)\big)=\Var \big(V_n^{a^+} Z(\frac {k} n) \big) - \Var \big(V_n^{a} Z(\frac {k} n)\big)
- 2\Cov \big(V_n^{a} Z(\frac k n ),V_n^{a} Z(\frac {k+j} n)\big)$ and using the facts established above. 
Then, (\ref{prop11}) is deduced from (\ref{CovVa}) and (\ref{VarVa}). \hfill $\Box$

\subsection{Asymptotic results}

Denote
\begin{equation}\label{mudef}
\mu := \frac{2(1-\alpha)(1 - \frac{\theta}{\gamma}) - \alpha}{4}, \qquad
\mu_1 := \frac{2(1-\alpha)(1 - \frac{\theta}{\gamma}) - \alpha -1}{4}.
\end{equation}
Note $\mu >0 $ (respectively, $\mu_1 >0$) is equivalent to $ \alpha < \frac {2\gamma - 2\theta}{3\gamma - 2\theta} $
(respectively,  $\alpha < \frac {\gamma - 2\theta}{3\gamma - 2\theta})$.

\medskip

\begin{thm}\label{Limitgeneral}

Let $Z= (Z(t))_{t \in (0,1)}$ be a zero-mean Gaussian process.
\smallskip

\noindent (i) If $Z$ satisfies condition {\bf (A)}$_{0}$  and  {\bf (B)$_{\alpha}$} with $m=2, \,
\alpha \in (0,1)$ then for any $t \in (0,1)$
$$
\widehat H^{(I\!R)}_{n,\alpha}(t) \limiteproban H(t).
$$

\noindent (ii) If $Z$ satisfies conditions {\bf (A)}$_{0}$  and  {\bf (B)$_{\alpha}$} with $m=2, \,
0< \alpha < \frac {\gamma - 2 \theta}{2(\gamma - \theta)}$
then for any $t \in (0,1)$
$$
\widehat H^{(I\!R)}_{n,\alpha}(t) \limitepsn H(t).
$$

\noindent (iii) If $Z$ satisfies conditions {\bf (A)}$_{{\kappa}}$ and   {\bf (B)$_{\alpha}$} with  $m=2, \, \kappa \ge \mu$ and
$ \frac {2\gamma - 2\theta}{3\gamma - 2\theta + 4\gamma (\eta \wedge 2)}
   \le  \alpha < \frac {2\gamma - 2\theta}{3\gamma - 2\theta} $ \,
then for any $\epsilon >0$
$$
\sup_{\epsilon <  t <  1-\epsilon}
\big|\widehat H^{(I\!R)}_{n,\alpha}(t) - H(t)\big| \ = \  O_p(n^{-\mu}).
$$

\noindent (iv) If $Z$ satisfies conditions {\bf (A)}$_{{\kappa}}$ and   {\bf (B)$_{\alpha}$} with   $m=2, \, \kappa \ge \mu_1$ and
$ \frac {\gamma - 2\theta}{3\gamma - 2\theta  + 4\gamma (\eta \wedge 2)}    \le  \alpha < \frac {\gamma - 2\theta}{3\gamma - 2\theta} $
then for any $\epsilon >0, \, \delta >0$
$$
\sup_{\epsilon <  t <  1-\epsilon}
\big|\widehat H^{(I\!R)}_{n,\alpha}(t) - H(t)\big| \ = \  O(n^{-(\mu_1-\delta) }) \quad\text{a.s.}
$$

\end{thm}


\begin{thm}\label{Limitgeneral2}
Let $Z= (Z(t))_{t \in (0,1)}$ be a zero-mean Gaussian process satisfying {\bf (A)}$_{{\kappa}}$ and {\bf (B)$_{\alpha}$}, with $m=2, \,
\alpha > \frac{1}{1+ 2 (\eta \wedge 2)}$, $\kappa \ge \frac{1-\alpha}{2}$ and
$\theta=0$.
Then for any $u \in \N^*$ and any $t_1, \cdots, t_u  \in (0,1), \, t_i \ne t_j \ (i \ne j)$
\begin{eqnarray}
\label{TLCIR}
\sqrt{2n^{1-\alpha}} \big (\widehat H^{(I\!R)}_{n,\alpha}(t_i)-H(t_i)\big )_{1\le i \le u}
&\limiteloin&\big (W^{(I\!R)}(t_i)\big )_{1\le i \le u},
\end{eqnarray}
where $W^{(I\!R)}(t_i), \, i=1, \cdots, u$ are independent centered Gaussian r.v.'s
such as, with
$\Lambda_2(\cdot)$ defined in (\ref{Lambda2}),
\begin{eqnarray}\label{IRvar}
&&\E [W^{(I\!R)}(t_i)]^2:=\Big [ \frac
{\partial}{\partial x}  (\Lambda_2)^{-1}(\Lambda_2(H(t_i)))\Big
]^2 \, \sigma^2(H(t_i)) \\
\label{sigma2}
\mbox{where} && \sigma^2(H):= \sum_{k\in {\Z}} \Cov\left( \psi \big(V^{a^*} B_{H}(0), V^{a^*} B_{H}(1)\big),
\psi \big(V^{a^*} B_{H}(k), V^{a^*} B_{H}(k+1)\big) \right)\quad\mbox{for $H\in (0,1)$}.
\end{eqnarray}
\end{thm}

The following Theorem 3 extends Theorems \ref{Limitgeneral} and \ref{Limitgeneral2} to
the estimator $\widehat H^{(QV)}_{n,\alpha}$. As noted in the Introduction, both estimators
$\widehat H^{(I\!R)}_{n,\alpha}(t)$ and $\widehat H^{(QV)}_{n,\alpha}$ are written in terms of sums of nonlinear functions of Gaussian vectors of the Hermite rank $2$ and their asymptotic analysis is very similar.
The proof of Theorem 3  can be deduced
{\em mutatis mutandis}
from the proofs of Theorems \ref{Limitgeneral} and \ref{Limitgeneral2} and is omitted.

\vskip.2cm

\begin{thm}
Let $Z$ be a zero-mean Gaussian process  and $a \in {\cal A}(m,q), \, m \ge 1$ a filter.  ~\\
\noindent (a)
The statements of Theorem \ref{Limitgeneral} (i)-(iv)
 hold for
$\widehat H^{(QV)}_{n,\alpha}(t)$ instead of $\widehat H^{(I\!R)}_{n,\alpha}(t)$
under respective assumptions (i)-(iv) with $m = 2 $ replaced by
$m\ge 1 $.

\smallskip

\noindent (b) Let $Z$ satisfy the assumptions of Theorem \ref{Limitgeneral2} with $m = 2 $ replaced by
$m\ge 1 $.  Then
for any $u \in \N^*$ and any $t_1, \cdots, t_u  \in (0,1), \, t_i \ne t_j \ (i \ne j)$
\begin{eqnarray}
\label{TLCQV}
\sqrt{2n^{1-\alpha}} \big (\widehat H^{(QV)}_{n,\alpha}(t_i)-H(t_i)\big )_{1\le i \le u}
&\limiteloin&\big (W^{(QV)}(t_i)\big )_{1\le i \le u},
\end{eqnarray}
where $W^{(QV)}(t_i), \, i=1, \cdots, u$ are independent centered Gaussian r.v.'s such as, with $ A$ defined in (\ref{defHQV}),
\begin{eqnarray}\label{QVvar}
\E [W^{(QV)}(t_i)]^2:=\frac {A^\intercal\, \Gamma(H(t_i))A }{4(A^\intercal A)^2}~~\mbox{and}~~
\Gamma(H):=\Big (  \frac 2 {i_1^{2H} i_2^{2H}} \,   \sum_{j\in {\Z}}\Big [\frac {\sum_{k_1,k_2=0}^q a^{(1)}_{k_1}a^{(1)}_{k_2}  |i_1k_1-i_2k_2+j|^{2H}}{\sum_{k_1,k_2=0}^q a^{(1)}_{k_1}a^{(1)}_{k_2}|k_1-k_2|^{2H}} \Big ] ^2\Big )_{1 \leq i_1,i_2 \leq p}.
\end{eqnarray}
\end{thm}

\section{Pseudo-generalized least squares versions of the IR and QV estimators}
\label{Impro}

Asymptotic and finite-sample performance of the estimators $\widehat H^{(I\!R)}_{n,\alpha}(t)$ and $\widehat H^{(QV)}_{n,\alpha}(t)$
can be improved by using their pseudo-generalized
least squares (PGLS) versions $\widehat H^{(I\!R2)}_{n,\alpha}(t)$ and  $\widehat H^{(QV2)}_{n,\alpha}(t)$  as defined below. \\
~\\
{\it Definition of $\widehat H^{(I\!R2)}_{n,\alpha}(t)$.}  Let $a^{(i*)} $ denote the $i$th dilatation of the filter $a^{*} = (1, -2, 1), \, i\in \N^*,$ so that
$V_n^{a^{(i*)}} Z(\ell/n)=Z((\ell+2i)/n)-2Z((\ell+i)/n)+Z(\ell/n)$. The corresponding
IR estimator of $H(t)$ is defined as
$$
\widehat H_{n,\alpha}^{(I\!R),i}(t):=\big [\Lambda_2^{(i)}\big ]^{-1} \Big ( \frac 1 {v_{n,\alpha}(t)}
\sum_{k\in {\cal V}_{n,\alpha}(t)}
\psi \big(V_n^{a^{(i*)}} Z(\frac k n), V_n^{a^{(i*)}} Z(\frac {k+1} n)\big) \Big),
$$
where the function $\Lambda_2^{(i)}: (0,1) \to (0,1)$ is defined as the expectation
\begin{eqnarray*}
\Lambda_2^{(i)}(H)
 &\hspace{-3mm} :=&\hspace{-3mm}\E \psi \big(V_n^{a^{(i*)}} B_H(0), V_n^{a^{(i*)}} B_H(\frac 1 n)\big)  =
\frac{1}{\pi}\arccos (-\rho_2 ^{(i)}(H)) +
\frac{1}{\pi}
\sqrt{\frac{1+\rho_2^{(i)} (H)}{1-\rho_2^{(i)} (H)}} \log \left(\frac{2}{1+\rho_2^{(i)} (H)}\right), \\
\label{rhoi}
\mbox{and}\qquad \rho_2^{(i)}(H)&\hspace{-3mm}:=&\hspace{-3mm}{\rm Cor}(V_n^{a^{(i*)}} B_H(0), V_n^{a^{(i*)}} B_H(\frac 1 n))  =  \frac{-|2i+1|^{2H}-|2i-1|^{2H}+4|i+1|^{2H}+4|i-1|^{2H}-6}{8-2^{2H+1}}.
\end{eqnarray*}
Under the assumptions of Theorem  \ref{Limitgeneral2}, these estimators satisfy the $p-$dimensional CLT:
\begin{eqnarray}\label{TLCHmulti}
n^{(1-\alpha)/2} \left (\widehat H_{n,\alpha}^{(I\!R),i}(t)-H(t)  \right)_{1\leq i \leq p}
&\limiteloin&{\cal N}_p \left  (0\, ,\,\Sigma^{(p)}(H(t)) \right )
\end{eqnarray}
with the limit covariance matrix $\Sigma^{(p)}(H(t))$, where for $H\in (0,1)$
\begin{eqnarray}\label{Sigmap}
\Sigma^{(p)}(H)&:=&\Big (  \Big [ \frac
{\partial}{\partial x}  (\Lambda_2^{(i)})^{-1}(\Lambda_2^{(i)}(H))\Big
] \, \sigma_{ij}(H) \, \Big [ \frac
{\partial}{\partial x}  (\Lambda_2^{(j)})^{-1}(\Lambda_2^{(j)}(H))\Big
]\Big) _{1\leq i,j \leq p}, \\
\label{sigmaij}
\sigma_{ij}(H)&:=&\sum_{k\in {\Z}} \Cov\left(
\frac {\big | V^{a^{(i*)}} B_{H}(0)+V^{a^{(i*)}} B_{H}(1)\big |}{\big |
V^{a^{(i*)}} B_{H}(0)\big |+\big |V^{a^{(i*)}} B_{H}(1)\big |}, \frac {\big | V^{a^{(j*)}} B_{H}(k)+V^{a^{(j*)}} B_{H}(k+1)\big |}{\big |
V^{a^{(j*)}} B_{H}(k)\big |+\big |V^{a^{(j*)}} B_{H}(k+1)\big |} \right)
\end{eqnarray}
(note that $V^{a^{(j*)}} B_{H}(k)=B_H(2j+k)-2B_H(j+k)+B_H(k)$ for $j=1,\cdots,p$ and $k\in \Z$). Now define
\begin{equation} \label{SigmaGLS}
\widehat \Sigma^{(p)}(t)\ := \ \Sigma^{(p)}\big (\widehat H^{(I\!R)}_{n,\alpha}(t) \big ),
\end{equation}
which is a consistent estimator of $\Sigma^{(p)}(H(t))$ since $\widehat H^{(I\!R)}_{n,\alpha}(t)$ is a consistent estimator of $H(t)$.
Then, with the vector-column $\1_p=(1,1,\cdots,1)^\intercal \in \R^p $, the PGLSE $\widehat H^{(I\!R2)}_{n,\alpha}(t)$ of $H(t)$ is defined by
\begin{eqnarray}\label{DefHIR2}
\widehat H^{(I\!R2)}_{n,\alpha}(t):= \big (\1_p^\intercal  \,  (\widehat \Sigma^{(p)}(t))^{-1} \, \1_p\big )^{-1}\,
\1_p^\intercal \,  (\widehat \Sigma^{(p)}(t))^{-1}  \, \Big ( \widehat H_{n,\alpha}^{(I\!R),i}(t)
\Big )_{1 \leq i \leq p}.
\end{eqnarray}
~\\
{\it Definition of $\widehat H^{(QV2)}_{n,\alpha}(t)$.}
Analogously  to (\ref{SigmaGLS}), define
\begin{eqnarray}\label{matrixGLS}
\widehat \Gamma(t)&:=&\Gamma \big (\widehat H^{(QV)}_{n,\alpha}(t)\big ),
\end{eqnarray}
which is a consistent estimator of $\Gamma(H(t))$ since $\widehat H^{(QV)}_{n,\alpha}(t)$ is a consistent estimator of $H(t)$.
The one-dimensional CLT  satisfied by $\widehat H^{(QV)}_{n,\alpha}(t)$ at a given point $t\in (0,1)$
(see  (\ref{TLCQV}))
is a consequence of the following $p-$dimensional CLT:
$$
\sqrt{2n^{1-\alpha}}  \Big (\widehat T_{n,\alpha}^{(p)}(t)-2H(t)\, \big ( \log(i/n) \big )_{1\leq i\le  p} -C(t)\1_p   \Big) \limiteloin {\cal N}_p \big (0\, ,\,  \Gamma
(H(t)) \big )
$$
with $C(t):=\log\Big (-\frac {c(t)} {2 } \, \sum_{k,k'=0}^q a^{(1)}_{k}a^{(1)}_{k'}|k'-k|^{2H(t)}\Big ) $
and
$$
\widehat T_{n,\alpha}^{(p)}(t):= \Big ( \log\Big (\frac 1 {v_{n,\alpha}(t)}
\sum_{k\in {\cal V}_{n,\alpha}(t)}\big |V_n^{a^{(i)}} Z(k/n)\big |^2 \Big ) \Big )_{1\leq i\le  p};
$$
see \cite{BSMBM} for details.
Then $2\widehat H^{(QV)}_{n,\alpha}(t)$ is the slope of the ordinary least squares regression of the vector $\widehat T_{n,\alpha}^{(p)}(t)$ onto the vector $\big ( \log(i/n) \big )_{1\leq i\le  p}$. The PGLS estimators $\widehat H^{(QV2)}_{n,\alpha}(t)$ and $\widehat C^{(QV2)}_{n,\alpha}(t) $ of $H(t)$ and $C(t)$ are obtained from the corresponding generalized least squares regression with the covariance matrix $\Gamma(H(t))$
replaced by its estimate  $\widehat \Gamma (t) $ (\ref{matrixGLS}):
\begin{eqnarray}\label{defd1}
\left ( \begin{array}{c} 2 \widehat H^{(QV2)}_{n,\alpha}(t)\\
\widehat C^{(QV2)}_{n,\alpha}(t) \end{array} \right )
&:=&\left (Z^{(p)}_n\, \widehat \Gamma(t)^{-1} \, (Z^{(p)}_n)^\intercal\right )^{-1} Z^{(p)}_n\, \widehat \Gamma(t)^{-1} \,  (\widehat T_{n,\alpha}^{(p)}(t))^\intercal, \\
\text{with} \qquad Z^{(p)}_n&:=&\left (\begin{array}{cccc}  \log (1/n) &\log (2/n) &\dots  &\log (p/n)\\
1 &1 &\dots &1
\end{array} \right).
\end{eqnarray}

\begin{cor} \label{cor1} (i)
Assume that for each $0< H < 1$ the matrix $\Sigma^{(p)}(H)$ of (\ref{Sigmap})
is non-degenerate and continuously depends on $H$.
Then all statements of Theorem \ref{Limitgeneral} (i)-(iv) hold for the estimator $\widehat H^{(I\!R2)}_{n,\alpha}$
instead of $\widehat H^{(I\!R)}_{n,\alpha}$

\smallskip

\noindent (ii) Assume that for each $0< H < 1$ the matrix  $\Gamma (H)$ in (\ref{QVvar})
is non-degenerate and continuously depends on $H$.
Then the statement of Theorem 3 (a) holds for the estimator
$\widehat H^{(QV2)}_{n,\alpha}$
instead of $\widehat H^{(QV)}_{n,\alpha}$.

\end{cor}

\begin{cor}\label{cor2} Let $Z$ satisfy the assumptions of Theorem \ref{Limitgeneral2} and  the matrices $\Sigma^{(p)}(H), \, \Gamma(H), \,
0< H < 1$ satisfy the conditions of Corollary \ref{cor1} (i)-(ii). Then for any $u \in \N^*$ and any $t_1, \cdots, t_u  \in (0,1), \, t_i \ne t_j \ (i \ne j)$,
\begin{eqnarray}
\label{TLCIR2}
\sqrt{2n^{1-\alpha}} \big (\widehat H^{(I\!R2)}_{n,\alpha}(t_i)-H(t_i)\big )_{1\le i \le u}
&\limiteloin&\big (W^{(I\!R2)}(t_i)\big )_{1\le i \le u}
\end{eqnarray}
and
\begin{eqnarray}
\label{TLCQV2}
\sqrt{2n^{1-\alpha}}  \left ( \begin{array}{c} 2 (\widehat H^{(QV2)}_{n,\alpha}(t_i)-H(t_i)) \\
\widehat C^{(QV2)}_{n,\alpha}(t_i) -C(t_i)\end{array} \right )_{1\le i \le u}   & \limiteloin & \big (W^{(QV2)}(t_i)\big )_{1\le i \le u}.
\end{eqnarray}
In (\ref{TLCIR2}), $W^{(QV2)}(t_i), \, 1\leq i \leq u $ are independent centered Gaussian r.v.'s with respective variances $\displaystyle \E [W^{(I\!R2)}(t_i)]^2:= \big(\1_p^\intercal  \,  (\Sigma^{(p)}(H(t_i)))^{-1} \, \1_p\big)^{-1}$.
In (\ref{TLCQV2}), $W^{(QV2)}(t_i), \, 1\leq i \leq u$ are independent centered Gaussian random vectors with respective $2\times 2$
covariance matrices
$\E \big [W^{(QV2)}(t_i)\,W^{(QV2)}(t_i)^\intercal  \big ] := \big (Z_1^{(p)} \, \Gamma(H(t_i))^{-1} \, (Z_1^{(p)})^\intercal\big )^{-1}$.

\end{cor}

\section{The case of multifractional Brownian motion}\label{Ya}

This section details Assumptions {\bf (A)}$_{{\kappa}}$ and   {\bf (B)$_{\alpha}$} for the MBM (\ref{MBMintro}) and
its generalization due to Stoev and Taqqu \cite{StoevTaqqu}.
As a consequence, the asymptotic behavior of the IR and QV estimators for the above process is established
under explicit conditions on the parameters $\alpha, \eta $ and the  function $H(\cdot)$.

The above mentioned generalization of the MBM in \cite{StoevTaqqu}
is given by stochastic integral representation
\begin{eqnarray}\label{Y0}
Y_{(a^+,a^-)}(t):= \frac{ \Gamma (H(t)+ \frac 1 2) }{\sqrt{2\pi}} \,  \int_{\R} \frac{{\e}^{{\i}t x} -
1}{|x|^{H(t)+ \frac 1 2}}\, \, U_{(a^+,a^-)}(H(t),x) \,  W(\d x), \qquad t \in \R,
\end{eqnarray}
where $(a^+,a^-)\in \R^2\setminus \{(0,0)\}$ are parameters, $W (\d x)$ is the same white noise as in (\ref{MBMintro}), and
\begin{equation}
U_{(a^+,a^-)}(H,x):= a^+ \e^{-\i \, \mbox{sign}(x)(H +\frac 1 2)\frac \pi 2} + a^- \e^{\i \, \mbox{sign}(x)(H+\frac 1 2)\frac \pi 2}, \qquad
H \in (0,1), \quad x \in \R.
\end{equation}
The above definition implies  $\E Y^2_{(a^+,a^-)}(t) =  \E Y^2_{(a^+,a^-)}(1) |t|^{2H(t)}$ and a covariance structure
of $Y_{(a^+,a^-)} $ similar to that of the MBM in (\ref{MBMintro}) (see \cite{StoevTaqqu}, Theorem 4.1). In particular,
 $Y_{(1,0)}$ gives a harmonizable representation of the time-domain MBM of \cite{lvmh1995}.
The so-called well-balanced case
case $a^+ = a^- = a \ne 0$, or
\begin{eqnarray*}
Y_{(a,a)}(t)&=&\frac{2a }{\sqrt{2\pi}} \cos \big((H(t) + \frac 1 2) \frac{\pi}{2}\big)  \Gamma (H(t) + \frac 1 2)   \int_{\R} \frac{{\e}^{{\i}t x} -
1}{|x|^{H(t)+ \frac 1 2}}\,  W(\d x)
\end{eqnarray*}
is more delicate, see (\cite{StoevTaqqu}, sec. 5),
since the function $H \mapsto \cos \big((H + \frac 1 2) \frac{\pi}{2}\big) $ vanishes and changes its sign at $H=1/2 $.
In particular, $Y_{(a,a)}(t)/(\E  Y^2_{(a,a)}(1))^{1/2}$ $ = -{\rm sign}(H(t) - \frac{1}{2}) X(t)$, where
$X $ is defined at (\ref{MBMintro}). Since $X$ has a.s. continuous trajectories under the condition (\ref{cond2})
on $H(\cdot)$,
see \cite{ayache2005},
this implies that the normalized process  $Y_{(a,a)}(t)/(\E  Y^2_{(a,a)}(1))^{1/2}$ is generally a.s. discontinuous at each
$t \ne 0$ with $H(t) = 1/2$.

\vskip.1cm

In order to avoid the above complication, we introduce the following
generalization $X_{(a^+,a^-)}$ of the MBM   in (\ref{MBMintro}):
\begin{eqnarray*}\label{Y}
X_{(a^+,a^-)}(t)&:=&\left\{\begin{array}{ll}X(t), &a^+ = a^- \ne 0, \\
Y_{(a^+, a^-)}(t)/{\rm Var}^{1/2}(Y_{(a^+,a^-)}(1)), &a^+ \ne a^-,
\end{array} \right.,
\end{eqnarray*}
where $X$ and $Y_{(a^+, a^-)}$ are defined at (\ref{MBMintro}) and (\ref{Y0}), respectively. Accordingly, the process
$X_{(a^+,a^-)}$ is defined by the following harmonizable representation:
\begin{eqnarray}\label{Y1}
X_{(a^+,a^-)}(t)&=&K(H(t)) \,  \int_{\R} \frac{{\e}^{{\i}t x} -
1}{|x|^{H(t)+ \frac 1 2}}\, \, {\cal U}_{(a^+,a^-)}(H(t),x) \,  W(\d x),  \qquad t \in (0,1),
\end{eqnarray}
where $K(H)$ is defined in (\ref{MBMintro}) and
\begin{equation}
{\cal U}_{(a^+,a^-)}(H,x):= \left \{ \begin{array}{ll}
\displaystyle 1,&a^+=a^- \ne 0,\\
\displaystyle \frac  {a^+ \e^{-\i \, \mbox{sign}(x)(H +\frac 1 2)\frac \pi 2} + a^- \e^{\i \, \mbox{sign}(x)(H+\frac 1 2)\frac \pi 2}} {\big ((a^+)^2+ (a^-)^2 -2 a^+a^- \sin(\pi \, H)\big )^{1/2} },  &a^+ \ne a^-.
\end{array} \right .
\end{equation}
Note that the function $H\mapsto {\cal U}_{(a^+,a^-)}(H,x)$ is analytic on $(0,1)$ for all $x\in \R$ and $(a^+,a^-)\in \R^2\setminus \{(0,0)\}$.
This definition of $X_{(a^+,a^-)}$ induces that for all $t\in (0,1)$ and  $(a^+,a^-)\in \R^2\setminus \{(0,0)\}$
\begin{eqnarray*}
\E X_{(a^+,a^-)}(t)=0\qquad\mbox{and}\qquad \E X^2_{(a^+,a^-)}(t)= t^{2H(t)}.
\end{eqnarray*}
The following Assumption  {\bf (C)}$_{\eta}$ is crucial for regularity properties of the MBM, see \cite{ayache2005}, \cite{ayache2004}.

\vskip.2cm

\noindent {\bf (C)}$_{\eta}$  \, There exists $\eta >0$ such that $H (\cdot) \in {\cal C}^\eta(0,1)$  and
\begin{eqnarray*}\label{cond2}
0 < \inf_{t\in (0,1)} H(t) \le  \sup_{t\in (0,1)} H(t) < \min(1,\eta).
\end{eqnarray*}

\noindent It is known that Assumption {\bf (C)}$_{\eta}$ guarantees that the
MBM $X=X_{(1,1)}$ of (\ref{MBMintro}) is locally asymptotically self-similar at each point $t \in (0,1)$ having a FBM $B_{H(t)}$ as its tangent
process at $t$ \cite{ben1997} and its pointwise Hölder exponent coincides with $H(t)$ \cite{ayache2004}.
A particular case ($a=(1,0,\cdots,0,-1)$ and $m=1$) of the following proposition shows
that a FBM $B_{H(t)}$ is also a tangent process at $t$ for $ X_{(a^+,a^-)}$ defined in (\ref{Y1}), see Remark \ref{tangent}
and Property \ref{popy1}.

\begin{prop}\label{propMBM1} Let  $H(\cdot)$ satisfy Assumption {\bf (C)}$_{\eta}$. Then, for any $a \in {\cal A}(m,q)$, with $m\geq 1$, $(a^+,a^-)\in \R^2\setminus \{(0,0)\}$  and $0<\varepsilon<1/2$, there exists a constant $C_\ell(\varepsilon) >0$ such that for $n \ge 1$,
\begin{eqnarray} \label{A1'}
\max _{[n\varepsilon]\le k \leq [(1-\varepsilon)n]} \, \Big |
 \frac{\displaystyle  \Var \big ( V_n^a X_{(a^+,a^-)}(\frac k n) \big)} 
{\displaystyle
\Var \big ( V_n^a B_{H(\frac k n)}(\frac k n)\big)
} - 1 \Big |\ \leq \ C_\ell (\varepsilon)\,
\Big (\frac { \log n} {n ^{(\eta \wedge 1)  }}+\frac { 1} {n ^{2((\eta\wedge m) -\sup_{t\in (0,1)} H(t))  }} \Big ).
\end{eqnarray}
\end{prop}

The next proposition provides an expansion of the covariance of the process
$V_n^{a^*}X_{(a^+,a^-)}$. For ease of writing, we consider here the case $a=a^*,   m=2$ only. The case of general filter $a \in {\cal A}(m,q)$ is discussed in Lemma \ref{lem0} (see Section \ref{proofsMBM}).

\begin{prop}\label{propMBM2}
Let  $H(\cdot)$ satisfy Assumption {\bf (C)}$_{\eta}$ and  $(a^+,a^-)\in \R^2\setminus \{(0,0)\}$.
Then for any $0<\varepsilon<1/2$, $n \in \N^*$, $[n \varepsilon] \le k <k' \leq [(1-\varepsilon)n]$, $k'-k>2q$,
\begin{multline}
 \Cov  \Big (V_n^{a^*}X_{(a^+,a^-)}(\frac k n), V_n^{a^*} X_{(a^+,a^-)}(\frac {k'} n)\Big)
\ = \   A\big (\frac k n ,\frac {k'}n \big )\, V_n^{a^*} H(\frac k n) \, V_n^{a^*}H(\frac {k'} n)\\
\hskip-2cm + \ B_2\big(\frac k n,\frac {k'} n\big)
 \frac {|k-k'|^{H(\frac k n) + H(\frac {k'}n) -4}}{n^{H(\frac k n) + H(\frac {k'} n)}}   +  \lambda_n(k,k'),  \label{lexp}
\end{multline}
where (recall) $V^{a^*} H (\frac k n)=H (\frac {k+2i} n)-2H (\frac {k+i} n)+H (\frac {k} n)$ are the second-order variations of $H(\cdot)$ and
$A(t,t'),\, B_2(t,t') $ are
defined in (\ref{Att}-\ref{Btt}). 
The remainder term  $\lambda_n(k,k') $ in (\ref{lexp}) satisfies the following bound:
for any $\delta>0$ there exist $n_0 \in \N^*, \, k_0\in \N^*$ such that for any $n > n_0,\, |k -k'| > k_0, \,  [n\varepsilon] \le  k ,k' \leq [(1-\varepsilon)n]$
\begin{eqnarray}  \label{lexp12}
|\lambda_n(k,k')|
&\leq&\delta \, \Big(\frac {1}{n^{2(\eta \wedge 2)}} +  \frac {|k'-k|^{H(\frac k n) + H(\frac {k'}n) -4}}
{n^{H(\frac k n) + H(\frac {k'} n)}} \Big).
\end{eqnarray}
\end{prop}

\begin{rem}\label{mis} {\rm From definitions  (\ref{Att}-\ref{Btt}) it immediately follows that
$A(t,t')$ and $B_2(t,t')$
are bounded on $(0,1) \times (0,1)$ and have finite and generally non-vanishing limits
\begin{eqnarray*}
\lim_{t'\to t} A(t,t')
=t^{2H(t)}\log^2 t  \qquad\mbox{and} \qquad
\lim_{t'\to t} B_2(t,t')
=- \frac {1} 2 \, \prod_{\ell=0}^{3} (2H(t)-\ell ).
\end{eqnarray*}
Thus, Proposition \ref{propMBM2} and Lemma \ref{lem0} allow to obtain lower bounds in (\ref{lexp})
and a decorrelation rate  of generalized variations of the MBM, see Corollary \ref{corolMBM}, (\ref{corMBM}),
which is close to optimal. In turn,  the decorrelation rate of generalized variations determines the
decay rate of the  variance and the 4th moment of the estimators written in terms of sums of nonlinear functions of Hermite rank 2
of normalized generalized variations (see the proof of Theorem   \ref{Limitgeneral2}, also (\ref{4ii}), (\ref{4iii})).
Proposition \ref{propMBM2} and Lemma \ref{lem0} also permit
to construct explicit
counter-examples to some earlier results pertaining to this issue.  In particular, we can show that
the expansions obtained in (\cite{coeur2005}, Lemma 1) and  (\cite{ben1998a}, Lemma 2) are erroneous. 
A detailed discussion of the above counter-examples is given in the extended version of this paper \cite{BSMBM}. Let us note also
that the Erratum in \cite{coeur2006} does not concern (\cite{coeur2005}, Lemma 1) but another error in the last paper.
}
\end{rem}

\begin{rem}{\rm Let us
note that in  \cite{Lacaux2004} and \cite{Begyn2005}, some correct bounds of the decorrelation rate of $ V^{a^*}_n X (\frac{k}{n})$ are obtained, but their bounds are  less sharp than the bound in (\ref{corMBM}). Another correct bound was also obtained in the preprint (\cite{ayache-elnouty2004}, Lemma 2.4). }
\end{rem}
Propositions \ref{propMBM1} and \ref{propMBM2} help to verify Assumptions {\bf (A)$_\kappa$} and {\bf (B)$_{\alpha}$} for the MBM
in (\ref{Y1}).

\begin{cor} \label{corolMBM} 
Let the conditions of Proposition \ref{propMBM2} hold.
Then for any $0<\varepsilon<1/2$, there exist
$C(\varepsilon)>0$ and $n_0(\varepsilon) \in \N^*$ such that for any $n \ge n_0(\varepsilon) $ and any $[n\varepsilon]\le  k <k' \leq [(1-\varepsilon)n]$,
\begin{equation} \label{corMBM}
\Big | {\rm Cor} \big (V_n^{a^*} X_{(a^+,a^-)}(\frac k n),  V_n^{a^*} X_{(a^+,a^-)}(\frac {k'} n)\big )
\Big | \  \le  \ C(\varepsilon)\, \Big(  \frac 1 {|k'-k|^{4-H(\frac k n) - H(\frac {k'} n)}} +\frac {1}{n^{2(\eta \wedge 2)-H(\frac k n) - H(\frac {k'} n)}} \Big).
\end{equation}
Moreover,

\smallskip

\noindent (i) $X_{(a^+,a^-)}$ satisfies Assumption {\bf (A)$_\kappa$} with $m=2$ and any $0\le \kappa < \min \big \{(\eta \wedge 1)\, , \, 2(\eta-\sup_{t\in (0,1)}H(t)) \big \}$;

\smallskip

\noindent (ii) 
$X_{(a^+,a^-)}$ satisfies  Assumption {\bf (B)$_{\alpha}$} with $m=2$ and any $0\le \theta < \gamma/2, \, 2 \ge \gamma > 1/2 $ satisfying the following condition:
\begin{equation}\label{iineq}
\gamma - \theta \ \le \  \min \big( \frac{2(\eta \wedge 2) - 2 H(t)}{1-\alpha}, \, 4 - 2H(t)\big), \qquad \text{for any} \quad t \in (0,1).
\end{equation}

\end{cor}

\medskip

\begin{rem} \label{tvFBM}
{\rm  The presence of the second terms on the r.h.s. of (\ref{A1'}) and (\ref{corMBM}) indicates that dependence properties of increments of the MBM are quite sensitive to the smoothness parameter $\eta $ of the Hurst function. These terms
have a negative effect on estimation of $H(\cdot)$ and the convergence rates,
by imposing restrictions on the bandwidth $\alpha$, see Proposition \ref{limitMBM} below.
\cite{sur2008} argued that the dependence properties of the MBM are rather peculiar and proposed a
different class of multifractional Gaussian processes
defined via nonhomogeneous fractional integration of white noise. These processes are locally self-similar in the sense of
(\ref{ass}) under less restrictive conditions on $H(\cdot)$ and
have better decorrelation properties than the MBM \cite{sur2008}. On the other hand,
the covariance
function 
of the above-mentioned time-varying fractionally integrated processes is not a local function
of $H(t)$ and $H(t')$ as in the case of the MBM (see (\ref{covY})) and
its study is more difficult. Extending Propositions  \ref{propMBM1}
and \ref{propMBM2} to time-varying fractionally integrated processes of \cite{sur2008}
is an interesting open problem.

 }
 \end{rem}


\begin{prop} \label{limitMBM}
Let $X_{(a^+,a^-)}$ be the MBM of \eqref{Y1} with $(a^+,a^-)\in \R^2\setminus \{(0,0)\}$ and $H(\cdot)$ satisfying Assumption
{\bf (C)}$_{\eta}, \, \eta >0$. Then, with $a=a^*$ and ${\cal E}=I\!R,\, I\!R2, \, QV$ or $QV2$,
\smallskip

\noindent (i) \,\,  For any $t \in (0,1)$ and any $\alpha \in (0,1), \, $ $
\widehat H^{({\cal E})}_{n,\alpha}(t) \limiteproban H(t).$


\noindent (ii) \, For any $t \in (0,1), \, $
$\widehat H^{({\cal E})}_{n,\alpha}(t) \limitepsn H(t),$ provided \ $H(t) < \eta - \frac{1}{8} $ and
$\max(0, 1 - 4((\eta \wedge 2) - H(t))) < \alpha < 1/2$  hold.

\noindent (iii) For any $\epsilon >0$, \,
$\sup_{\epsilon <  t < 1-\epsilon}
\big|\widehat H^{({\cal E})}_{n,\alpha}(t) - H(t)\big| \ = \  O_p(n^{-(2-3\alpha)/4}), $ provided $\alpha $ and $\eta $ satisfy
\begin{eqnarray}\label{MBM1}
&&\alpha \ >\   \sup_{t\in (0,1)} (1 - 4((\eta \wedge 2) - H(t))), \\
 \label{MBM2}
&&\sup_{t \in (0,1)} H(t)\, <\, \eta - \frac{1}{12} \qquad \mbox{and} \qquad \frac{2}{3 + 4(\eta \wedge 2)}\ \le \  \alpha\ <\ \frac{2}{3}.
\end{eqnarray}

\noindent (iv) For any $\epsilon >0 $ and  $\delta >0, $
$\sup_{\epsilon < t < 1-\epsilon} \big|\widehat H^{({\cal E})}_{n,\alpha}(t) - H(t)\big| \ = \  O(n^{-(1-3\alpha - \delta)/4}) \ \text{a.s.},$ provided $\alpha$ and $\eta$ satisfy (\ref{MBM1}),
\begin{equation} \label{MBM3}
\sup_{t \in (0,1)} H(t) < \eta - \frac{1}{6} \qquad \text{and} \qquad      \frac{1}{3 + 4(\eta \wedge 2)} \ \le \  \alpha < \frac{1}{3}.
\end{equation}

\noindent (v) For any $u \in \N^*$ and $t_1, \cdots, t_u  \in (0,1), \, t_i \ne t_j \ (i \ne j)$ the multidimensional CLTs   (\ref{TLCIR}), (\ref{TLCQV}),  (\ref{TLCIR2})  and (\ref{TLCQV2}) hold, provided  the matrix-valued functions $\Sigma^{(p)}$ and $\Gamma$ satisfy the condition of Corollary \ref{cor1} and the parameters $\alpha$ and $\eta$ satisfy (\ref{MBM1}) and
\begin{equation}\label{MBM4}
\alpha\ >\  \max \Big\{ \frac{1}{1 + 2(\eta \wedge 2)}\, , \, 1 - 4\big ((\eta \wedge 2) - \sup_{t\in (0,1)} H(t)\big )\Big\}.
\end{equation}

\end{prop}

\medskip

\begin{rem} \label{}
 {\rm Roughly speaking, conditions in (ii) - (v) require that $\eta $ and $H(t)$ are sufficiently separated, or the difference $\eta - H(t)$
is large enough. Else, it may happen that $\alpha $ satisfying these conditions does not exist.
However, if $\eta \ge (3 + \sqrt{41})/8 \simeq 1.175 $ then $\alpha $  satisfying (ii) - (iv) exists
for any $H(\cdot) \in {\cal C}^\eta(0,1)$ with $\sup_{t \in (0,1)} H(t) < 1 $, while  condition \eqref{MBM4} in (v) reduces to
$\alpha >  \frac{1}{1 + 2(\eta \wedge 2)}$.
 }
 \end{rem}

\begin{rem} {\rm  In the case of estimators ${\cal E}= QV$ or $QV2$,
Proposition \ref{limitMBM} can be easily extended to  a general filter $a \in {\cal A}(m,q), \, m \ge 1 $.  For $m=1$,
the corresponding results in (ii) - (v) can be proved under weaker conditions on $\alpha, \eta $ but the
convergence rates are worse.
On the other hand, if $a \in {\cal A}(m,q), \, m\geq 3$, the corresponding
conditions on $\alpha, \eta $  in (ii) - (iv)  
are the same as in the case $a=a^*$ except that $\eta \wedge 2$ in \eqref{MBM1} and \eqref{MBM3} can be replaced by $\eta \wedge m$. In (v), condition  \eqref{MBM4} holds. Note also that this proposition also holds for a process
$(\sigma(t)\, X_{(a^+,a^-)}(t))_{t\in (0,1)}$ when $\sigma \in {\cal C}^\eta(0,1)$. }
\end{rem}

\begin{rem} \label{remlacaux}
{\rm For a class of harmonizable L\'evy processes,  Lacaux  (\cite{Lacaux2004}, Theorem 4.1) obtained the weak consistency rate of the
QV estimator
$\widehat H^{(QV)}_{n,\alpha}(t)-H(t)=O_{p}(n^{-\rho} \log n)$ with
$\rho \le \eta - H(t), \rho < \frac \eta {1+2\eta}$ in the case when  $H(\cdot ) \in C^\eta (0,1), \, 0<\eta \le 1.$
For the MBM the above result follows from the CLT in
Proposition \ref{limitMBM} (v) with $\rho =  \frac{1-\alpha}{2} $ satisfying $\rho < 2(\eta - H(t))$ and $\rho < \frac \eta {1+2\eta}. $
This means that concerning the weak consistency of the point-wise QV estimator of $H(t)$, our results for
the MBM are more accurate than those in \cite{Lacaux2004}. Moreover, the last paper  does not
discuss strong and uniform consistency rates and the CLT as in Proposition \ref{limitMBM}. }

\end{rem}

\begin{rem} \label{rembegyn}
 {\rm
Begyn \cite{Begyn2007} obtained strong consistency and asymptotic normality of sums of quadratic variations for general non-stationary Gaussian processes
under different assumptions on the covariance function which exclude the case of MBM. In particular, if $\eta < 2 $ then the covariance
$R(t,t') = \E X(t) X(t') $ in (\ref{covY}) is {\it not} twice differentiable in $t$ or $t'$ 
and the derivative $\partial^4 R(t,t')/\partial^2 t \, \partial t'^2 $ for $t \ne t'$ does not exist,
contrary to what is assumed
in (\cite{Begyn2007}, Theorem 1 and 2). In the case of quadratic variations localized in a neighborhood ${\cal V}_{n,\alpha}(t)$ of given point
$t \in (0,1)$ as in the estimator $\widehat H^{(QV)}_{n,\alpha}(t)$ (\ref{defHQV}),
the imposed conditions in (\cite{Begyn2007}, eqs. (5), (21), (23)) exclude the appearance
of the first term of the order $O(n^{-2(\eta \wedge 2)})$ on the r.h.s. of (\ref{lexp})  which cannot be ignored and which plays a crucial role in
our asymptotic results.


}
\end{rem}
\begin{rem} \label{pratic}
 {\rm Since $\alpha$ depends on $\eta$ and $H(t)$ that are not available, it is a problem to apply the CLTs of Proposition \ref{limitMBM} in concrete situations. However, when $\eta \geq 2$ it is possible to select any $\alpha>1/5$ and even to use $n^{1/5}\log^2 n$ instead of $n^\alpha$; then the convergence rate of the CLT is $n^{2/5}$  up to a logarithm term. Note that  \cite{coeur2005} proposed a procedure by minimization of the MISE to select an optimal $\widehat \alpha$ and then an adaptive estimator of $H(t)$ is $\widehat H^{(QV)}_{n,\widehat \alpha}(t)$.
 }
 \end{rem}
\section{Simulations}\label{simu}
All the softwares  used in this Section are available with a free access on {\tt
http://samm.univ-paris1.fr/-Jean-Marc-Bardet} (in Matlab language).\\
We use the original version of the MBM $X$ in \eqref{MBMintro}. 
Since trajectories of the MBM in our simulation study are generated (for a given Hurst function $H(\cdot)$) using
the Choleski decomposition of the covariance matrix, the number of observation points limited to $n=6000$. Although this data length may appear rather small in the present context, some interesting features can nevertheless be noted.
Three cases are considered: Case 1: $H(\cdot)$ is a smooth function; Case 2: $H(\cdot)$ is a trajectory of an integrated FBM with Hurst parameter $0<h<1$ and independent of $X$, therefore $H(\cdot) \in {\cal C}^{\eta-}$ with $\eta=1+h \in (1,2)$; and Case 3: $H(\cdot)$ is a trajectory of a
FBM with Hurst parameter $0<\eta<1$, independent of $X$.\\
For each local Hurst function $H(\cdot)$, Monte-Carlo experiments are realized from $100$ independent replications of observed paths $(X(1/n),X(2/n),\cdots,X(1))$ of the MBM \eqref{MBMintro} for the following choices of parameters:
\begin{itemize}
\item $\alpha=0.2,\,0.3, \,0.4$ and $0.5$;
\item $p=5$ ($=$ the number of dilatations of $\widehat H^{(QV)}_{n,\alpha}$ and $\widehat H^{(QV2)}_{n,\alpha}$, $=$ the regression length providing $\widehat H^{(I\!R2)}_{n,\alpha}$) in all cases;
\item $a=a^*=(1,-2,1)$ and therefore $m=2$ in all cases.
\end{itemize}
Each estimator of  $H(\cdot)$ is computed for $t=\{n^{-\alpha},n^{-\alpha}+0.01,\cdots, \min(1-n^{-\alpha},n^{-\alpha}+0.99) \}$ and therefore an approximation of $\sqrt{M\!I\!S\!E}=\big ( \int_0^1 \E  (\widehat H_n(t)-H(t))^2\d t \big)^{1/2}$ can be computed. Here there are the results of simulations:

\subsubsection*{Case 1: $H(\cdot)\in {\cal C}^\infty$}

\begin{center}
\begin{figure}
\[
\epsfxsize 9.5cm \epsfysize 8cm  \epsfbox{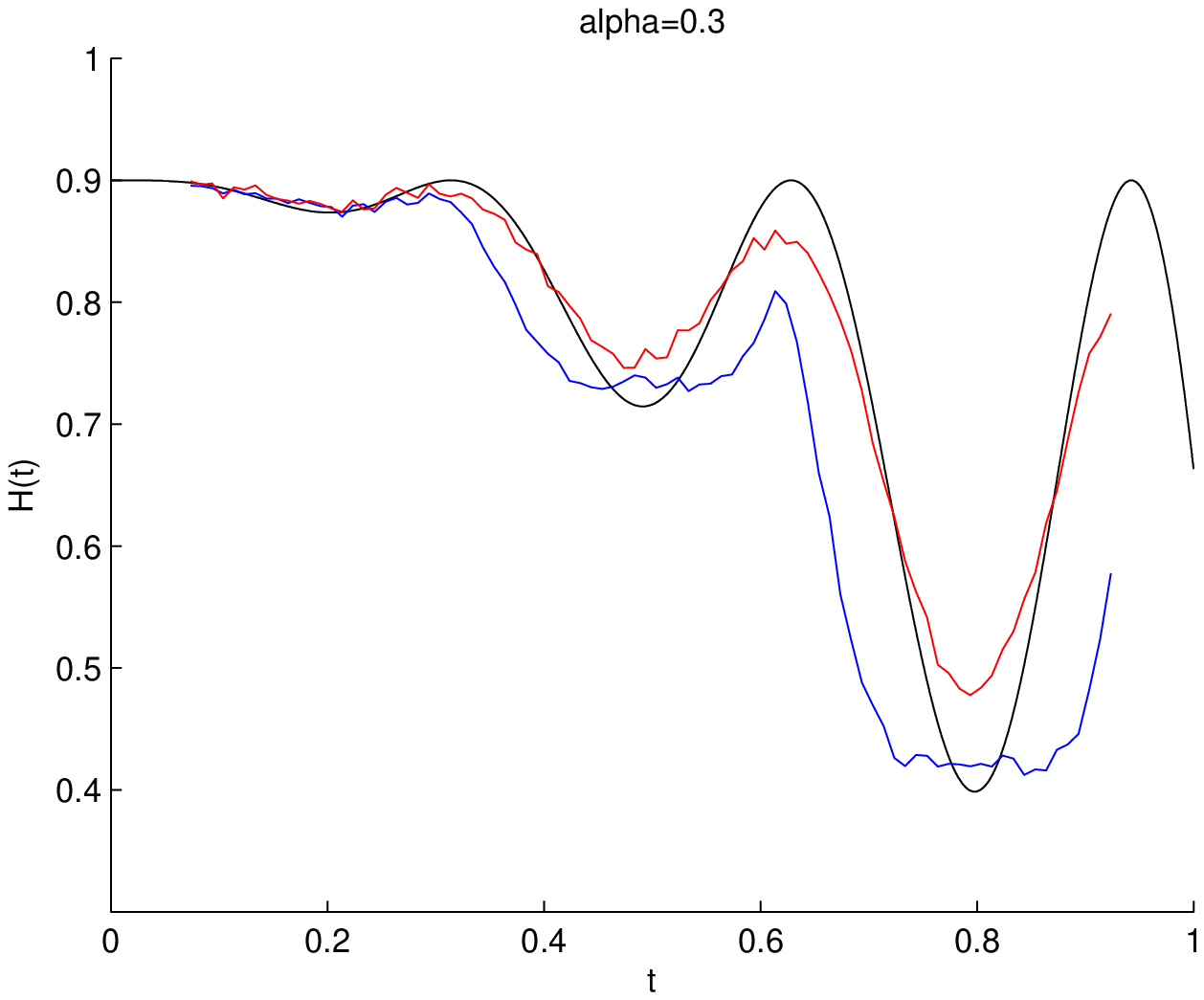} \hspace{-0.3cm}
\epsfxsize 9.5cm \epsfysize 8cm  \epsfbox{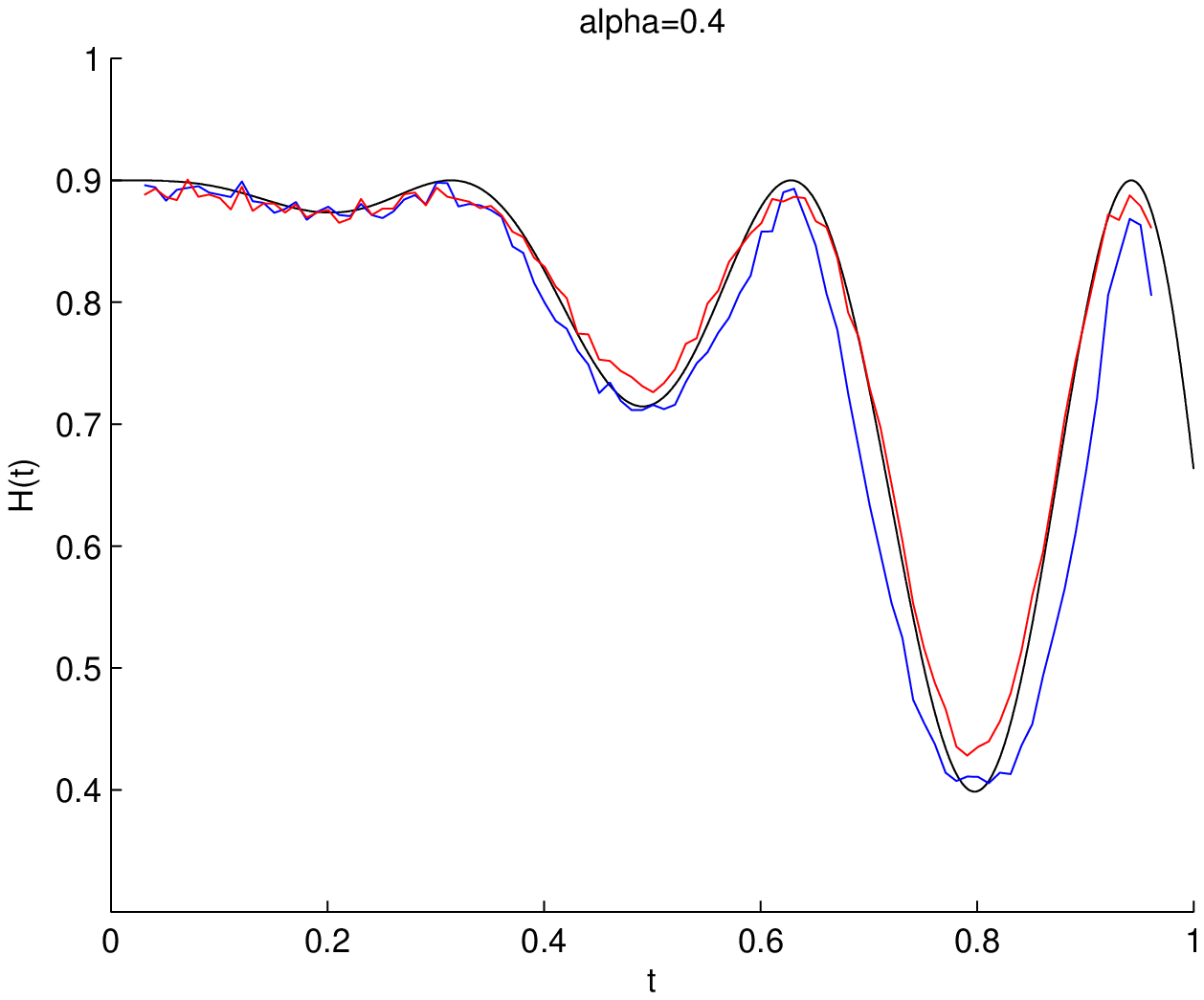}
\]
\[
\epsfxsize 9.5cm \epsfysize 8cm  \epsfbox{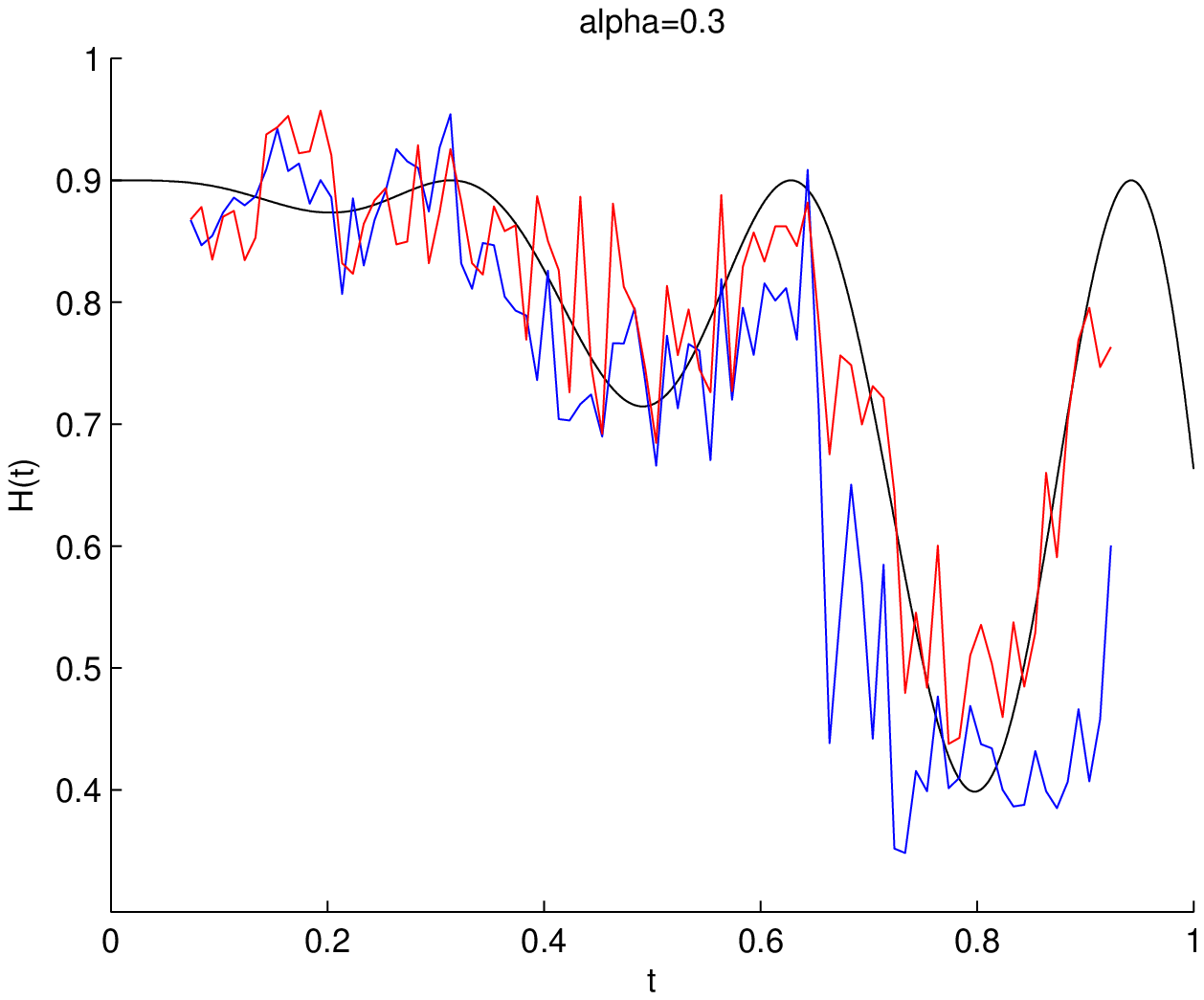} \hspace{-0.3cm}
\epsfxsize 9.5cm \epsfysize 8cm  \epsfbox{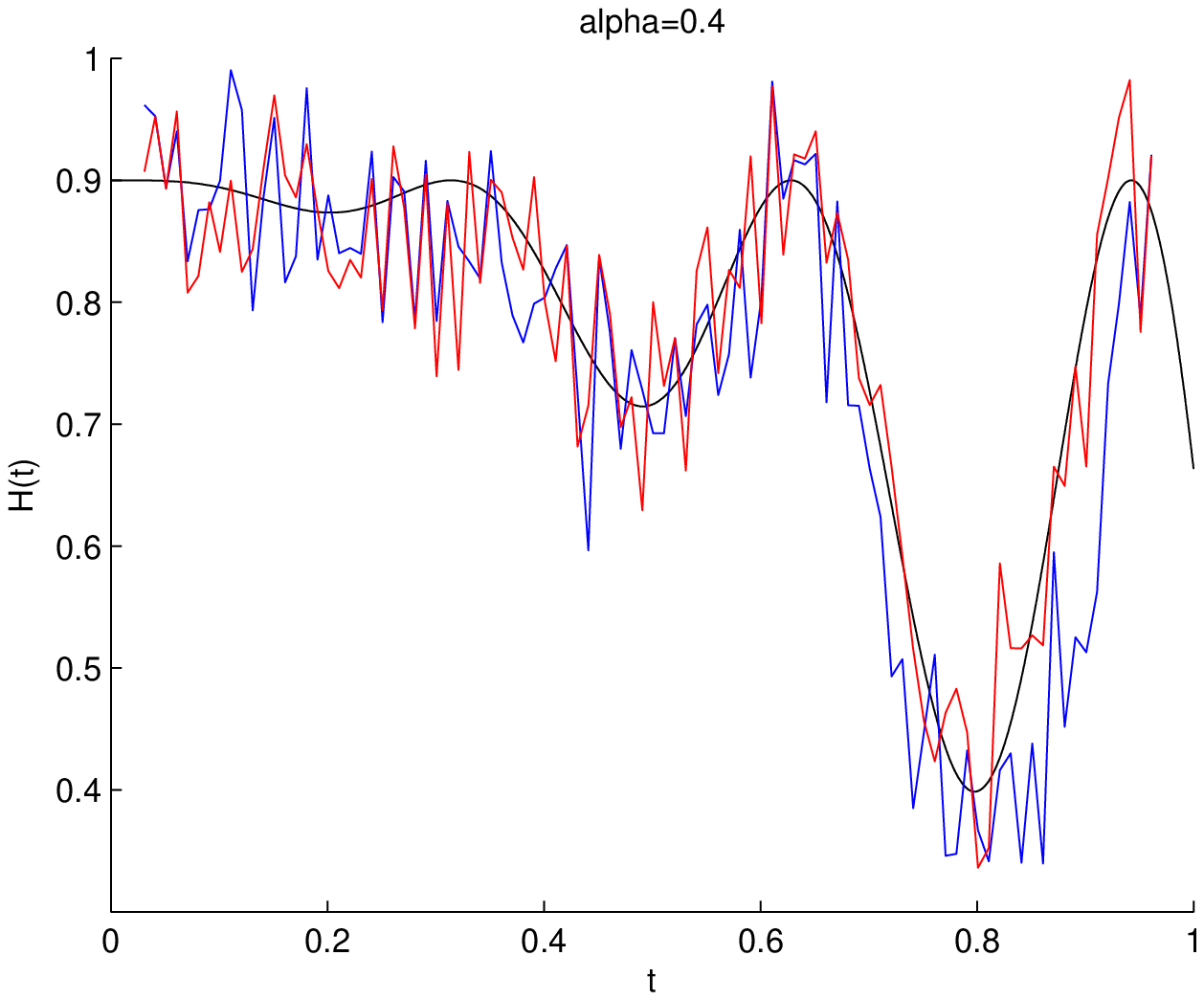}
\]
\caption{\it Estimates of the function $H_4(t)=0.1+0.8(1-t) \sin^2(10 t))$ with $t\in (0,1)$ for $n=6000$ and $\alpha=0.3$ and $0.4$
(from left to right). The top row represents the mean trajectories of $\widehat H_{n,\alpha}^{(QV2)}(t)$ and
$\widehat H_{n,\alpha}^{(IR2)}(t)$ obtained from
$100$ independent replications of MBM with the above function $H(\cdot)$. The bottom row represents a trajectory of $\widehat H_{n,\alpha}^{(QV2)}(t)$ and $\widehat H_{n,\alpha}^{(IR2)}(t)$ obtained from a trajectory of a MBM with the same function $H(\cdot)$. The graphs of $H(t)$,
$\widehat H_{n,\alpha}^{(QV2)}(t),$ and $\widehat H_{n,\alpha}^{(IR2)}(t)$ are in black, blue and red, respectively.}\label{Fig11}
\end{figure}
\end{center}

We have chosen $4$ different functions $t\mapsto H(t)$. 
These functions are:
\begin{itemize}
\item $H_1(t)=0.6$ for any $t\in (0,1)$.
\item $H_2(t)=0.1+0.8t$ for any $t\in (0,1)$.
\item $H_3(t)=0.5+0.4\sin(5 t)$ for any $t\in (0,1)$.
\item $H_4(t)=0.1+0.8(1-t)\sin^2(10 t)$ for any $t\in (0,1)$.
\end{itemize}
From $H_1(\cdot)$ to $H_4(\cdot)$, the functions $H(\cdot)$ display more and more ample fluctuations and less regularity, even if these function are all ${\cal C}^\infty(0,1)$ functions.
Figure \ref{Fig11} provides two kinds of graphs (the mean trajectory and single trajectory) 
of the two  estimators $\widehat H_{n,\alpha}^{(QV2)}$ and $\widehat H_{n,\alpha}^{(I\!R2)}$ of the function $H_4$ for $n=6000$ and two values of $\alpha$ ($\alpha = 0.3$ and $0.4$).
The corresponding MISE are given in Table \ref{table1}.
\begin{table}

\begin{tabular}{|c|c|c|c|c|c|}
\hline
$H_1(t)=0.6$ &$\alpha$ & $0.2$ & $0.3$ & $0.4$ & $0.5$  \\ \hline
$n=2000$ & $\sqrt{\widehat{M\!I\!S\!E}}$ for $\widehat H_{n,\alpha}^{(QV)}$ &$0.044$ &  $0.055$ &   $0.073$ &    $0.104$ \\
& $\sqrt{\widehat{M\!I\!S\!E}}$ for $\widehat H_{n,\alpha}^{(QV2)}$ &$0.041$ &  $0.051$ &   $0.069$ &    $0.096$ \\
& $\sqrt{\widehat{M\!I\!S\!E}}$ for $\widehat H_{n,\alpha}^{(I\!R)}$ & $0.111$ &   $0.137$ &   $0.186$  &  $0.260$  \\
& $\sqrt{\widehat{M\!I\!S\!E}}$ for $\widehat H_{n,\alpha}^{(I\!R2)}$ & $0.061$ &   $0.077$ &   $0.106$  &  $0.145$  \\
\hline
$n=6000$ & $\sqrt{\widehat{M\!I\!S\!E}}$ for $\widehat H_{n,\alpha}^{(QV)}$ &$0.026$ &  $0.035$ &   $0.053$ &    $0.079$  \\
& $\sqrt{\widehat{M\!I\!S\!E}}$ for $\widehat H_{n,\alpha}^{(QV2)}$ &$0.025$ &  $0.033$ &   $0.050$ &    $0.074$  \\
& $\sqrt{\widehat{M\!I\!S\!E}}$ for $\widehat H_{n,\alpha}^{(I\!R)}$ & $0.065$ &   $0.091$ &   $0.128$  &  $0.202$  \\
& $\sqrt{\widehat{M\!I\!S\!E}}$ for $\widehat H_{n,\alpha}^{(I\!R2)}$ & $0.037$ &   $0.049$ &   $0.076$  &  $0.115$  \\
\hline \hline
$H_2(t)=0.1+0.8t$ &$\alpha$ & $0.2$ & $0.3$ & $0.4$ & $0.5$  \\ \hline
$n=2000$ & $\sqrt{\widehat{M\!I\!S\!E}}$ for $\widehat H_{n,\alpha}^{(QV)}$ &$0.170$ &  $0.076$ &   $0.075$ &    $0.101$ \\
& $\sqrt{\widehat{M\!I\!S\!E}}$ for $\widehat H_{n,\alpha}^{(QV2)}$ &$0.170$ &  $0.073$ &   $0.072$ &    $0.096$ \\
& $\sqrt{\widehat{M\!I\!S\!E}}$ for $\widehat H_{n,\alpha}^{(I\!R)}$ & $0.115$ &   $0.143$ &   $0.184$  &  $0.247$  \\
& $\sqrt{\widehat{M\!I\!S\!E}}$ for $\widehat H_{n,\alpha}^{(I\!R2)}$ & $0.059$ &   $0.071$ &   $0.098$  &  $0.135$  \\
\hline
$n=6000$ & $\sqrt{\widehat{M\!I\!S\!E}}$ for $\widehat H_{n,\alpha}^{(QV)}$ &$0.115$ &  $0.045$ &   $0.051$ &    $0.074$  \\
& $\sqrt{\widehat{M\!I\!S\!E}}$ for $\widehat H_{n,\alpha}^{(QV2)}$ &$0.114$ &  $0.044$ &   $0.048$ &    $0.070$  \\
& $\sqrt{\widehat{M\!I\!S\!E}}$ for $\widehat H_{n,\alpha}^{(I\!R)}$ & $0.070$ &   $0.094$ &   $0.134$  &  $0.195$  \\
& $\sqrt{\widehat{M\!I\!S\!E}}$ for $\widehat H_{n,\alpha}^{(I\!R2)}$ & $0.036$ &   $0.046$ &   $0.069$  &  $0.103$  \\
\hline \hline
$H_3(t)=0.5+0.4\sin(5 t)$ &$\alpha$ & $0.2$ & $0.3$ & $0.4$ & $0.5$  \\ \hline
$n=2000$ & $\sqrt{\widehat{M\!I\!S\!E}}$ for $\widehat H_{n,\alpha}^{(QV)}$ &$0.362$ &  $0.125$ &   $0.084$ &    $0.102$ \\
& $\sqrt{\widehat{M\!I\!S\!E}}$ for $\widehat H_{n,\alpha}^{(QV2)}$ &$0.362$ &  $0.123$ &   $0.080$ &    $0.096$ \\
& $\sqrt{\widehat{M\!I\!S\!E}}$ for $\widehat H_{n,\alpha}^{(I\!R)}$ & $0.129$ &   $0.133$ &   $0.171$  &  $0.229$  \\
& $\sqrt{\widehat{M\!I\!S\!E}}$ for $\widehat H_{n,\alpha}^{(I\!R2)}$ & $0.093$ &   $0.071$ &   $0.091$  &  $0.124$  \\
\hline
$n=6000$ & $\sqrt{\widehat{M\!I\!S\!E}}$ for $\widehat H_{n,\alpha}^{(QV)}$ &$0.260$ &  $0.078$ &   $0.056$ &    $0.077$  \\
& $\sqrt{\widehat{M\!I\!S\!E}}$ for $\widehat H_{n,\alpha}^{(QV2)}$ &$0.260$ &  $0.077$ &   $0.052$ &    $0.072$  \\
& $\sqrt{\widehat{M\!I\!S\!E}}$ for $\widehat H_{n,\alpha}^{(I\!R)}$ & $0.078$ &   $0.089$ &   $0.125$  &  $0.180$  \\
& $\sqrt{\widehat{M\!I\!S\!E}}$ for $\widehat H_{n,\alpha}^{(I\!R2)}$ & $0.057$ &   $0.047$ &   $0.065$  &  $0.097$  \\
\hline \hline
$H_4(t)=0.1+0.8(1-t)\sin^2(10 t)$ &$\alpha$ & $0.2$ & $0.3$ & $0.4$ & $0.5$  \\ \hline
$n=2000$ & $\sqrt{\widehat{M\!I\!S\!E}}$ for $\widehat H_{n,\alpha}^{(QV)}$ &$0.321$ &  $0.165$ &   $0.121$ &    $0.120$ \\
& $\sqrt{\widehat{M\!I\!S\!E}}$ for $\widehat H_{n,\alpha}^{(QV2)}$ &$0.320$ &  $0.164$ &   $0.117$ &    $0.112$ \\
& $\sqrt{\widehat{M\!I\!S\!E}}$ for $\widehat H_{n,\alpha}^{(I\!R)}$ & $0.178$ &   $0.138$ &   $0.160$  &  $0.210$  \\
& $\sqrt{\widehat{M\!I\!S\!E}}$ for $\widehat H_{n,\alpha}^{(I\!R2)}$ & $0.165$ &   $0.098$ &   $0.091$  &  $0.112$  \\
\hline
$n=6000$ & $\sqrt{\widehat{M\!I\!S\!E}}$ for $\widehat H_{n,\alpha}^{(QV)}$ &$0.251$ &  $0.136$ &   $0.074$ &    $0.084$  \\
& $\sqrt{\widehat{M\!I\!S\!E}}$ for $\widehat H_{n,\alpha}^{(QV2)}$ &$0.251$ &  $0.135$ &   $0.071$ &    $0.078$  \\
& $\sqrt{\widehat{M\!I\!S\!E}}$ for $\widehat H_{n,\alpha}^{(I\!R)}$ & $0.158$ &   $0.088$ &   $0.115$  &  $0.164$  \\
& $\sqrt{\widehat{M\!I\!S\!E}}$ for $\widehat H_{n,\alpha}^{(I\!R2)}$ & $0.148$ &   {\bf $0.062$} &   $0.067$  &  $0.091$  \\
\hline
\end{tabular} \centering
\label{table1}
\caption{Values of the (empirical) MISE for estimators $\widehat H_{n,\alpha}^{(QV)}$, $\widehat H_{n,\alpha}^{(QV2)}$, $\widehat H_{n,\alpha}^{(I\!R)}$ and $\widehat H_{n,\alpha}^{(I\!R2)}$
of ${\cal C}^\infty$-Hurst functions for $n \in \{2000, \, 6000\}, \alpha \in \{ 0.2, 0.3, 0.4, 0.5\}$ and $a = (1,-2,1)$.
 }

\end{table}

\subsubsection*{Case 2: $H(\cdot)\in {\cal C}^{\eta-}, \, 1<\eta<2$}
\begin{center}
\begin{figure}
\[
\epsfxsize 9.5cm \epsfysize 8cm \epsfbox{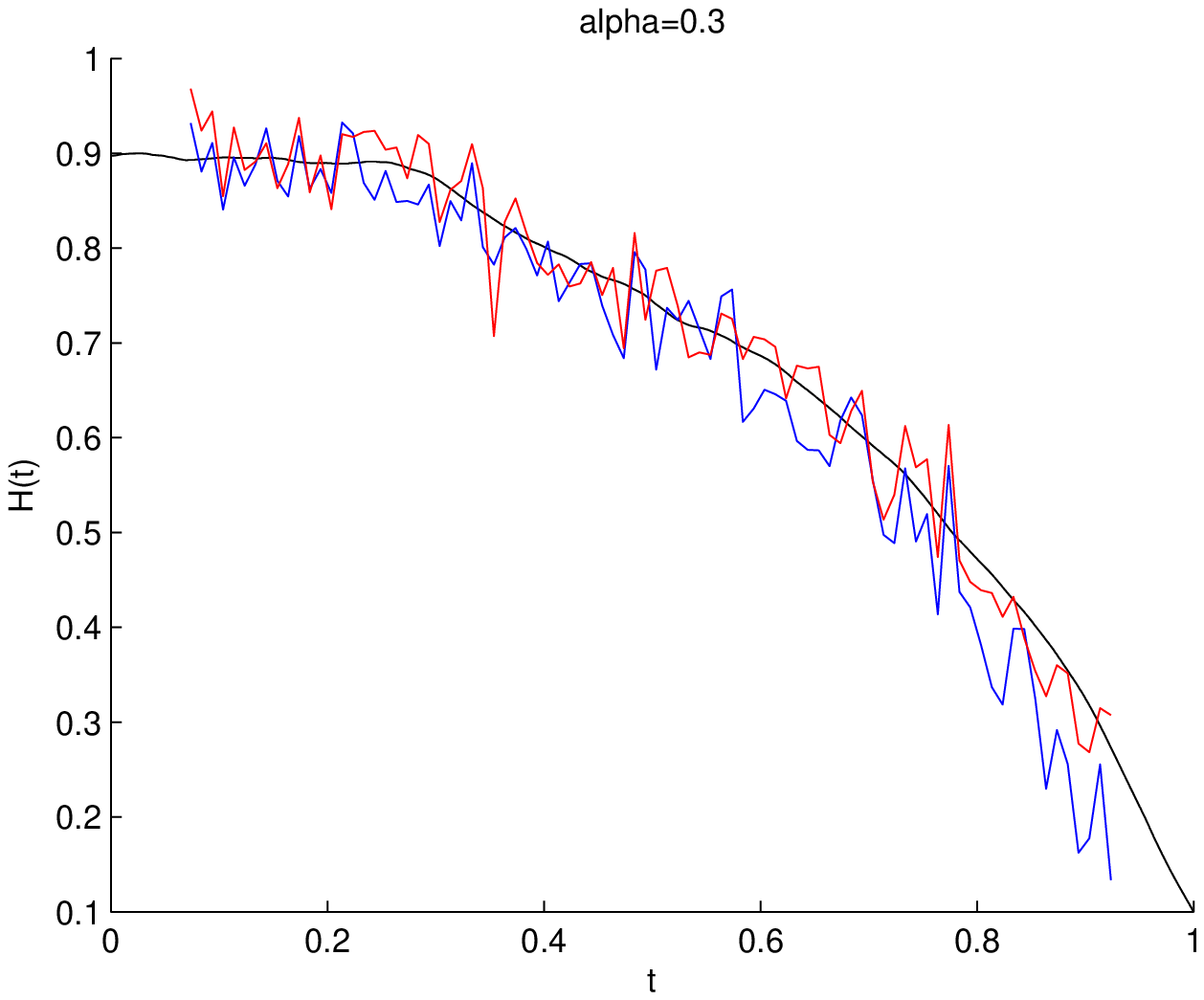} \hspace{-0.3cm}
\epsfxsize 9.5cm \epsfysize 8cm \epsfbox{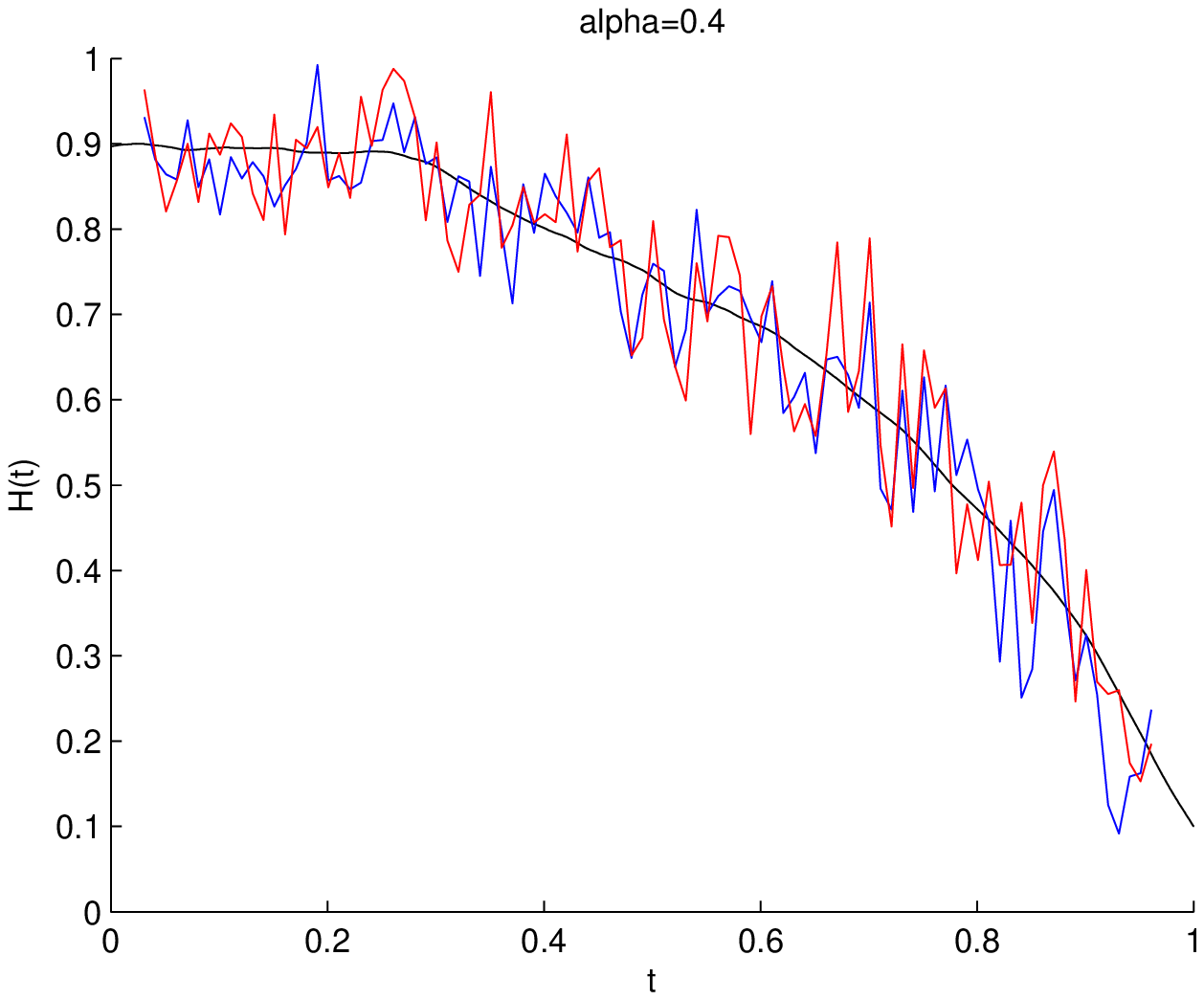} \]
\caption{\it Trajectories of $\widehat H_{n,\alpha}^{(QV2)}(t)$ and $\widehat H_{n,\alpha}^{(I\!R2)}(t)$ for one of the $50$
differentiable Hurst functions $H(\cdot) \in {\cal C}^{1.5-}$ for $n=6000$ and  $\alpha=0.3,\, 0.4$ (from left to right). The graphs of $H(t)$, $\widehat H_{n,\alpha}^{(QV2)}(t)$ and $\widehat H_{n,\alpha}^{(I\!R2)}(t)$ are respectively in black, blue and red.}\label{Fig12}
\end{figure}
\end{center}

We have chosen $\eta=1.5$ but other simulations with different values of $\eta$ lead to similar conclusions. To ensure impartial results, we chose to compute the estimators for $50$ different functions $H(\cdot)$ generated from $50$ trajectories of integrated FBM.
Then $100$ independent replications of processes generated with each function $H(\cdot)$. Table \ref{table2} contains the empirical MISE computed from these $50\times 100$ processes.
Figure \ref{Fig12} provides the trajectories of $\widehat H_{n,\alpha}^{(QV2)}(t)$ and $\widehat H_{n,\alpha}^{(I\!R2)}(t)$ for one of the $50$
differentiable Hurst functions $H(\cdot) \in {\cal C}^{1.5-}$ for $n=6000$ and  $\alpha=0.3,\, 0.4$.
\begin{table}
\begin{tabular}{|c|c|c|c|c|c|}
\hline
$H\in {\cal C}^{1.5}(0,1)$  &$\alpha$ & $0.2$ & $0.3$ & $0.4$ & $0.5$  \\
\hline
$n=2000$ & $\sqrt{\widehat{M\!I\!S\!E}}$ for $\widehat H_{n,\alpha}^{(QV)}$ &$0.261$ &  $0.113$ &   $0.088$ &    $0.103$ \\
& $\sqrt{\widehat{M\!I\!S\!E}}$ for $\widehat H_{n,\alpha}^{(QV2)}$ &$0.261$ &  $0.112$ &   $0.085$ &    $0.098$ \\
& $\sqrt{\widehat{M\!I\!S\!E}}$ for $\widehat H_{n,\alpha}^{(I\!R)}$ & $0.139$ &   $0.141$ &   $0.175$  &  $0.233$  \\
& $\sqrt{\widehat{M\!I\!S\!E}}$ for $\widehat H_{n,\alpha}^{(I\!R2)}$ & $0.098$ &   $0.077$ &   $0.093$  &  $0.128$  \\
\hline
$n=6000$ & $\sqrt{\widehat{M\!I\!S\!E}}$ for $\widehat H_{n,\alpha}^{(QV)}$ &$0.164$ &  $0.067$ &   $0.055$ &    $0.074$  \\
& $\sqrt{\widehat{M\!I\!S\!E}}$ for $\widehat H_{n,\alpha}^{(QV2)}$ &$0.164$ &  $0.066$ &   $0.053$ &    $0.070$ \\
& $\sqrt{\widehat{M\!I\!S\!E}}$ for $\widehat H_{n,\alpha}^{(I\!R)}$ & $0.084$ &   $0.094$ &   $0.131$  &  $0.186$  \\
& $\sqrt{\widehat{M\!I\!S\!E}}$ for $\widehat H_{n,\alpha}^{(I\!R2)}$ & $0.054$ &   $0.047$ &   $0.066$  &  $0.098$  \\
\hline
\hline
$H\in {\cal C}^{0.6}(0,1)$ &$\alpha$ & $0.2$ & $0.3$ & $0.4$ & $0.5$  \\ \hline
$n=2000$ & $\sqrt{\widehat{M\!I\!S\!E}}$ for $\widehat H_{n,\alpha}^{(QV)}$ &$0.140$ &  $0.087$ &   $0.083$ &    $0.096$ \\
& $\sqrt{\widehat{M\!I\!S\!E}}$ for $\widehat H_{n,\alpha}^{(QV2)}$ &$0.140$ &  $0.086$ &   $0.081$ &    $0.094$ \\
& $\sqrt{\widehat{M\!I\!S\!E}}$ for $\widehat H_{n,\alpha}^{(I\!R)}$ & $0.148$ &   $0.156$ &   $0.192$  &  $0.249$  \\
& $\sqrt{\widehat{M\!I\!S\!E}}$ for $\widehat H_{n,\alpha}^{(I\!R2)}$ & $0.088$ &   $0.078$ &   $0.096$  &  $0.135$  \\
\hline
$n=6000$ & $\sqrt{\widehat{M\!I\!S\!E}}$ for $\widehat H_{n,\alpha}^{(QV)}$ &$0.129$ &  $0.067$ &   $0.057$ &    $0.074$  \\
& $\sqrt{\widehat{M\!I\!S\!E}}$ for $\widehat H_{n,\alpha}^{(QV2)}$ &$0.130$ &  $0.067$ &   $0.056$ &    $0.071$  \\
& $\sqrt{\widehat{M\!I\!S\!E}}$ for $\widehat H_{n,\alpha}^{(I\!R)}$ & $0.096$ &   $0.106$ &   $0.143$  &  $0.201$  \\
& $\sqrt{\widehat{M\!I\!S\!E}}$ for $\widehat H_{n,\alpha}^{(I\!R2)}$ & $0.066$ &   $0.052$ &   $0.067$  &  $0.103$  \\
\hline\end{tabular} \centering
\label{table2}
\caption{Values of the (empirical) MISE for estimators $\widehat H_{n,\alpha}^{(QV)}$, $\widehat H_{n,\alpha}^{(QV2)}$, $\widehat H_{n,\alpha}^{(I\!R)}$ and $\widehat H_{n,\alpha}^{(I\!R2)}$
of ${\cal C}^{\eta-}$-Hurst functions for $\eta=1.5$ (up) and $\eta=0.6$ (down), $n \in \{2000, \, 6000\}, \alpha \in \{ 0.2, 0.3, 0.4, 0.5\}$ and $a = (1,-2,1)$.
 }
\end{table}

\subsubsection*{Case 3: $H(\cdot) \in {\cal C}^{\eta-}, \, 0<\eta<1$}
\begin{center}
\begin{figure}
\[
\epsfxsize 9.5cm \epsfysize 8cm \epsfbox{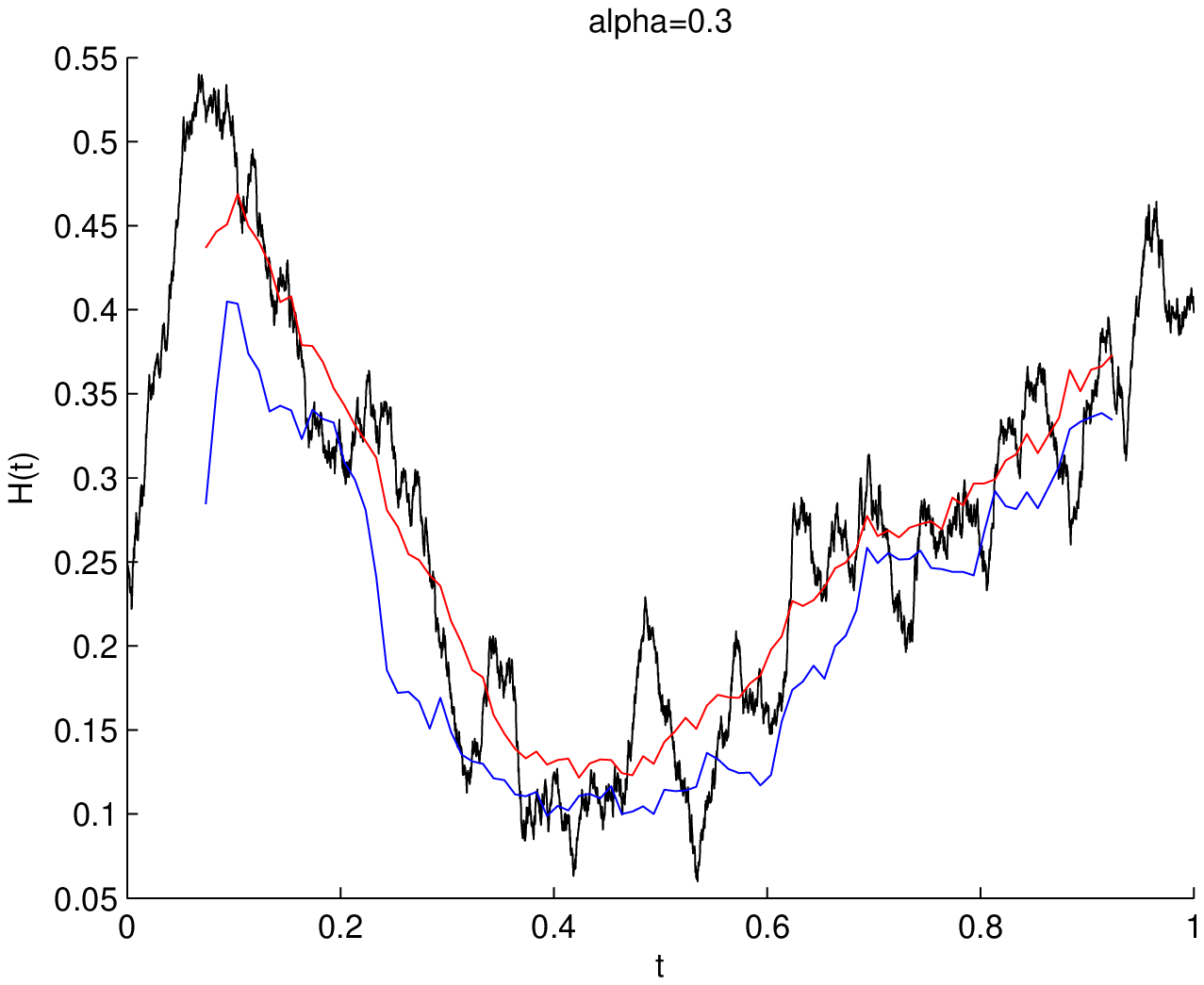} \hspace{-0.3cm}
\epsfxsize 9.5cm \epsfysize 8cm  \epsfbox{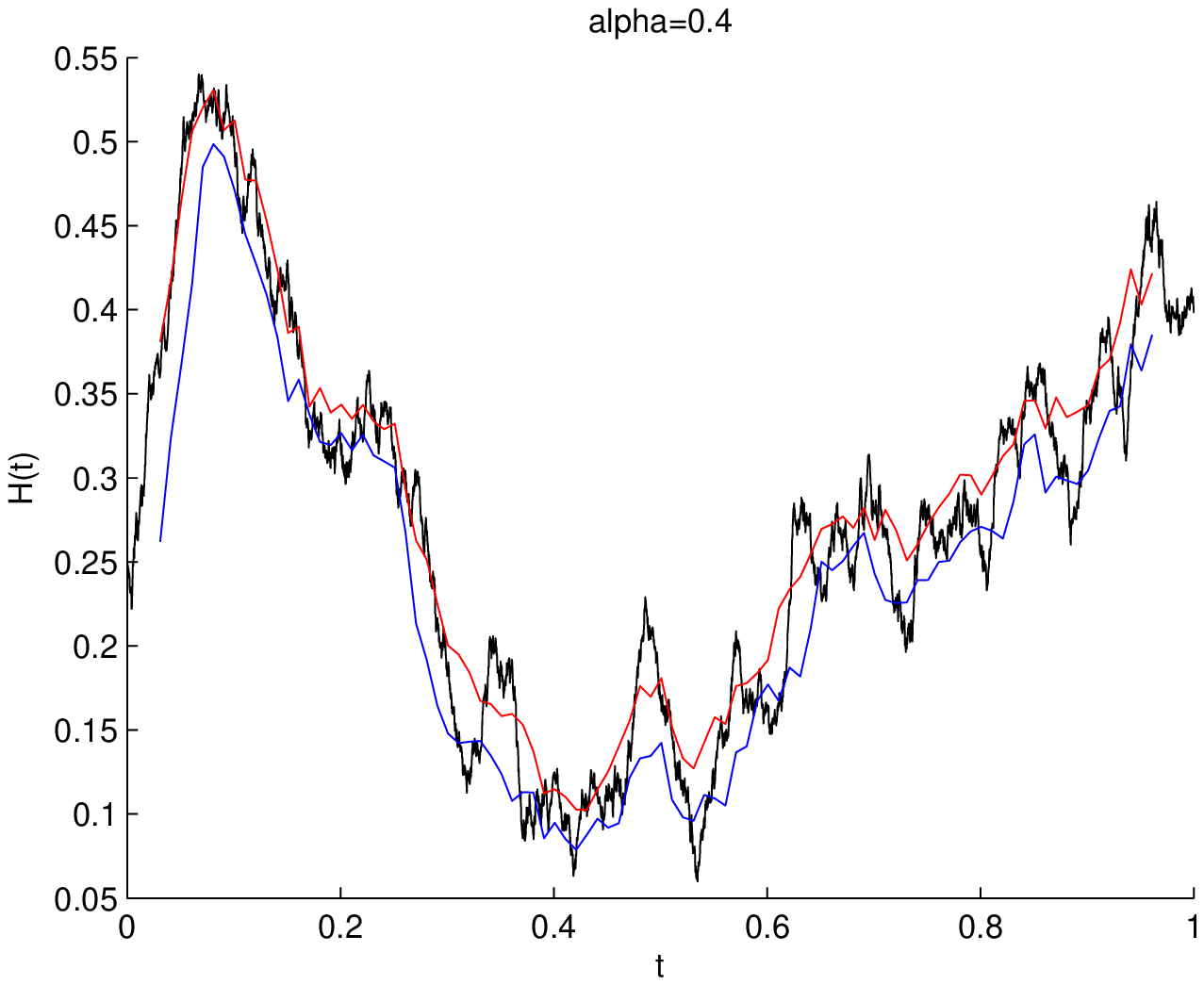}
\]
\caption{\it Mean trajectory (from $100$ independent replications) of $\widehat H_{n,\alpha}^{(QV)}(t)$ and
$\widehat H_{n,\alpha}^{(I\!R)}(t)$ for one of the $50$
differentiable Hurst functions $H(\cdot) \in {\cal C}^{0.6-}$ for $n=6000$ and  $\alpha=0.3,\, 0.4$ (from left to right). The graphs of $H(t)$, $\widehat H_{n,\alpha}^{(QV2)}(t)$ and $\widehat H_{n,\alpha}^{(I\!R2)}(t)$ are respectively in black, blue and red.}\label{Fig13}
\end{figure}
\end{center}

We have chosen here $\eta=0.6$. As in  the case $\eta=1.5^-$, we generated $50$ different trajectories $H(\cdot)$ obtained from trajectories of FBM with parameter $H=0.6$. Then for each function $H(\cdot)$, $100$ independent replications of MBM with Hurst function $H(\cdot)$ are generated. Table \ref{table2} contains the empirical MISE computed from these $50 \times 100$ processes.
An example of a graph of a function $H(\cdot)\in {\cal C}^{0.6-}$ and the mean trajectories (obtained from the $100$ replications) of both the estimators are drawn in Figure \ref{Fig13}.

{\bf Conclusions of simulations:}
\begin{enumerate}
\item In any cases of function $H(\cdot)$, for $n=2000$ or $n=6000$ and $\alpha=0.2,\, 0.3, \, 0.4$ or $0.5$, the empirical $M\!I\!S\!E$ of $\widehat H_{n,\alpha}^{(QV2)}(t)$ and $\widehat H_{n,\alpha}^{(I\!R2)}(t)$ are smaller than the ones of $\widehat H_{n,\alpha}^{(QV)}(t)$ and $\widehat H_{n,\alpha}^{(I\!R)}(t)$, respectively. This theoretically corresponds to the Gauss-Markov Theorem which is also empirically satisfied. Note that PGLSE $\widehat H_{n,\alpha}^{(I\!R2)}(t)$ allows for certain cases to divide by $2$ the empirical $\sqrt{M\!I\!S\!E}$ of $\widehat H_{n,\alpha}^{(I\!R)}(t)$, while the gain obtained with $\widehat H_{n,\alpha}^{(QV2)}(t)$ instead of $\widehat H_{n,\alpha}^{(QV)}(t)$ is more limited (less than $10\%$).
\item In agreement with the theory
(see the remark at the end of Section \ref{Ya}) in the case $H(\cdot)\in {\cal C}^\eta$ with $\eta= 1.5^{-}$ or $\eta=\infty$, our simulation suggests
to select $\alpha$ greater than $(1+2(\eta\wedge 2))^{-1}$ for ``optimal'' results. More precisely, it seems that a good choice is $\alpha\simeq 0.3$ for $\widehat H_{n,\alpha}^{(I\!R)}$ and $\widehat H_{n,\alpha}^{(I\!R2)}$, and $\alpha\simeq 0.4$ for $\widehat H_{n,\alpha}^{(QV)}$ and $\widehat H_{n,\alpha}^{(QV2)}$. This empirical rule seems also valid for $H(\cdot)\in {\cal C}^{0.6^-}$ even if the theoretical choice should be $\alpha\geq  \max\big \{5/11\, ,\, 2(2\sup_{t\in(0,1)} H(t)-0.7)\big \}$. However a bias clearly appears for $\alpha=0.3$ when $H(t)$ is close to $\eta$, while this bias is reduced when $\alpha=0.4$ (see Figure \ref{Fig13}).
The tables also confirm that, with the above choice of the bandwidth $\alpha$,  the accuracy of both estimators increases according
to the theoretical convergence rate $n^{(1-\alpha)/2} ~(\simeq n^{1/3}$)
as $n$ increases from $n=2000$ to $n=6000$ (for instance for $\alpha=0.3$ the empirical $\sqrt{M\!I\!S\!E}$ of $\widehat H_{n,\alpha}^{(I\!R2)}(\cdot)$ seems to be divided by $1.5\simeq 3^{0.35}$).
\item {\em Mutatis mutandis}, the values of the empirical $\sqrt{M\!I\!S\!E}$ are quiet the same for any functions $H(\cdot)\in {\cal C}^\eta$ with $\eta\geq 1$. However, even if this observation seems to be able to be extended for $\eta=0.6$, this does not hold  since in this case the function $H(\cdot)$ is required to vary from $0$ to $0.6$ while for $\eta\geq 1$ it varies from $0$ to $1$.
\item Finally, in all cases, the best estimator is clearly $\widehat H_{n,\alpha}^{(I\!R2)}$ with the choice $\alpha=0.3$ (even when $H$ is a constant, the empirical $\sqrt{M\!I\!S\!E}$ of $\widehat H_{n,\alpha}^{(I\!R2)}$ for $\alpha=0.3$ is smaller than the one of $\widehat H_{n,\alpha}^{(QV2)}$ for $\alpha=0.4$). However the asymptotic variance of $\widehat H_{n,\alpha}^{(I\!R2)}(t)$  is a little larger than the one of $\widehat H_{n,\alpha}^{(QV2)}$. But, as it can be seen in Figure \ref{Fig11} and Figure \ref{Fig13}, $\widehat H_{n,\alpha}^{(I\!R2)}$ is nearly unbiased when $\alpha=0.3$ (and sometimes when $\alpha=0.2$) while the one of $\widehat H_{n,\alpha}^{(QV2)}$ is still biased for $\alpha=0.3$. Then $\widehat H_{n,\alpha}^{(I\!R2)}$  can be used when $\alpha=0.3$ with an asymptotic variance varying with $n^{-(1-\alpha)/2}= n^{-0.35}$ while $\widehat H_{n,\alpha}^{(QV2)}$ requires $\alpha=0.4$ and an asymptotic variance varying with $ n^{-0.3}$.
\end{enumerate}

\section{Proofs of Theorems \ref{Limitgeneral}, \ref{Limitgeneral2} and Corollaries \ref{cor1}, \ref{cor2} }\label{proofs}

The proof of Theorem \ref{Limitgeneral} uses the moment inequality in Lemma \ref{lemgauss}, below, which is a particular case
of a more general moment bound in \cite{bs2009}.

Let $({\mbf Y}_{\!1}, \cdots, {\mbf Y}_{\!N})$ be a collection of
Gaussian vectors ${\mbf Y}_{\!t} = (Y^{(1)}_{t}, \cdots,
Y^{(\nu)}_{t}) \in {\R}^{\nu} $ with zero mean $\E {\mbf Y}_{\!t} = 0 $ and non-degenerated covariance matrices
$\Sigma_t = \big({\rm Cov}\big(Y^{(u)}_{t}, Y^{(v)}_{t}\big)\big)_{1\le u,v\le \nu}, $
having a joint Gaussian distribution
in ${\R}^{\nu N}$. Let $\varepsilon \in [0,1]$ be a fixed number.
Call $({\mbf Y}_{\!1}, \cdots, {\mbf Y}_{\!N})$
{\it $\varepsilon-$correlated} if $\big|{\rm Cor}\big(Y^{(u)}_{t}, Y^{(v)}_{s}\big)\big| \le
\varepsilon $ for any $t \neq s, 1\le t,s\le N$ and any $1\le u, v
\le \nu$. Finally, $\sum'$ denotes the sum over all distinct integers $1 \le t_1, \cdots, t_p \le N, \ t_i \ne t_j \ (i \ne j)$.

Let $|{\mbf x}| = \big(\sum_{i=1}^\nu (x^{(i)})^2 \big)^{1/2} $ denote the Euclidean norm
in $\R^\nu$,
$A = (a_{ij}) $  a real $\nu\times \nu$ matrix, $A^\intercal $ the transposed matrix,
$I$ the unit matrix,
and
$\|A \| := \sup_{|{\mbf x}| =1} |A {\mbf x}|$  the matrix spectral norm, respectively. Lemma \ref{lemgauss}
is a particular case of the moment bound established in Bardet and Surgailis (2011, Corollary 1).

\smallskip

\begin{lem}\label{lemgauss}
Let $({\mbf Y}_{\!1}, \cdots, {\mbf Y}_{\!N})\in {\R}^{\nu N}$ be an $\varepsilon-$correlated
Gaussian vector such that
\begin{equation}\label{cmax}
\max_{1\le t \le N} \|\Sigma_t^{-1} \| \ \le \  c_{\max}
\end{equation}
for some constant $c_{\max} >0.$
Let $G_{j,t}: \R^\nu \to \R,
1\le j \le p \  (p\ge 2), 1\le t \le N $ be measurable functions such that $
\|G_{j,t} \|^2 :=  \E|G_{j,t}({\mbf Y}_{t})|^2 < \infty$.
Assume that for some integer $ 0\le \alpha \le p $,
\begin{equation}\label{Hrank2}
\E \big[ G_{j,t}( {\mbf Y}_{\!t})\big] = 0, \qquad
 \E \big[{\mbf Y}_{\! t}  G_{j,t}( {\mbf Y}_{\!t})\big] = 0, \qquad \text{for all } \quad
1 \le j \le \alpha, \quad 1 \le t \le N
\end{equation}
and that
$\varepsilon \nu^2 c_{\max} < \frac{1}{\nu p-1}. $ Then
\begin{eqnarray}\label{Mombdd}
\sum\nolimits^\prime  \big|\E [G_{1,t_1}({\mbf Y}_{\!t_1}) \cdots
G_{p,t_p}({\mbf Y}_{\!t_p})] \big| &\le& C( \varepsilon, p, \alpha, c_{\max}) K N^{p- \frac{\alpha}{2}} Q_N^{\frac{\alpha}{2} },
\end{eqnarray}
where
\begin{eqnarray}
Q_N&:=&\max_{1\le t \le N} \sum_{1\le s\le N, s\neq t} \max_{1\le
u,v \le \nu} |\E Y^{(u)}_{t} Y^{(v)}_{s}|^2, \label{Q_N}
\end{eqnarray}
$K := \prod_{j=1}^p \max_{1\le t \le N}  \|G_{j,t}\|$, and
the constant $C(\varepsilon, p, \alpha, c_{\max}) $ depends
on $\varepsilon, p, \alpha, c_{\max}$ only.
\end{lem}

\noi In what follows, we use Lemma \ref{lemgauss} with $2 \le p \le 4, \nu =2 $ only. Let
\begin{equation}\label{defSn}
S_n(t):= \frac 1 {v_{n,\alpha}(t)}
\sum_{k\in {\cal V}_{n,\alpha}(t)} \eta_n(k), \qquad \eta_n(k) :=
\frac{\left|V_n^{a^*} Z(k/n) + V_n^{a^*} Z((k+1)/n)\right|}
{|V_n^{a^*} Z(k/n)| + |V_n^{a^*} Z((k+1)/n)| }.
\end{equation}
Thus, $\widehat H^{(I\!R)}_{n,\alpha}(t) = \Lambda_2^{-1} (S_n(t))$.

\medskip
\begin{lem}\label{lemEsp}
Let $Z= (Z(t))_{t \in (0,1)}$ be a zero-mean Gaussian process satisfying
{\bf (A)}$_{{\kappa}}$ with $m=2$. Then for any $0<\varepsilon<1/2$,
\begin{eqnarray} \label{S2rate}
\sup_{\epsilon < t < 1-\epsilon}\big|\E S_n(t)-  \Lambda_2(H(t))\big|&=&o(n^{-\kappa}) + O(n^{-1}+n^{-\alpha (\eta \wedge 2) }).
\end{eqnarray}
\end{lem}


\noindent {\it Proof of Lemma \ref{lemEsp}.} The first relation in (\ref{S2rate}) follows from
\begin{eqnarray} \label{S3rate}
&&\max_{[n\varepsilon/2]\le k \le [(1-\varepsilon/2)n]} \big|\E \eta_n (k) -  \Lambda_2(H(k/n))\big|\ =\  o(n^{-\kappa}), \\
&&\sup_{\epsilon < t < 1-\epsilon}
v^{-1}_{n,\alpha}(t) \Big|
\sum_{k\in {\cal V}_{n,\alpha}(t)} \big\{\Lambda_2\big(H(k/n)\big) - \Lambda_2(H(t))\big\}
\Big|\ \le \  Cn^{-\alpha (\eta \wedge 2) } \label{S33rate}
\end{eqnarray}
To show (\ref{S3rate}), write $\E \eta_n (k)=\E \big [\psi \big (V_n^{a^{*}} Z(k/n)\,, \,V_n^{a^{*}} Z((k+1)/n) \big )\big ]=\Lambda\big (\Cor(V_n^{a^{*}} Z(k/n)\,, \,V_n^{a^{*}} Z((k+1)/n)  \big )$ where $\Lambda(\rho) :=
\frac{1}{\pi}\arccos (-\rho) +
\frac{1}{\pi}
\sqrt{\frac{1+\rho}{1-\rho}} \log \big(\frac{2}{1+\rho}\big),$ see (\cite{bs2010}, p.777), also (\ref{Lambda2}). From Property \ref{popy1},
$$
\max_{[n\varepsilon/2] \leq k \leq [(1-\varepsilon/2)n]} n^\kappa \, \big | \Cor(V_n^{a^{*}} Z(k/n)\,, \,V_n^{a^{*}} Z((k+1)/n)  \big )-\rho_2(H(k/n)) \big | \limiten 0.
$$
Since $\Lambda_2(H) = \Lambda(\rho_2(H))$ and the function $\rho\in (-1,1) \mapsto \Lambda(\rho)$ is analytic we deduce  \eqref{S3rate}. \\
To show (\ref{S33rate}), write
\begin{eqnarray*}
\Lambda_2\big(H(k/n)\big) - \Lambda_2(H(t)) &=& \Lambda'_2(H(t)) \big(H(k/n) - H(t) \big) +
O( \big(H(k/n) - H(t) \big)^2) \\
&=& \Lambda'_2(H(t)) \big(H(k/n) - H(t) \big) + O\big( \big(\frac{k}{n} -t \big)^{2 (\eta \wedge 1)} \big)\\
&= & \Lambda'_2(H(t)) \big(H(k/n) - H(t) \big) + O\big(n^{-2\alpha  (\eta \wedge 1)} \big)
\end{eqnarray*}
 for $k \in {\cal V}_{n,\alpha}(t)$.
Hence
\begin{eqnarray*}
v^{-1}_{n,\alpha}(t) \Big|
\sum_{k\in {\cal V}_{n,\alpha}(t)} \big\{\Lambda_2\big(H(k/n)\big) - \Lambda_2(H(t))\big\}
\Big|
&\le&Cv^{-1}_{n,\alpha}(t)  \Big|
\sum_{k\in {\cal V}_{n,\alpha}(t)} \big\{H(k/n) - H(t)\big\} \Big| + C n^{-2\alpha  (\eta \wedge 1)} \\
&\le&C\big(n^{-\alpha (\eta \wedge 2)} + n^{-1}  + n^{-2\alpha  (\eta \wedge 1)}\big),
\end{eqnarray*}
where for $\eta > 1$ we used $H(k/n) - H(t) =
H'(t) ((k/n) -t) + O(n^{-\alpha (\eta \wedge 2) }) $ and
$v^{-1}_{n,\alpha}(t)\sum_{k\in {\cal V}_{n,\alpha}(t)}\big\{(k/n) - t\big\} = O(n^{-1})$.
This proves (\ref{S33rate}) and Lemma \ref{lemEsp}, too. \hfill $\Box$ \\

\noindent {\it Proof of Theorem \ref{Limitgeneral}.} (i) Recall $\widehat H^{(I\!R)}_{n,\alpha}(t) = \Lambda_2^{-1} (S_n(t))$. Note that the map $H \mapsto \Lambda_2(H)$ is continuous and strictly increasing on $[0,1]$ and the inverse map $s \mapsto \Lambda_2^{-1}(s) $ is continuously differentiable on $s \in [\Lambda_2(0), \Lambda_2(1)]$. Therefore, since
and $S_n(t) \in [0,1]$ a.s., so part (i) follows from $ S_n(t) \limiteproban \Lambda_2(H(t))$ which is a consequence of
\begin{eqnarray} \label{Si}
\E S_n(t) \limiten \Lambda_2(H(t)) \qquad \text{and} \qquad S_n(t) - \E S_n(t)  \limiteproban 0..
\end{eqnarray}
The first term of (\ref{Si}) is obtained from Lemma \ref{lemEsp}. The second member of Lemma \ref{lemEsp} can be deduced from the following bound:
\begin{eqnarray} \label{4ii}
\E \big(S_n(t) - \E S_n(t))^4&<&Cn^{-\rho}, \qquad \text{for some } \rho:= 2(1-\alpha)(1 - \frac{\theta}{\gamma}) >0,
\end{eqnarray}
from Assumption {\bf (B)}$_\alpha$ with $\alpha \in (0,1)$. \\

We shall prove a slightly stronger bound: there exists a constant $C<\infty $ such that for any  integers $0\le s < t < n, \,
|t-s| \le 2n^{1-\alpha}$
\begin{eqnarray} \label{4iii}
\E \big(\sum_{k=s+1}^t \tilde \eta_n(k)\big)^4&<&C(t-s)^2
n^{\frac{2\theta (1-\alpha)}{\gamma}}, \qquad \text{where} \quad
\tilde \eta_n(k) := \eta_n(k) - \E \eta_n (k).
\end{eqnarray}
Relation (\ref{4iii}) implies (\ref{4ii}) with $\rho :=  4(1-\alpha) - 2(1-\alpha) - \frac{2\theta (1-\alpha)}{\gamma}
= 2(1-\alpha)(1 - \frac{\theta}{\gamma})$.\\
In order to apply Lemma \ref{lemgauss} to prove (\ref{4iii}), we rewrite the ratios
$\eta_n(k) = \psi({\mbf Y}_{\!\!n}(k)) $ as functions of Gaussian vectors ${\mbf Y}_{\!\!n}(k) =
(Y^{(1)}_n(k), Y^{(2)}_n(k)) \in \R^2, $ where
\begin{eqnarray}
Y^{(1)}_{n}(k)
&:=&\frac {V_n^{a^{*}} Z(k/n)} {\sigma_{n} (k)}, \qquad Y^{(2)}_{n}(k)\ := \ \frac {V_n^{a^{*}} Z((k+1)/n)} {\sigma_{n} (k+1)}.
\label{Ydef}
\end{eqnarray}
with
\begin{equation}
\label{sigman}
\sigma^2_{n} (k):=\Var \big (V_n^{a^{*}} Z(k/n) \big).
\end{equation}
Also define
\begin{eqnarray}
W^{(1)}_n(k)&:=&\frac{V^{a*} B_{H(k/n)}(0)}{\sqrt{{\rm Var}(V^{a*} B_{H(k/n)}(0))}}, \qquad W^{(2)}_n(k)\ :=\
\frac{V^{a*} B_{H(k/n)}(1)}{\sqrt{{\rm Var}(V^{a*} B_{H(k/n)}(1))}}. \label{Wdef}
\end{eqnarray}
Note $\psi\big (-y^{(1)},-y^{(2)}\big ) = \psi\big (y^{(1)},y^{(2)}\big )$, implying
$\E \big[{\mbf Y}_{\!\!n}(k)  \psi({\mbf Y}_{\!\!n}(k))\big] = 0$; c.f. condition (\ref{Hrank2}) of Lemma
\ref{lemgauss}. Property \ref{popy1} implies that the covariance
matrices $\Sigma_n(k) $ of the ${\mbf Y}_{\!\!n}(k)$'s approach as $n \to \infty$ uniformly in $[n\varepsilon] \le k \le [(1-\varepsilon)n]$ 
the corresponding non-degenerate covariance matrices $\tilde \Sigma_n (k) :=
\left ( \begin{array}{cc} 1 &\rho_2(H(k/n))  \\  \rho_2(H(k/n))&1 \end{array} \right) $
of Gaussian vectors ${\mbf W}_{\!\!n}(k) = \big( W^{(1)}_n(k), W^{(2)}_n(k)\big)$
and therefore these matrices  satisfy the bound in (\ref{cmax}), viz.,
\begin{equation} \label{Sigmabdd}
\max_{1\le k \le n-2} n^\kappa \|\Sigma_n(k) - \tilde \Sigma_n(k) \| \to 0, \qquad \max_{1\le k \le n-2}  \| \Sigma_n(k)^{-1}\| \ \le \ c_{\max}
\end{equation}
with some $c_{\max} <\infty $ independent of $n$. 

Next, we ``decimate'' the sum on the l.h.s. of (\ref{4iii}) so that
the ``remaining'' vectors  ${\mbf Y}_{\!\!n}(k)$'s are $\varepsilon'-$correlated, as follows.
Let $\ell = \ell_n $
be the sequence of integers increasing to $\infty$ at a rate
$o(n^{1-\alpha})$ which is specified in (\ref{ellrate}) below.
Write
\begin{eqnarray}
R_n(s,t)&:=&\sum_{k=s+1}^t \tilde \eta_n(k) \ = \ \sum_{j=0}^{\ell-1} R_{n, \ell}(j),  \label{defRn} \qquad
R_{n, \ell}(j)\ :=\  \sum_{s < k \le t:\, k = j (\!\!\!\!\!\! \mod \ell) } \tilde \eta_n(k).
\end{eqnarray}
Then
\begin{equation}
\E \left(R_n(s,t)\right)^4   \le \ell^4 \max_{0\le j <\ell} \E
\big(R_{n, \ell}(j) \big)^4~~\mbox{and}~~\E \big(R_{n, \ell}(j)\big)^4 \le C \, \Big(\mbox{$\sum_4$} + \mbox{$\sum_3$} +
\mbox{$\sum_{2,1}$}+   \mbox{$\sum_{2,2}$}\Big), \label{Sigma}
\end{equation}
where
\begin{eqnarray*}
\mbox{$\sum_4$}&:=&\sum\nolimits^\prime \big |\E  \big [ \tilde \eta_{n}(k_1) \cdots \tilde \eta_{n}(k_4) \big ] \big |,  \qquad
\mbox{$\sum_3$}\ :=\  \sum\nolimits^{\prime} \big |\E  \big [\tilde  \eta_{n}(k_1)
\tilde \eta_{n}(k_2) \tilde \eta_{n}^2(k_3)\big ]\big |, \label{Sigma3}\\
\mbox{$\sum_{2,1}$}&:=&\sum\nolimits^{\prime} \big |\E  \big [ \tilde \eta_{n}^3(k_1) \tilde \eta_{n}(k_2)\big ] \big |, \hskip1cm
\mbox{$\sum_{2,2}$}\ :=\  \sum \E  \big [\tilde \eta_{n}^2(k_1)
\tilde \eta_{n}^2(k_2)\big ], \nonumber
\end{eqnarray*}
where $ \sum\nolimits^{\prime}$ (respectively, $\sum $) stands for
the sum over all {\it different} integers $ s < k_p \le t, k_p = u \,(\!\!\!\! \mod \ell),
k_p \neq k_q (p\ne q) $ (respectively, over {\it all} $s < k_1, k_2 \le t, k_1 = u \, (\!\!\!\! \mod \ell), k_2 = u \, (\!\!\!\! \mod \ell)
$). But,
 since $|\tilde \eta_{n}(k)| \leq 1$ for all $k, n,$
\begin{equation}
\mbox{$\sum_{2,1}$} + \mbox{$\sum_{2,2}$}
 \le  C \, (|t-s|/\ell)^2. \label{2bdd}
\end{equation}
Let us estimate $\sum_3$ and $ \sum_4$. Let 
\begin{equation}
Q_{n, \ell}(j):= \max_{s <k_1 \le t, k_1 = j (\!\!\!\!\!\! \mod \ell)}
\sum_{s < k_2 \le t, k_2 = j (\!\!\!\!\!\! \mod \ell), k_2 \ne k_1 } \bar \rho^2_{n}(k_1,k_2),  \label{Q_nl}
\end{equation}
where
$$
\bar \rho_{n}(k_1,k_2) := \max_{u,v=1,2} \big |\E \big [Y^{(u)}_{n}(k_1)
Y^{(v)}_{n}(k_2)\big ] \big|
$$
and $Y^{(1)}_n(k), Y^{(2)}_n(k)$ are defined in (\ref{Ydef}). Assumptions {\bf (B)$_{\alpha}$} and {\bf (A)}$_{0}$
imply the following bound. There exists a constant $C_1 < \infty $ such that for any $n \ge 1 $ and any $[n\varepsilon] \le k_1< k_2 \le [(1-\varepsilon)n]$
\begin{equation}\label{barrho}
\bar \rho_{n}(k_1,k_2)\, \leq \, C_1\, n^{\theta (1-\alpha)} (|k_1-k_2|\wedge n^{1-\alpha})^{-\gamma};
\end{equation}
see  also \cite{bs2009}.  Now, given $C_1>0, $
define
\begin{equation}\label{ellrate}
C_2 := C_1/\varepsilon'  \qquad \text{and} \qquad
\ell  :=  [C_2 n^{(\theta/\gamma)(1-\alpha)}].
\end{equation}
Relations (\ref{barrho}) and (\ref{ellrate}) imply that for any $0\le j < \ell$, the Gaussian vector
$\big({\mbf Y}_{\!n}(k), \ k= j \, (\!\!\!\mod \ell), \, k =1, \cdots, n \big)$ is $\varepsilon'-$correlated. We choose
$\varepsilon' < 1/ (28 c_{\max}), $ thus guaranteeing condition $\varepsilon' < 1/(c_{\max} \nu^2 (\nu p-1)$ of Lemma \ref{lemgauss}
for $\nu =2 $ and $p=3,4$.

Using (\ref{barrho}), ({\ref{ellrate}) and $\gamma > 1/2, \, |t-s| \le 2n^{1-\alpha}$  we obtain
\begin{eqnarray}
\nonumber \max_{0\le j <\ell} Q_{n,\ell}(j)&\leq& C \ \big(\frac{n^{\theta (1-\alpha)}}{\ell^\gamma}\big)^2   \max_{1\le k_1 \le |t-s|/\ell}
\sum_{k_2 =1, k_1 \ne k_2}^{|t-s|/\ell} |k_1-k_2|^{-2\gamma}  \\
\label{maxQ} & \leq & 2 C \, \big(\frac{n^{\theta (1-\alpha)}}{\ell^\gamma}\big)^2 \sum_{k'=1}^{\infty} |k'|^{-2\gamma} \ \leq \ C
\end{eqnarray}
Thus, the sums
$\sum_4 $ and $\sum_3 $ can be bounded using Lemma \ref{lemgauss} with $N = [|t-s|/\ell]$.
Accordingly, for $\sum_4 $ with $p=\alpha =4$ we have 
\begin{equation}
\mbox{$\sum_{4}$}  \le C \, (|t-s|/\ell)^{4-\frac{4}{2}}
\big (\max_{0\le j <\ell} Q_{n,\ell}(j)\big )^{\frac{4}{2}}   \le  C \, (|t-s|/\ell)^2.  \label{4bdd}
\end{equation}
Similarly, for $\sum_3 $, with $ p=3$ and $\alpha =2$ we get
\begin{equation}
\mbox{$\sum_{3}$} \le  C \, (|t-s|/\ell)^{3-\frac{2}{2}} \big (\max_{0\le j <\ell} Q_{n,\ell}(j)\big ) ^{\frac{2}{2}}
 \le  C \, (|t-s|/\ell)^{2}. \label{3bdd}
\end{equation}
Relation (\ref{4iii}) follows by combining (\ref{Sigma}), (\ref{2bdd}), (\ref{4bdd}) and (\ref{3bdd}).

\medskip

\noi (ii) Similarly as in the proof (i), it suffices to show
\begin{eqnarray} \label{Sii}
S_n(t) - \E S_n(t)  \limitepsn 0.
\end{eqnarray}
In turn, (\ref{Sii}) is a consequence of (\ref{4ii}) with $\rho>1$ (since $\alpha < \frac{\gamma - 2\theta}{2(\gamma - \theta)}$) and  Chebyshev's inequality imply that for any $0<\delta < (\rho - 1)/4$
\begin{eqnarray*}
\sum_{n=1}^\infty \P \big(|S_n(t) - \E S_n(t)|> n^{-\delta})&<&\sum_{n=1}^\infty n^{4\delta} \E \big(S_n(t) - \E S_n(t))^4 \ \le \
C\sum_{n=1}^\infty n^{4 \delta -\rho}  < \infty.
\end{eqnarray*}
Hence, (\ref{Sii}) follows by the Borel-Cantelli lemma.

\medskip

\noi (iii) We first estimate
$\sup_{\epsilon < t < 1-\epsilon}\big|S_n(t))-  \Lambda_2(H(t))\big|, $ where $S_n(t)$ is defined in (\ref{defSn}). Split
$S_n(t)-  \Lambda_2(H(t)) = (S_n(t)-  \E S_n(t))  + (\E S_n(t)-  \Lambda_2(H(t)))$. Since $v_{n} \equiv v_{n,\alpha}(t) \sim 2n^{1-\alpha}$ do not depend
on $t\in (\epsilon, 1-\epsilon)$ for $n$ large enough, so
\begin{eqnarray}
\sup_{\epsilon < t < 1-\epsilon}\big|S_n(t)-  \E S_n(t)\big|
&\le&v_n^{-1} \max_{0\le s < n-v_n} |R_n(s,s + v_n)| \ \le \  v_n^{-1} \sum_{1 \le i < n/v_n} \max_{(i-1)v_n \le s < iv_n} |R_n(s,s+v_n)|,
\end{eqnarray}
with $R_n(s,t)$ as defined in (\ref{defRn}).
Therefore for any $a >0$
\begin{eqnarray}\label{SRbdd}
\P \big(\sup_{\epsilon < t < 1-\epsilon}\big|S_n(t)-  \E S_n(t)\big| > a\big)
&\le&\sum_{1 \le i < n/v_n} \P\big(\max_{(i-1)v_n \le s < iv_n} |R_n(s,s+v_n)| > av_n \big).
\end{eqnarray}
Since $|R_n(s,s+v_n)| = |R_n((i-1)v_n, s+v_n) - R_n((i-1)v_n,s)| \le |R_n((i-1)v_n, s+v_n)| + |R_n((i-1)v_n,s)|,$ therefore
\begin{eqnarray*}
\P\big(\max_{(i-1)v_n \le s < iv_n} |R_n(s,s+v_n)| > av_n \big)
&\le&2\P\big(\max_{(i-1)v_n \le s < (i+1)v_n} |R_n((i-1)v_n,s)| > av_n/2 \big).
\end{eqnarray*}
We will use the following fact from Billingsley (\cite{Bil}, Theorem 12.1, see also (12.5), (12.10)). {\it Let
$\xi_j, \, 1 \le j \le m$ be arbitrary r.v.'s.
Suppose there exist nonnegative numbers $u_j \ge 0, \, 1 \le j \le m$ such that
\begin{equation} \label{Bill1}
\E \big(\sum_{j = t+1}^s \xi_j\big)^4  \ \le \  \big(\sum_{j=t+1}^s u_j \big)^2, \qquad 0\le t \le s \le m.
\end{equation}
Then for any $\epsilon >0$
\begin{equation} \label{Bill2}
\P\big( \max_{1 \le s \le m}\big|\sum_{j=1}^s \xi_j \big| \ge  \epsilon\big) \ \le \  K \epsilon^{-4} \big(\sum_{j=1}^m u_j \big)^2
+  \P\big( \big|\sum_{j=1}^m \xi_j \big| \ge  \epsilon/2\big),
\end{equation}
where $K>0$ is an absolute constant. }

We shall apply (\ref{Bill2}) with
$$
\xi_j := \tilde \eta_n ((i-1)v_n + j), \qquad u_j := C^{1/2} n^{\frac{\theta (1-\alpha)}{\gamma}}, \qquad
m := v_n \sim 2n^{1-\alpha},
$$
where $C>0$ is the same as in (\ref{4iii}).
The validity of condition (\ref{Bill1}) follows from the bound (\ref{4iii}). Hence from  (\ref{Bill2}), (\ref{4iii}) and Chebyshev's inequality  we obtain
\begin{eqnarray*}
\P\big(\max_{(i-1)v_n < s \le (i+1)v_n} |R_n((i-1)v_n,s)| > av_n/2 \big)
&\le&C(av_n/4)^{-4} v_n^2   n^{\frac{2\theta (1-\alpha)}{\gamma}} +  \P\big(|R_n((i-1)v_n, iv_n)| > av_n/4 \big) \\
&\le&Ca^{-4} v_n^{-2}   n^{\frac{2\theta (1-\alpha)}{\gamma}} +   C(a v_n)^{-4} \E |R_n((i-1)v_n, iv_n)|^4  \\
&\le&Ca^{-4} v_n^{-2}   n^{\frac{2\theta (1-\alpha)}{\gamma}} \\
&\le&
C a^{-4} n^{-2(1-\alpha)( 1- \frac{\theta}{\gamma})}.
\end{eqnarray*}
Hence and from (\ref{SRbdd}) we obtain
\begin{eqnarray}\label{S00rate}
\P \big(\sup_{\epsilon < t < 1-\epsilon}\big|S_n(t)-  \E S_n(t)\big| > a\big)
&\le&Ca^{-4} n^{\alpha -2(1-\alpha)( 1- \frac{\theta}{\gamma})},
\end{eqnarray}
where $C >0$ does not depend on $a$ and $n$. Taking $a = \delta n^{-\mu}$ with $\mu $ given in (\ref{mudef})
results in $\P \big(\sup_{\epsilon < t < 1-\epsilon}\big|S_n(t)-  \E S_n(t)\big| > \delta n^{-\mu} \big)
\le C\delta^{-4}, \, \forall \, \delta >0, $ or
\begin{eqnarray}\label{S1rate}
\sup_{\epsilon < t < 1-\epsilon}\big|S_n(t)-  \E S_n(t)\big|&=&O_p(n^{-\mu}).
\end{eqnarray}
Next, from Lemma \eqref{lemEsp}  we have
\begin{eqnarray} \label{S2rate2}
\sup_{\epsilon < t < 1-\epsilon}\big|\E S_n(t)-  \Lambda_2(H(t))\big|&=&o(n^{-\kappa}) + O(n^{-\alpha (\eta \wedge 2) }) \, = \,
O(n^{-\mu}),
\end{eqnarray}
since $2 \alpha<1$. Finally, (\ref{S2rate2}) and (\ref{S1rate}) imply
\begin{eqnarray}\label{S0rate}
\sup_{\epsilon < t < 1-\epsilon}\big|S_n(t)-  \Lambda_2(H(t))\big|&=&O_p(n^{-\mu}).
\end{eqnarray}
The result in (iii) follows from (\ref{S2rate2}) and the properties of the inverse map $s \mapsto \Lambda_2^{-1}(s) $
mentioned in the proof of (i).

\medskip

\noi (iv) Follows from the Borel-Cantelli lemma,
(\ref{S00rate}) with $a = n^{-(\mu_1 - \delta)} $
 and (\ref{S2rate}) with $\mu$ replaced by $\mu_1$ (which in turn follows from
(\ref{S3rate}), and (\ref{S33rate})). Here, we
use the facts that  $2(1-\alpha)(1- (\theta/\gamma)) - \alpha - 4(\mu_1 - \delta) = 1+ 4\delta  > 1$ and $\alpha (\eta \wedge 2) \ge \mu_1$.
Theorem \ref{Limitgeneral} is proved.  \hfill $\Box$ \\

\medskip


The proof of Theorem \ref{Limitgeneral2}, below,  uses the following lemma which is a particular case
of a more general result in \cite{bs2009} (see also \cite{bs2010}, Theorem A.1).

\smallskip

\begin{lem}\label{tlcgauss2}
Let ${\mbf
Y}_{\!k,N} =  (Y^{(1)}_{k,N}, \cdots, Y^{(\nu)}_{k,N}) \in \R^\nu, \, 1 \le k \le N, \, N \in \N $ be a  triangular
array of jointly Gaussian vectors,   with zero mean $\E {\mbf Y}_{\!k, N}= 0$
and non-generated covariance matrices $\Sigma_{k,N} =  \big({\rm Cov}(Y^{(u)}_{k,N}, Y^{(v)}_{k,N}\big)\big)_{1\le u,v\le \nu}, \,
  1 \le k \le N, \, N \in \N.$
Assume that there exists a function $\rho:\Z\to \R$ such that
\begin{eqnarray}
\forall (j, k)\in \{1,\cdots,N\}^2,~~~~ \max_{1\le u, v \le \nu} \left|
\E Y^{(u)}_{j,N} Y^{(v)}_{k,N} \right|  &\le& |\rho(j-k)|~~~~\mbox{with} ~~\sum_{j\in {\Z}}
|\rho(j)|^2 < \infty. \label{rho_dom2}
\end{eqnarray}
Moreover, assume that for any $\tau \in (0,1)$
any $J \in \N^*$,
\begin{equation}
\left({\mbf Y}_{\!j + [N\tau], N}
\right)_{-J \le j \le J} \
\limiteloin \  \left({\mbf W}_{\!j}\right)_{-J \le j \le J},   \label{Y_W2}
\end{equation}
where ${\mbf W}_{\!j}=  (W^{(1)}_j, \cdots, W^{(\nu)}_j),  \, j\in {\Z}$ is a
stationary Gaussian process taking values in ${\R}^\nu$ and such that the covariance matrix
$\Sigma := \big({\rm Cov}(W^{(u)}_0, W^{(v)}_0\big)\big)_{1\le u,v\le \nu} $ is non-degenerated: $\det \Sigma \ne 0$. In addition,
assume that
\begin{equation} \label{approxS}
\max_{1\le k \le N}\|\Sigma_{k,N} - \Sigma \| \ \to \  0 \qquad (N \to \infty).
\end{equation}
Let $G: \R^\nu \to \R$ be a bounded measurable function such that
\begin{equation} \label{rank22}
\E \left[ {\mbf Y}_{\!k,N} G\big({\mbf
Y}_{\!k,N}\big)\right] \ = \  0, \qquad \forall \, 1\le k \le N, \, N \in \N.
\end{equation}
Then, with $\displaystyle
\sigma^2 :=\  \sum_{j\in {\Z}} {\rm Cov}\big(
G({\mbf W}_0), G({\mbf
W}_j)\big) <\infty, $
\begin{equation}
N^{-1/2} \sum_{k=1}^N \Big \{G\big({\mbf Y}_{\!k,N}\big) - \E G\big({\mbf Y}_{\!k,N}\big)  \Big\}
\limiteloiN {\mathcal N}(0, \sigma^2). \label{CLT2}
\end{equation}
\end{lem}

\medskip

\noi {\it Proof of Theorem \ref{Limitgeneral2}.} For notational simplicity, we restrict the proof
of the $p-$dimensional CLT in (\ref{TLCIR}) to the case $p=2$. Let $0< t_1 < t_2 < 1 $ be arbitrary time moments
and $S_n(t)$ be defined as in  (\ref{defSn}).
First, we prove the following bivariate CLT for $(S_n(t_1), S_n(t_2))$:
\begin{eqnarray} \label{cltS1}
\sqrt{2 n^{1-\alpha}}\big(S_n(t_1)- \Lambda_2(H(t_1)), S_n(t_2)- \Lambda_2(H(t_2))\big)\limiteloin{\mathcal N}(0, \Gamma),
\qquad \Gamma :=  \left( \begin{array}{cc}
\sigma^2(H(t_1)) &0 \\
0 &\sigma^2(H(t_2))
\end{array} \right),
\end{eqnarray}
where $2 n^{1-\alpha} \sim v_n = \# {\cal V}_{n, \alpha}(t_i), i =1,2, $ for all $n $ large enough. Relation
(\ref{cltS1}) follows from $\E S_n(t_i) - \Lambda_2(H(t_i)) = o(n^{(1-\alpha)/2})$ (which is a consequence of (\ref{S2rate}) and the
inequalities for $\kappa $ and $\alpha $ in the statement of Theorem \ref{Limitgeneral2}) and
\begin{eqnarray} \label{cltS2}
&&(2 n^{1-\alpha})^{1/2}(S_n - \E S_n)\ \limiteloin\  {\mathcal N}(0, \sigma^2), \\
&&S_n := c_1 S(t_1) + c_2 S_n(t_2), \qquad \sigma^2 := c^2_1 \sigma^2(H(t_1)) + c^2_2 \sigma^2(H(t_2)), \label{cltS3}
\end{eqnarray}
where $c_1, c_2 $ are arbitrary numbers.

To show (\ref{cltS2}) we use Lemma \ref{tlcgauss2}, with $N := v_n, \, \nu = 4$. Let $ \psi(y^{(1)},y^{(2)}) =
\frac{|y^{(1)} + y^{(2)}|}{|y^{(1)}| + |y^{(2)}|}$  as above and rewrite
\begin{eqnarray*}
S_n&=&N^{-1} \sum_{k=1}^{N} G \left({\mbf Y}_{\! k,N}\right), \qquad G({\mbf y}):= c_1\psi(y^{(1)},y^{(2)}) + c_2\psi(y^{(3)},y^{(4)}),
\qquad {\mbf y} = (y^{(1)},y^{(2)},y^{(3)},y^{(4)}) \in \R^4,\\
{\mbf Y}_{\! k,N}
&:=&\left(Y^{(1)}_{k,N}, \cdots, Y^{(4)}_{k,N}\right), \qquad j_{ni}:= [nt_i - n^{1-\alpha}], \quad  i=1,2,\\
Y^{(1)}_{k,N}&:=&\frac{V^{a^*}_n Z\big(\frac{k+j_{n1}}{n}\big)}{\sigma_n(k+j_{n1})}, \quad
Y^{(2)}_{k,N}\ :=\ \frac{V^{a^*}_n Z\big(\frac{k+j_{n1}+1}{n}\big)}{\sigma_n(k+j_{n1})}, \quad
Y^{(3)}_{k,N}\ := \ \frac{V^{a^*}_n Z\big(\frac{k+j_{n2}}{n}\big)}{\sigma_n(k+j_{n2})}, \quad
Y^{(4)}_{k,N}\  :=\ \frac{V^{a^*}_n Z\big(\frac{k+j_{n2}+1}{n}\big)}{\sigma_n(k+j_{n2})}. \\
\end{eqnarray*}
Similarly to (\ref{Wdef}) define
\begin{eqnarray*}
{\mbf W}_{\! k}
&:=&\left(W^{(1)}_{k}, \cdots, W^{(4)}_{k}\right), \quad k \in \Z, \qquad \sigma^{*2} (H) :=  {\rm Var}(V^{a^*} B_H (0)) = 4 - 4^H,   \\
W^{(1)}_{k}&:=&\frac{V^{a^*} B_{1,H(t_1)}(k)}{\sigma^* (H(t_1))}, \quad
W^{(2)}_{k}\ :=\ \frac{V^{a^*} B_{1,H(t_1)}(k+1)}{\sigma^* (H(t_1))}, \quad
W^{(3)}_{k}\ := \ \frac{V^{a^*} B_{2,H(t_2)}(k)}{\sigma^* (H(t_2))}, \quad
W^{(4)}_{k}\  :=\ \frac{V^{a^*} B_{2,H(t_2)}(k)}{\sigma^* (H(t_2))}, \\
\end{eqnarray*}
where $(B_{i,H(t_i)}(t)), i=1,2 $ are independent FBM's with respective parameters $H(t_i)$. By definition,
$\big({\mbf W}_{\! k}\big)_{k \in \Z}$ is a stationary Gaussian process with ${\mbf W}_{\! 0} $ having a non-degenerated covariance matrix
$\Sigma = \left(\E W^{(u)}_{0} W^{(v)}_{0} \right)_{1\le u,v \le 4}, \, \E (W^{(u)}_{0})^2 = 1, \, u=1, \cdots, 4, \, \E [W^{(1)}_{0} W^{(2)}_{0}] =
\rho_2(H(t_1)), \, \E [W^{(3)}_{0} W^{(4)}_{0}] = \rho_2(H(t_2)), \, \E [W^{(u)}_{0} W^{(v)}_{0}] = 0, \, u \in \{1,2\}, \, v \in \{3,4\}$.

\noindent By definition, $\E (Y^{(i)}_{k,N})^2 = 1, \, i= 1,3. $  From Property \ref{popy1}, it  follows that
$\max_{1\le k \le N}|\E (Y^{(i)}_{k,N})^2 - 1| \to 0, \, i=2,4 $ and
\begin{eqnarray}\label{prop14}
\max_{1\le k,j \le N}
\Big| {\rm Cov}(Y^{(u)}_{k,N}, Y^{(v)}_{j,N})
- {\rm Cov}(W^{(u)}_k,  W^{(u)}_j)\Big|&\to&0
\end{eqnarray}
for $u,v \in \{1,2\}$ and $u,v \in  \{3,4 \}$.  For $u \in \{1, 2\}$ and $v \in \{3,4\}$ relation (\ref{prop14}) follows from
{\bf (B)} and   (\ref{prop12}), implying  for the above $u, v$ that
$\max_{1\le k,j \le N}
\big| {\rm Cov}(Y^{(u)}_{k,N}, Y^{(v)}_{j,N}) \big| \le C (n|t_2-t_1|)^{-\gamma} \to 0. $  Relation
(\ref{prop14}) guarantees the finite-dimensional convergence in (\ref{Y_W2}) for the Gaussian vectors
$\left({\mbf Y}_{\! k,N}\right)$ and
$\left({\mbf W}_{\!k}\right) $ defined above, including the convergence (\ref{approxS}). The covariance bound in
(\ref{rho_dom2}) also follows from {\bf (B)} since $\gamma > 1/2$. Condition (\ref{rank22}) follows from
symmetries ${\mbf Y}_{\! k,N} \egallaw - {\mbf Y}_{\! k,N}$ and $G({\mbf y}) = G(-{\mbf y})$. Finally, the expression
for $\sigma^2 $ in
(\ref{cltS3}) follows by independence of $B_{1,H(t_1)} $ and $B_{2,H(t_2)}$ and the definitions
of $G({\mbf y})$ and $\sigma^2 $ in (\ref{CLT2}).  This proves (\ref{cltS2}) and (\ref{cltS1}) as well.
Then, the clt of Theorem \ref{Limitgeneral2} (\ref{TLCIR}) for $p=2$, together with the expression (\ref{IRvar}) for the limit variance,
follows from (\ref{cltS1}) by an application of the Delta-method. Theorem \ref{Limitgeneral2} is proved. \hfill $\Box$ \\


\noi {\it Proof of Corollary \ref{cor1}.} By definition,
$$
\widehat H^{(I\!R2)}_{n,\alpha}(t) = \sum_{i=1}^p b_i (\widehat H^{(I\!R)}_{n,\alpha}(t)) \widehat H^{(I\!R),i}_{n,\alpha}(t),
\qquad \text{where} \quad b_i(H) := \frac{\1_p^\intercal \,  (\Sigma^{(p)}(H))^{-1} {\mbf e}_i}{ \1_p^\intercal  \,(\Sigma^{(p)}(H))^{-1} \, \1_p},
$$
where
${\mbf e}_i := (\underbrace{0, \cdots, 0}_{i-1}, 1, 0, \cdots, 0)^\intercal \in \R^p$  is a unit vector. Using the fact that
$\sum_{i=1}^p b_i(H) = 1$, we can decompose
\begin{equation}\label{IR2decom}
\widehat H^{(I\!R2)}_{n,\alpha}(t) -H \ =\  W_n(t) + U_n(t),
\end{equation}
where
$$
W_n(t)\ =\  \sum_{i=1}^p b_i(H)\big(\widehat H^{(I\!R),i}_{n,\alpha}(t)- H\big),
\qquad  U_n(t)\ := \ \sum_{i=1}^p \big(b_i (\widehat H^{(I\!R)}_{n,\alpha}(t))- b_i(H)\big)
\big(\widehat H^{(I\!R),i}_{n,\alpha}(t)- H\big)
$$
If view of the assumptions on $\Sigma^{(p)}(H)$, the functions $ H \mapsto b_i(H), \, 1\le i \le p$ are continuous on $(0,1)$ implying
that $|U_n(t)| \le \big(\max_{1\le i \le p}\big|\widehat H^{(I\!R),i}_{n,\alpha}(t)- H\big| \big) \delta( |\widehat H^{(I\!R)}_{n,\alpha}(t))- b_i(H)|)$, where $\delta (x) = o(x), \, x \to 0$.
As a consequence, in all cases (i) - (iv)
the term $U_n(t)$ is negligible w.r.t. the term $W_n(t)$.
From the proof of Theorem 1 it is evident that all statements of Theorem 1 hold  with  $\widehat H^{(I\!R)}_{n,\alpha}$ replaced by
any estimator $\widehat H_{n,\alpha}^{(I\!R),i}, \, 1 \le i \le p$. This proves  the corollary. \hfill $\Box$ \\


\noi {\it Proof of Corollary \ref{cor2}.} From the decomposition (\ref{IR2decom}),
\begin{equation}\label{IR2clt}
\sqrt{2n^{1-\alpha}}\big( \widehat H^{(I\!R2)}_{n,\alpha}(t) -H\big) \ =\  \sqrt{2n^{1-\alpha}}\big(W_n(t) + o_p(1)).
\end{equation}
Clearly, it suffices to prove the corollary with  $\widehat H^{(I\!R2)}_{n,\alpha}(t) -H$ replaced by
$W_n(t)$. In turn, this follows from the $(p \times m)-$dimensional CLT
\begin{eqnarray}
\label{TLCIR22}
\sqrt{2n^{1-\alpha}}\Big( \big (\widehat H^{(I\!R),i}_{n,\alpha}(t_j)-H(t_j)\big )_{1\le i \le p}\Big)_{1\le j \le m}
&\limiteloin&\Big(\big (W^{(I\!R),i}(t_j)\big )_{1\le i \le p}\Big)_{1\le j \le m},
\end{eqnarray}
where $\big (W^{(I\!R),i}(t_j)\big )_{1\le i \le p}, j=1, \cdots, m$ are independent Gaussian vectors with zero mean and respective
covariance matrices $\Sigma^{p)}(H(t_j))$. The proof of (\ref{TLCIR22}) mimics that of Theorem \ref{Limitgeneral2} and
we omit the details. \hfill $\Box$

\section{Proofs of Section \ref{Ya}}\label{proofsMBM}

This section contains the proofs of Propositions \ref{propMBM1} and \ref{propMBM2}
concerning the asymptotic behavior of (co)variances of generalized variations of the general
MBM $X_{(a^+,a^-)}$ for $(a^+,a^-) \in \R^2\setminus \{(0,0)\}$ with local Hurst function $H(\cdot)$ satisfying
Assumption ${\bf (C)}_\eta$. We discuss the general case of filters $a:=(a_0,\cdots,a_q) \in \R^{q+1}$ satisfying condition (\ref{moment}) for some $m \in \N^*$.\\
First introduce the following notations. Using Theorem 4.1 of \cite{StoevTaqqu}, for $t,t' \in (0,1)$
\begin{eqnarray} \label{covY}
\E X_{(a^+,a^-)}(t) X_{(a^+,a^-)}(t')&=&Q( H(t),H(t'),t,t'),
\end{eqnarray}
where $Q$ is a function of $4$ variables $ t, t' \in (0,1), \, H, H' \in (0,1) $ satisfying  $Q( H,H',t,t') = Q( H',H,t',t) $ and defined
for $t \le t' $  by
\begin{eqnarray}
Q( H,H',t,t') :=\left \{ \begin{array}{lll}
\displaystyle   L(H,H')\Big (\frac 1 \pi \,\sin \big(\Delta \beta(H,H')\big)\big (  t' \log t' - t \log t  - (t'-t) \log (t'-t)\big ) & &\\
 & \hspace{-4cm}  +  \cos\big(\Delta \beta(H,H')\big) \, t  \Big )& \mbox{if $H+H'=1$} \\
 \displaystyle \frac 1 2 \, \Big(L_{11}(H,H')\, t^{H+H'} +L_{22}(H,H')\,t'^{H+H'} - L_{12}(H,H')\, (t'-t)^{H+H'}\Big) & & \mbox{if $H+H'\neq 1$},
 \end{array} \right . \label{QQ}
\end{eqnarray}
with
\begin{eqnarray}
L_{11}(H,H')&:=&\frac{L(H,H') \cos \big(\Delta \beta(H,H') - \frac{\pi}{2}(H+ H')\big)}{\cos \big(\frac{\pi}{2}(H+ H')\big)}, \label{L1} \\
L_{22}(H,H')&:=&\frac{L(H,H') \cos \big(\Delta \beta(H,H') + \frac{\pi}{2}(H+ H')\big)}{\cos \big(\frac{\pi}{2}(H+ H')\big)}  \ = \
L_{12}(H,H'),\label{L2}
\end{eqnarray}
 and
where
\begin{equation}
L(H,H')\ :=  \ \frac {\big (K(H)\, K(H')\big)^{1/2}}{K(\frac{H + H'}{2})  },
\end{equation}
$K(\cdot)$ is defined at (\ref{MBMintro}), $\Delta \beta (H,H') := \beta(H') - \beta(H),$
\begin{equation}
\beta (H) \ :=\  \left \{ \begin{array}{ll} \mbox{Arg} \big (a^+ \e^{-\i(H+\frac 1 2 )\frac \pi 2 } +a^- \e^{\i(H+\frac 1 2 )\frac \pi 2 }\big ), &a^+ \ne a^-, \\
0, &a^+ = a^- \ne 0,
\end{array}\right.
\end{equation}
where $\mbox{Arg}(z) \in [0, 2\pi)$ is the principal argument of the complex number $z \ne 0$ (note that
$L(H,H) = 1$ and $L_{ij}(H,H',t,t') = L(H,H') $
for $a^+ = a^- \ne 0$).
In particular,  for any $a \ne 0$
\begin{equation} \label{covX}
\E X_{(a,a)}(t) X_{(a,a)}(t') \  = \ \E X(t) X(t') \ = \    \frac{1}{2}L(H,H')  \big(t^{H+H'} + t'^{H+H'} -  (t'-t)^{H+H'}\big), \qquad
t \le t'.
\end{equation}
Note that the functions
$(H,H') \mapsto  L(H,H'), \, (H,H') \mapsto \Delta \beta(H,H') $ are analytic on $(0,1)^2 $
and hence
$(H,H') \mapsto  Q( H,H',t,t')$ is also
analytic on $(0,1)^2$ for $t' \neq t$ with possible exception on $H + H' = 1$ where
the denominator $\cos \big(\frac{\pi}{2}(H+ H')\big)$ in (\ref{L1})-(\ref{L2}) vanishes. To show the analycity of
$(H,H') \mapsto  Q( H,H',t,t')$ on $H + H' = 1$, let $\omega := H+H' -1 \ne 0$, $0<t<t'<1$ and rewrite $ Q( H,H',t,t')$ in (\ref{QQ})  as
\begin{eqnarray}  \label{Qomega}
Q( H,H',t,t')
&=&\frac{1}{2}L(H,H')\cos (\Delta \beta(H,H')) \big(t^{1 + \omega} + t'^{1+ \omega} - (t'-t)^{1+ \omega}\big) \\
&-&\frac{1}{2} L(H,H') \sin (\Delta \beta(H,H'))  \frac{ \omega}{ \tan(\frac{\omega \pi}{2})}
\Big( t\,\big ( \frac {t^\omega -1}{\omega}\big ) - t'\, \big (\frac {t'^\omega -1}\omega\big ) + (t'-t)\, \big (\frac{(t'-t)^\omega -1}\omega\big ) \Big).  \nonumber
\end{eqnarray}
Then, the first term on the r.h.s. of  (\ref{Qomega}) clearly is analytic in $\omega $ at $\omega =0$
and the same conclusion holds for the second term because of the analycity of the functions
$\omega \mapsto  \frac{ \omega}{ \tan(\omega \pi/2)} $ and
$\omega \mapsto \omega^{-1}(t^\omega -1) \, (t>0)$ in a neighborhood  of $\omega = 0$. It is easy to check
that the limit of (\ref{Qomega}) as $\omega \to 0$ coincides with the definition of $Q(H,H',t,t')$
in (\ref{QQ}) for $H+H'=1$.
As a function of four variables,
$Q( H,H',t,t')$ is analytic on $(0,1)^4$ except for  $t = t'$ but $Q$ is continuous on $(0,1)^4$. \\
For $0< t < t ' < 1 $, introduce the $2\times 2$-matrix of the second order partial derivatives:
\begin{eqnarray}  \label{R2}
R^{(2)} (t,t')&:=&\begin{pmatrix}
R^{(2)}_{H,H'}(t,t')  &R^{(2)}_{H,t'}(t,t')\\ \\
R^{(2)}_{t,H'}(t,t')  &R^{(2)}_{t,t'}(t,t')
\end{pmatrix},
\end{eqnarray}
whose elements are the corresponding partial derivatives of the function $Q $  computed at
$H := H(t),\, H' := H(t')$,
$$
 R^{(2)}_{H,H'}(t,t'):=\frac {\partial^2 Q}{\partial H \partial H'}( H(t),H(t'),t,t'),~
 \cdots \ , ~ R^{(2)}_{t,t'}(t,t'):=\frac {\partial^2 Q}{\partial t \partial t'}(H(t),H(t'),t,t').
 $$
In particular, for all $0<t<t'<1$,
\begin{equation} \label{Qtt}
R^{(2)}_{t,t'}(t,t')= \left \{ \begin{array}{ll}
\frac 1 \pi \, L(H(t),H(t')) \, \sin \big(\beta(H(t')) - \beta(H(t))\big)\, \big ( t'-t \big )^{-1} & \mbox{if $H(t)+H(t') = 1$} \\
 -\frac 1 2 \, L_{12}(H(t),H(t'))
\big (H(t)+H(t')\big )\, \big (H(t)+H(t')-1\big )  \, \big (t'-t\big )^{H(t)+H(t')-2}& \mbox{if $H(t)+H(t') \neq 1$}
\end{array}
\right .
\end{equation}

\bigskip

Finally, define also for $0<t,t'<1$,
\begin{eqnarray}\label{Att}
\hspace{-5mm}A(t,t')&\hspace{-3mm}:=&\hspace{-3mm} \frac {\partial^2 Q}{\partial H \partial H'}( H(t),H(t'),t,t') \\
\label{Btt}
\hspace{-5mm} B_m(t,t')&\hspace{-3mm}:=&\hspace{-3mm}\left \{ \begin{array}{ll}\displaystyle  (-1)^{m+1} \, \frac {(2m-2)! \, \big (\sum_{i=0}^q i^{m}a_i \big )^2} {\pi (m!)^2} L\big(H(t),H(t')\big)\, \sin \big(\beta(H(t')) - \beta(H(t))\big) & \mbox{if $H(t)+H(t') =1$} \\
\displaystyle  (-1)^{m+1} \, \frac {\big (\sum_{i=0}^q i^{m}a_i \big )^2} {2(m!)^2} L_{12}\big(H(t),H(t')\big)  \prod_{\ell=0}^{2m-1}
\big (H(t)+H(t')-\ell\big)& \mbox{if $H(t)+H(t') \neq 1$},
\end{array} \right .
\end{eqnarray}
with $\beta=0$ if $a^+=a^-$. Remark that  $A(t,t')$ and $B_m(t,t')$ are bounded on $(0,1)^2$.

\begin{lem}\label{lem0} Under Assumptions of Proposition \ref{propMBM2} and for any $0<\varepsilon<1/2$,  $[n\varepsilon]\le  k <k' \leq [(1-\varepsilon)n]$, $k'-k\geq 2q$,
\begin{multline}
\label{J}
\Cov\big (V_n^a X_{(a^+,a^-)}(\frac k n)\, , \,  V_n^a X_{(a^+,a^-)}(\frac {k'} n)\big )
=\big (V_n^a H(\frac k n), V_n^a \frac k n \big ) \times R^{(2)} \big (\frac k n ,\frac {k'} n \big )
\times \big (V_n^a H(\frac {k'} n), V_n^a \frac {k'} n \big )^\intercal  \\
+\1_{\{m\geq 2\}} \,B_m\big (\frac {k} n, \frac {k'} n\big )  \frac {\big |k'-k|^{H(\frac k n)+H(\frac {k'}n)-2m}}{n^{H(\frac k n)+H(\frac {k'} n)}} + \lambda_n(k,k').
\end{multline}
The remainder term $\lambda_n(k,k') $ in (\ref{J}) satisfies the following bound: for any $\delta>0$ there exist
$n_0 \in \N^*, \, k_0\in \N^*$ such that for any $n > n_0,\, |k -k'| > k_0, \, [n\varepsilon]\le  k <k' \leq [(1-\varepsilon)n]$
\begin{eqnarray}  \label{lexp1}
|\lambda_n(k,k')|
&\leq&\delta \, \Big(\frac {1}{n^{2(\eta \wedge m)}} +  \frac {|k'-k|^{H(k/n) + H(k'/n) -2m}}
{n^{H(k/n) + H(k'/n)}} \Big).
\end{eqnarray}

\end{lem}

\begin{rem} {\rm
In (\ref{J}) $\displaystyle V_n^a \frac k n = \frac 1 n (\sum_{i=1}^q i a_i) $ or $0$ depending on whether $m=1$ or $m \ge 2$ holds.
In particular, for $a = a^*$ the quadratic form on the r.h.s. of (\ref{J}) coincides with the term $A\big (\frac k n ,\frac {k'}n \big )\, V_n^{a^*} H(\frac k n) \, V_n^{a^*}H(\frac {k'} n)$ in the expansion (\ref{lexp}). }
\end{rem}

\noindent {\it Proof of Lemma \ref{lem0}:} For $k,k' \in \{1,2,\cdots,n-q-1\}$ and $i,j \in \{0,1,\cdots,q\}$ denote
\begin{eqnarray} \label{deltaNot}
&&\delta_i := H(\frac k n+\frac i n)-H(\frac k n), \quad  \delta'_j := H(\frac {k'} n+\frac j n)-H(\frac {k'} n), \quad
\|\delta \| := \max_{1\le i \le q} |\delta_i|, \quad \|\delta' \| :=
\max_{1\le j \le q} |\delta'_j|.
\end{eqnarray}
Assumption $H(\cdot) \in {\cal C}^\eta(0,1)$ implies that there exists $C >0 $ not depending on $k$ and $n$
and such that $\max \{|\delta_i|,|\delta'_j|\} \leq C \,n^{-{\eta\wedge 1}}$.
Using the definition in (\ref{covY}) and the notation in (\ref{deltaNot}),
\begin{eqnarray}\label{VarQ1}
\E \big [V_n^a X_{(a^+,a^-)}(\frac k n) V_n^a X_{(a^+,a^-)}(\frac {k'} n)\big ]
= \sum_{i,j=0}^q a_i \, a_j \, Q\Big (H(\frac k n)+\delta_i, H(\frac {k'} n)+\delta_j', \frac k n+\frac i n,\frac {k'} n+\frac j n \Big ).
\end{eqnarray}
Similarly to (\ref{R2}), for any collection ${\mbf p} = (p_1, \cdots, p_4) \in \N^4 $ of integers and $s \ne s'$,
introduce  partial derivatives of order $|{\mbf p}| := p_1 + \cdots + p_4$ of $Q$ in (\ref{QQ}):
\begin{equation}\label{Rp}
R^{(|{\mbf p}|)}_{H^{p_1}, H'^{p_2},s^{p_3}, s'^{p_4}}(s,s') := \frac{\partial^{|{\mbf p}|} Q}
{\partial H^{p_1} \, \partial H'^{p_2}\,  \partial s^{p_3} \, \partial s'^{p_4} }  (H(s), H(s'), s,s').
\end{equation}
Without loss of generality, assume $k'- k > 2q $ in the sequel.
Using the Taylor expansion of order $2m$ of $Q$ in $(H,H',s,s')$, one obtains:
\begin{eqnarray}\label{taylor}
\E \big [V_n^a X_{(a^+,a^-)}(\frac k n) V_n^a X_{(a^+,a^-)}(\frac {k'} n)\big ]
&=&\sum_{0\le |{\mbf p}|\le 2m}J_n({\mbf p}) + \sum_{|{\mbf p}|= 2m+1} \tilde J_n({\mbf p}),
\end{eqnarray}
where
\begin{eqnarray*}
J_n({\mbf p})&:=&\frac{1}{p_1! \cdots p_4!}
R^{(|{\mbf p}|)}_{H^{p_1}, \cdots, s'^{p_4}}(\frac{k}{n},\frac{k'}{n})\sum_{i,j=0}^q a_i a_j \delta_i^{p_1} (\delta_j')^{p_2} \big(\frac{i}{n}\big)^{p_3} \big(\frac{j}{n}\big)^{p_4}, \\
\widetilde J_n({\mbf p})&\le&C\sup_{I(\frac{k}{n}, \frac{k'}{n})}\big| R^{(|{\mbf p}|)}_{H^{p_1}, \cdots, s'^{p_4}}(s,s')\big|
\sum_{i,j=0}^q \big|\delta_i^{p_1} (\delta_j')^{p_2} \big(\frac{i}{n}\big)^{p_3} \big(\frac{j}{n}\big)^{p_4}\big|
\end{eqnarray*}
and where the set
$I(\frac k n, \frac {k'} n) \subset (0,1)^4 $ is defined by
$$
I(\frac k n, \frac {k'} n) := \Big\{( H, H', s,s'):~ \big |H- H(\frac{k} n)\big | \le \|\delta\|,\,
\big|H'- H(\frac{k'} n)\big| \le \|\delta'\|, \, \big|s - \frac k n\big| \le \frac q n, \, \big|s'- \frac {k'} n\big| \le \frac q n, \, s  <  s' \Big\}.
$$
Note terms in (\ref{taylor}) corresponding to $|{\mbf p}|=0$ and $|{\mbf p}|=1$ vanish because of the moment condition (\ref{moment}).
In particular,
$$
\sum_{|{\mbf p}|=1} J_n({\mbf p}) = \sum_{i,j=0}^q a_i \, a_j \, \Big ( \delta_i \,
R^{(1)}_H\big (\frac k n,\frac {k'} n \big ) + \cdots + \frac j n \, R^{(1)}_{s'}\big ( \frac k n,\frac {k'} n \big ) \Big ) = 0.
$$
In a similar way,
\begin{eqnarray}
\sum_{|{\mbf p}|=2} J_n({\mbf p})&=&\frac 1 2 \Big ( 2  \big (\sum_{i=0}^q a_i \delta_i\big )  \big (\sum_{j=0}^q a_j\delta'_j \big ) R^{(2)}_{H, H'}\big (\frac k n,\frac {k'} n \big )+ \cdots + 2\big (\sum_{i=0}^q a_i \frac i n \big )  \big (\sum_{j=0}^q a_j\frac j n \big ) R^{(2)}_{s,s^{'}}\big (\frac k n,\frac {k'} n \big ) \Big )\nonumber \\
\label{mainR}
&=&\big (V_n^a H(\frac k n),  V_n^a \frac k n \big) \times R^{(2)} \big (\frac k n,\frac {k'} n \big )\times
\big (V_n^a H(\frac {k'} n),  V_n^a \frac {k'} n\big )^\intercal,
\end{eqnarray}
since $\sum_{i=0}^q a_i \delta_i=V_n^a H(\frac k n),\cdots,  \sum_{j=0}^q a_j\frac j n=V_n^a \frac k n$. Note that the partial derivatives in  (\ref{Rp})  involving arguments $H $ and $H'$ contributes to a factor $\log \big ( \frac {k'-k} n\big ) $ at each differentiation while
each differentiation in $s$ or $s'$ contributes a factor $|\frac {k'-k} n \big |^{-1}$; in other words,
\begin{eqnarray}\label{supRbdd}
\sup_{I(\frac k n, \frac {k'} n)} \Big |R^{(|{\mbf p}|)}_{H^{p_1},\cdots,
s'^{p_4}  }(s,s') \Big|&\le&C
\Big(1 + \big |\frac {k'-k} n \big |^{H(\frac k n)+H(\frac {k'} n)-p_3-p_4 } \times \big | \log \big ( \frac {k'-k} n\big ) \big |^{p_1+p_2}\Big).
\end{eqnarray}
Also note
\begin{equation}\label{orderV}
|V_n^a H(\frac k n)| \, \le \, C n^{-(\eta\wedge m)},  \qquad  |V_n^a \frac k n| \, \le \, C n^{-m}  \qquad 1 \le k < n-q.
\end{equation}
From (\ref{orderV}) and (\ref{supRbdd}) we obtain
\begin{eqnarray}\label{mainRbdd}
\big|\sum_{|{\mbf p}|=2} J_n({\mbf p})\big|&\le&C\Big(\frac{1}{n^{2(\eta \wedge m)}} + \frac{|k'-k|^{H(k/n)+H(k'/n)-2}}{
n^{H(k/n)+H(k'/n)}} \1_{\{m=1\}} \Big)
\end{eqnarray}
since $V_n^a \frac k n = 0$ for $m> 1$.  We also have from \eqref{QQ} (with $\beta=0$ if $a^+=a^-$),
\begin{eqnarray}
\nonumber J_n(0,0,m,m)\1_{\{m\ge 2\}}
&\hspace{-3mm} =&\hspace{-3mm} (-1)^{m+1} \frac{(2m-2)! \big (\sum_{i=0}^q a_i i^m \big )^2}{2(m!)^2}  L\big(H(\frac k n),H(\frac {k'} n)\big ) \sin \big (\beta(H(\frac {k'} n)) - \beta(H(\frac k n)) \big )
\, \frac {|k-k'|^{1 -2m}} {n} \\
\nonumber && \hspace{9cm} \mbox{if $H(\frac k n)+H(\frac {k'} n)=1$};\\
\nonumber &\hspace{-3mm} =&\hspace{-3mm} (-1)^{m+1} \frac{\big (\sum_{i=0}^q a_i i^m \big )^2}{2(m!)^2}  L_{12}\big(H(\frac k n),H(\frac {k'} n)\big )
\prod_{\ell=0}^{2m-1} \big (H(\frac k n)+H(\frac {k'} n)-\ell \big)\, \frac {|k-k'|^{H(\frac k n)+H(\frac {k'} n) -2m}} {n^{H(\frac k n)+H(\frac {k'} n)}} \\
\label{2m} && \hspace{9cm} \mbox{if $H(\frac k n)+H(\frac {k'} n) \neq 1$}
\end{eqnarray}
which
coincides with the second expansion term on the r.h.s. of (\ref{J}).
It remains to prove that all the other terms in the Taylor expansion (\ref{taylor}) are negligible with
respect to the main terms in (\ref{J}), in other words, to show the bound  (\ref{lexp1}) on the remainder
term $\lambda_n(k,k')$.
Define
$$
\Lambda_n:=
\frac{\log^2 n}{n^{2(\eta \wedge m)}} + \frac{|k'-k|^{H(k/n) + H(k'/n)-2m}}{n^{H(k/n) + H(k'/n)}}.
$$
We shall prove below that there exist $\beta_1> 0, \, \beta_2 >0 $ independent of $n, k, k' $ and such that
for any $0<\varepsilon<1/2$, $n  \in \N^*$ and $[n\varepsilon] \le k< k' \leq  [(1-\varepsilon)n]$,
\begin{eqnarray} \label{Jnbdd}
J_n({\mbf p})
&\le&\delta_n\Lambda_n, \qquad 2< |{\mbf p}|\le 2m, \quad {\mbf p} \ne (0,0,m,m), \\
\widetilde J_n({\mbf p})
&\le&\delta'_n\Lambda_n, \qquad |{\mbf p}| = 2m+1.  \label{tiJnbdd}
\end{eqnarray}
where
\begin{eqnarray*}\label{deltabdd1}
\delta_n(k',k)&:=&(\frac{1}{n})^{\beta_1} + (\frac{|k'-k|}{n})^{\beta_2}, \qquad \delta'_n(k',k) := \delta_n(k',k) + \frac 1 {|k'-k|}.
\end{eqnarray*}
The statement of the lemma including (\ref{lexp1})
follows from (\ref{taylor}), (\ref{mainR}) and (\ref{Jnbdd}), (\ref{tiJnbdd}).\\
Let us prove (\ref{Jnbdd}). Define 
\begin{eqnarray}\label{Udef}
U_n({\mbf p})&:=&\sum_{i=0}^q a_i \delta_i^{p_1} \big(\frac i n \big)^{p_3}, \qquad
U'_n({\mbf p}) \ := \ \sum_{j=0}^q a_j (\delta'_j)^{p_2} \big(\frac j n\big)^{p_4}.
\end{eqnarray}
We claim that for any $\eta >0$ and $m \ge 1$
\begin{eqnarray}\label{Uclaim}
|U_n({\mbf p})|&\le&C\Big(\frac{1}{n^{m\vee p_3}} + \frac{1}{n^{\eta + p_3}}\Big)
\end{eqnarray}
and  a similar bound holds for $U'_n({\mbf p})$. Let us check that (\ref{Uclaim}) and (\ref{supRbdd}) imply (\ref{Jnbdd}).
Indeed, from the above relations and the definition of $J_n({\mbf p})$ we obtain
\begin{eqnarray*}
|J_n({\mbf p})|&\le&
\big|R^{(|{\mbf p}|)}_{H^{p_1}, \cdots, t'^{p_4}}(\frac{k}{n},\frac{k'}{n})\big|\, |U_n({\mbf p})| \,|U'_n({\mbf p})| \\
&\le&C \, (\log n)^{p_1+p_2} \, \big |\frac {k'-k} n \big |^{H(k/n) + H(k'/n)-p_3-p_4}\Big(\frac{1}{n^{m\vee p_3}} + \frac{1}{n^{\eta + p_3}}\Big)
\Big(\frac{1}{n^{m\vee p_4}} + \frac{1}{n^{\eta + p_4}}\Big) \\
& \le & C \, (\log n)^{p_1+p_2} \,  \big ( I_{11} + I_{12} + I_{21} + I_{22}\big ),
\end{eqnarray*}
where
\begin{eqnarray*}
I_{11}&:=&\big |\frac {k'-k} n \big |^{H(k/n)+ H(k'/n)-p_3-p_4}\frac{1}{n^{2m}}, \qquad
I_{12}\ :=\ \big |\frac {k'-k} n \big |^{H(k/n)+ H(k'/n)-p_3-p_4}\frac{1}{n^{(m\vee p_4) +\eta + p_3}}, \\
I_{21}&:=&\big |\frac {k'-k} n \big |^{H(k/n)+ H(k'/n)-p_3-p_4}\frac{1}{n^{(m\vee p_3) +\eta + p_4}}, \qquad
I_{22}\ :=\ \big |\frac {k'-k} n \big |^{H(k/n)+H(k'/n)-p_3-p_4}\frac{1}{n^{2\eta + p_3+p_4}}.
\end{eqnarray*}
It suffices to prove (\ref{Jnbdd}) for $I_{ij}, \, i,j=1,2$.

Let $H := (H(k/n) + H(k'/n))/2, \, j = |k' -k|, \,  p:= p_3+p_4 $.  For $I_{11}$, (\ref{Jnbdd}) follows from
\begin{equation}\label{delta11}
\big(\frac {j} n \big)^{2H -p}\frac{1}{n^{2m}} = \big(\frac {j} n \big)^{2H} \big(\frac {n} j \big)^{p}\frac{1}{n^{2m}}
\ \le \ \big(\frac {j} n \big)^{2H} \big(\frac {n} j \big)^{2m-1}\frac{1}{n^{2m}}
\ = \
(\frac{j^{2H-2m}}{n^{2 H}})( \frac{j}{n}) \ \le \  \Lambda_n \delta_{11n}, \quad \delta_{11,n} := \frac{j}{n},
\end{equation}
which is immediate from $0< p \le  2m-1 $ and $ 1\le j < n$.

Next, consider $I_{22}$. Assume first that $\eta < m$.
Observe that in this case,
$\Lambda_n  \ge n^{-2\eta}$ if $j\ge n^{(\eta - H)/(m- H)}$ and
$\Lambda_n \ge j^{2H-2m}n^{-2H}$ if $j\le n^{(\eta - H)/(m- H)}$.
Thus in the case $j \ge n^{(\eta - H)/(m- H)}$, the bound (\ref{Jnbdd}) for $I_{22}$ translates
to
\begin{equation}\label{delta221}
(\frac{j}{n})^{2 H - p} \frac{1}{n^{2\eta + p}} \le  \frac{1}{n^{2\eta}} \delta_{22,n}, \quad \text{with}  \quad \delta_{22,n} :=
\frac{j^{2 H - p}}{n^{2H}} = \frac{1}{j^p} (\frac{j}{n} )^{2H} \le (\frac{j}{n} )^{2 H}
\end{equation}
since $j \ge 1 $.
Next, let  $1\le j \le n^{(\eta - H)/(m- H)}, p \ge 1$, then the bound (\ref{Jnbdd}) for $I_{22}$ follows by
\begin{eqnarray}\label{delta222}
(\frac{j}{n})^{2H - p} \frac{1}{n^{2\eta + p}} &\le& \delta_{22,n} \frac{j^{2H - 2m}}{n^{2H}}, \quad \text{with}  \quad
\delta_{22,n} := \frac{j^{2m-p}}{n^{2\eta}}\
= \
\big(\frac{j}{n^{\frac{\eta - H}{m - H}}}\big)^{2m-p} n^{-\beta_{22}} \le n^{-\beta_{22}},
\end{eqnarray}
where $\beta_{22} := \frac{2H(m - \eta) + p(\eta - H)}{m - H} >0.$
Next, consider $I_{22}$ for $\eta \ge m$. Then the first term in the definition of
$\Lambda_n$ is negligible with respect to the second term and the corresponding relation reduces to
\begin{equation}\label{delta223}
\frac{j^{2H-p}}{n^{2\eta + 2H}} \ \le  \ \delta_{22,n} \frac{j^{2H-2m}}{n^{2H}}, \quad \text{with}  \quad
\delta_{22,n} := \frac{j^{2m-p}}{n^{2\eta}}\ \le \ (\frac{j}{n})^{2m}.
\end{equation}

Next, consider $I_{12}$. Accordingly, we need to show
\begin{equation} \label{delta12}
\big(\frac {j} n \big)^{2H-p_3-p_4}\frac{1}{n^{(m\vee p_4) +\eta + p_3}}\ \le \  \delta_n \big(\frac{1}{n^{2(\eta \wedge m)}}+
\frac{j^{2H - 2m}}{n^{2H}}\big),
\end{equation}
where $\delta_n $ satisfies the bound in (\ref{deltabdd1}). Let $\eta < m $ then (\ref{delta12}) becomes
\begin{equation} \label{delta12m}
\big(\frac {j} n \big)^{2H-p_3-p_4}\frac{1}{n^{(m\vee p_4) +\eta + p_3}}\ \le \  \delta_n \big(\frac{1}{n^{2\eta}}+
\frac{j^{2H - 2m}}{n^{2H}}\big).
\end{equation}
Let  $j > n^{\frac{\eta -H}{m-H}}$. The first term on the r.h.s. of (\ref{delta12m})
dominates the second one and (\ref{delta12}) for $2H \le p_3 + p_4 $ follows from
\begin{eqnarray*} \label{delta12m1}
\big(\frac {j} n \big)^{2H-p_3-p_4}\frac{1}{n^{(m\vee p_4) +\eta + p_3}}\ \le \  \delta'_{12,n} \frac{1}{n^{2\eta}}, \quad \text{with}  \quad
\delta'_{12,n} := n^{-\beta'_{12}}, \qquad \beta'_{12} >0,
\end{eqnarray*}
where 
\begin{eqnarray*}
\beta'_{12}&:=&2H - \eta + (m\vee p_4) - p_4 + \frac{\eta -H}{m-H} (p_3+p_4 - 2H)\\
&=&(m-H)^{-1} \times \begin{cases} (p_3+p_4 -m)(\eta - H) + H(m - \eta), & p_4 \ge m, \\
(m- p_4 + H)(m - \eta) + p_3(\eta - H), & p_4 < m.
\end{cases}
\end{eqnarray*}
Thus $\beta'_{12} \ge  H(m- \eta)/(m-H) >0$.  On the other hand, if $p_3 + p_4 < 2H < 2$ then the first inequality
in (\ref{delta12m}) holds with $\delta_n=\delta'_{12,n}$ where
\begin{eqnarray*} \label{delta12m2}
\delta'_{12,n} :=  \big(\frac {j} n \big)^{2H-p_3-p_4}  \frac{1}{n^{(m\vee p_4) - \eta + p_3}}\ \le \  n^{-(m - \eta)},
\end{eqnarray*}
with $m - \eta >0$, since $m\ge 2, \, p_4 < 2 $ and  $(m\vee p_4) - \eta + p_3 \ge m- \eta > 0$.

\noindent Next, let $j \le n^{\frac{\eta -H}{m-H}}$. Then the second term on the r.h.s. of (\ref{delta12m}) dominates the first one
and (\ref{delta12}) follows from
\begin{eqnarray*} \label{delta12m3}
\big(\frac {j} n \big)^{2H-p_3-p_4}\frac{1}{n^{(m\vee p_4) +\eta + p_3}}\ \le \  \delta''_{12,n}
\frac{j^{2H - 2m}}{n^{2H}}, \quad \text{with} \quad \delta''_{12,n} := \frac{j^{2m - p_3-p_4}}{n^{\eta + (m \vee p_4) - p_4}} \le n^{-\beta''_{12}},  \end{eqnarray*}
where
\begin{eqnarray*}
\beta''_{12}&:=&\eta + (m \vee p_4) - p_4 -  (\frac{\eta -H}{m-H})(2m - p_3-p_4)\\
&=&(m-H)^{-1} \times \begin{cases} (p_3+p_4 -m)(\eta - H) + H(m - \eta), & p_4 \ge m, \\
(m- p_4 + H)(m - \eta) + p_3(\eta - H), & p_4 < m,
\end{cases}
\end{eqnarray*}
and therefore  $\beta''_{12} = \beta'_{12}$, implying $\beta''_{12} >0$ as above.

It remains to consider $I_{12}$ for $\eta \ge m$. In this case the second term on the r.h.s. of (\ref{delta12m}) dominates the first one
for all $1\le j \le n$
and (\ref{delta12}) follows from
\begin{eqnarray*} \label{delta12m4}
\big(\frac {j} n \big)^{2H-p_3-p_4}\frac{1}{n^{(m\vee p_4) +\eta + p_3}}\ \le \  \delta'''_{12,n}
\frac{j^{2H - 2m}}{n^{2H}}, \quad \text{with} \quad \delta'''_{12,n}\ := \ \frac{j^{2m - p_3-p_4}}{n^{\eta + (m \vee p_4) - p_4}}\ \le \
(\frac{j}{n})^{2m - p_3-p_4} \ \le \  (\frac{j}{n})
\end{eqnarray*}
since $p_3 + p_4 < 2m $ and $\eta + (m\vee p_4) - 2m + p_3 = (\eta - m) + ((m\vee p_4) - m) + p_3 \ge 0$.  This proves
(\ref{delta12}) or
(\ref{Jnbdd}) for $J_{12}$. Since  consideration of $I_{21}$ is completely analogous, the proof of (\ref{Jnbdd}) is now complete.

\smallskip

Let us prove the claim (\ref{Uclaim}). Write $\delta_i = \delta_{0i} + \Gamma_i$, where for $i\in\{0,\cdots,q\}$ and $r=\{0,\cdots,[\eta]\}$,
\begin{eqnarray*}
\Gamma_i:=\sum_{r=1}^{[\eta]}\frac {H_r}{r !} \big (\frac i n \big )^r , \qquad
\qquad H_r := H^{(r)}(\frac k n ).
\end{eqnarray*}
By Assumption C($\eta$),
\begin{equation} \label{deltabdd}
\max_i(|\delta_{0i}|) \leq C \, n^{-\eta}, \qquad \max_r (|H_r|) \leq C.
\end{equation}
From binomial expansion,
$$
U_n({\mbf p}) \ = \
\sum_{\tau_0+\cdots+\tau_{[\eta]}=p_1}    {p_1 \choose \tau_{0}, \cdots, \tau_{[\eta]}}
\Big(\prod_{r=1}^{[\eta]} \frac{ H^{\tau_{r}}_r}{r!} \Big)
\sum_{i=0}^q a_i  \delta_{0i}^{\tau_{0}}
\big(\frac{i}{n}\big)^{p_3+ \sum_{r=1}^{[\eta]} r \tau_{r}},
$$
where $\tau_{i} \in \N, \, i= 0,1 \cdots, [\eta]$
are nonnegative integers. 
According to (\ref{moment}), for any $s=0,1, \cdots $
$$
\Big|\sum_{i=0}^q a_i \big(\frac{i}{n}\big)^{s}\Big|\ \le
\ C\,n^{-s} \1(s \ge m).
$$
Therefore from
(\ref{deltabdd}) it follows that
\begin{eqnarray*}
\Big|\sum_{i=0}^q a_i  \delta_{0i}^{\tau_{0}}
\big(\frac{i}{n}\big)^{p_3+ \sum_{r=1}^{[\eta]} r\tau_{r}}\Big|
&\le&C\begin{cases}n^{-\eta \tau_0 - p_3- \sum_{r=1}^{[\eta]}r\tau_{r}}, &\tau_0 \ge 1, \\
n^{-p_3- \sum_{r=1}^{[\eta]}r\tau_{r}}, & \tau_0= 0 \ \text{and} \ p_3+ \sum_{r=1}^{[\eta]}r\tau_{r} \ge m, \\
0, &\text{otherwise}.
\end{cases}
\end{eqnarray*}
Note that for $\tau_0 =0$, we have $\sum_{r=1}^{[\eta]}r\tau_{r}
\geq p_1$. Hence,
$$
n^{-p_3- \sum_{r=1}^{[\eta]}r\tau_{r}}\1\big(\tau= 0, p_3+ \sum_{r=1}^{[\eta]}r\tau_{r} \ge m\big)
\ \le \ (n^{-m}) \wedge (n^{-p_1-p_3}) \  \le \ n^{-(m\vee p_3)}.
$$
Therefore
\begin{eqnarray*}
\Big|\sum_{i=0}^q a_i  \delta_{0i}^{\tau_{0}}
\big(\frac{i}{n}\big)^{p_3+ \sum_{r=1}^{[\eta]} r\tau_{r}}\Big|
\ \le \ C\big(n^{-\eta  - p_3} + n^{-(m \vee p_3)}\big),
\end{eqnarray*}
proving (\ref{Uclaim}) and (\ref{Jnbdd}).

It remains to prove (\ref{tiJnbdd}). Let $\widetilde U_n({\mbf p}):= \sum_{i=0}^q \big|a_i  \delta_i^{p_1} \big(\frac i n \big)^{p_3}\big|.$ Since $\|\delta\| \le C\, n^{-(\eta \wedge 1)}, $ see (\ref{deltabdd}), it  follows that
\begin{equation}\label{tiUbdd}
\widetilde U_n({\mbf p}) \ \le  C n^{-(\eta \wedge 1)p_1 - p_3}
\end{equation}
and a similar bound holds for $\widetilde U^\prime_n({\mbf p}):=
\sum_{j=0}^q \big|a_j  (\delta'_j)^{p_2} \big(\frac j n \big)^{p_4}\big|$. Relations (\ref{tiUbdd}) and
(\ref{supRbdd}) imply
\begin{eqnarray}\label{supJbdd}
\widetilde J_n({\mbf p})&\le&
C\Big(1 + \big |\frac {k'-k} n \big |^{2H-p}\times  \big |\log \big ( \frac {k'-k} n \big )\big |^{2m+1-p} \Big)\, n^{-(\eta \wedge 1)(2m+1 - p) - p}
\end{eqnarray}
where $p= p_3+ p_4 \in \{0,1, \cdots, 2m+1\}$ and $|{\mbf p}| = 2m+1$. Relation (\ref{tiJnbdd}) now can be  verified
directly, by inspecting the cases $\eta < 1, \eta \ge 1, p=0,1, 2m+1 $ and $ 1< p < 2m+1$ separately. This ends the proof
of Lemma \ref{lem0}.
\hfill $\Box$  \\
~\\

\noindent {\it Proof of Proposition \ref{propMBM1}.} Note first that the result of Lemma \ref{lem0} cannot be used for establishing  Proposition \ref{propMBM1} because Proposition \ref{propMBM1} treats the case $k$ close to $k'$ which is excluded by the proof of Lemma \ref{lem0} (the second derivative $R_{t,t'}$ does not exists for $t=t'$). \\
For any  $[n\varepsilon] \leq k\leq [(1-\varepsilon)n]$ similarly to (\ref{VarQ1}) and using the notation in (\ref{deltaNot}) with $k' = k+\ell $
we rewrite the variance ${\cal V}_n(k) :=\Var \big ( V_n^a X_{(a^+,a^-)}(\frac k n) \big)$ as
\begin{eqnarray}
\label{VarQ}
{\cal V}_n(k)&=& \frac 1 2 \, \sum_{i,j=0}^q a_i \, a_j \, Q\Big (H(\frac k n)+\delta_i,H(\frac {k} n)+\delta_j,
\frac k n+\frac i n,\frac k n+\frac {j} n\Big ).
\end{eqnarray}
Next, we use a Taylor expansion of $Q$ given by two analytic expressions in \eqref{QQ}, see also
(\ref{Qomega}). In order to avoid cumbersome notation, we shall assume in the rest of the proof that the function
$H(\cdot)$ is separated from $1/2$, i.e., $\sup_{t \in (0,1)} |H(t) - 1/2| >0$. Then the second expression
for $Q$ in  \eqref{QQ} applies, viz., $Q(H,H', t,t' ) = (\frac 1 2)(Q_{11}(H,H',t) + Q_{22}(H,H', t') - Q_{12}(H,H',t,t')) $, where
$$
Q_{11}(H,H',t) := L_{11}(H,H') t^{H+H'}, \quad Q_{22}(H,H',t) := L_{22}(H,H') t'^{H+H'}, \quad
Q_{12}(H,H',t,t') := L_{12}(H,H') |t'-t|^{H+H'}.
$$
Accordingly, ${\cal V}_n(k) = \frac 1 2 \big({\cal V}_{n11}(k) +  {\cal V}_{n22}(k)-  {\cal V}_{n12}(k)\big)$, where
$ {\cal V}_{nij}$ are defined as in (\ref{VarQ}) with $Q$ replaced
by $Q_{ij}, \, 1\le i\le j\le 2$. Note that $(H,H',t)\in (0,1)^3 \mapsto Q_{11}(H,H',t)$ is analytic except for $H+H' =  1$ and the last
case is excluded by our assumption that $H(\cdot)$ is separated from $1/2$.
Using the last fact, the Taylor expansion of $Q_{11}$ as in (\ref{taylor}) in the proof of Lemma \ref{lem0} and the bounds (\ref{orderV}),
it follows that
for all $n$ large enough,
\begin{eqnarray}
\label{cte4}\sup_{[n\varepsilon] \leq k\leq [(1-\varepsilon)n]} \ \big|{\cal V}_{n11}(k)\big|
\leq C(\varepsilon) \, \frac 1 {n^{2(\eta\wedge m)}}
\end{eqnarray}
with $C(\varepsilon)\geq 0$ not depending on $n$. Clearly, the same bound (\ref{cte4}) holds for ${\cal V}_{n22}$.

It remains to estimate ${\cal V}_{n12}$. For $i,j \in\{0,\cdots,q\}$, we have
$$
\displaystyle \big |\frac {j -i} n\big |^{2H(\frac k n)+\delta_i+\delta_j}\ = \ \big |\frac {j-i} n\big |^{2H(\frac k n)} \exp\Big \{(\delta_i+\delta_j)\log \big |\frac {j-i} n\big | \Big \}\ = \ \big |\frac {j-i} n\big |^{2H(\frac k n)} + O\big (\frac {\log n} {n^{2H(\frac k n)+\eta\wedge 1}} \big ).
$$
Here and in the rest of the proof, we denote $O(\cdot)$ a bounded function depending on $H(\cdot)$, $\varepsilon$ and $q$ only, and not depending on $k$ or $n$. Next,
since $(H,H') \mapsto L_{12}(H,H')$ is analytic for $H +H'\neq 1$, so
$$
L_{12}\Big ( H(\frac k n)+\delta_i\, ,\,  H(\frac k n)+\delta_j\Big ) = 1 + O\big (\frac {1} {n^{\eta\wedge 1}} \big ).
$$
Therefore, 
\begin{eqnarray}\label{expL12}
{\cal V}_{n12}(k)
&=&C(H(k/n),a) \, n^{-2H(\frac k n)}
 + O\big (\frac {\log n} {n^{2H(\frac k n)+\eta\wedge 1}} \big )
\end{eqnarray}
where
$C(H,a) = -  n^{2H} \sum_{i,j=0}^q a_i \, a_j\, \ \big |\frac {j-i} n\big |^{2H}
 =  2 \Var \big (  V^a B_{H}(0)\big )$ 
does not depend on $n$ and does not vanish for $H \in (0,1)$, see Sec. 2.2. 
~\\
Using \eqref{cte4} and  \eqref{expL12}, we obtain:
\begin{eqnarray*}
{\cal V}_n(k) 
&=&-\frac 1 2 \, C(H(k/n),a) \, n^{-2H(\frac k n)}  + O\Big ( \frac {1} {n^{2(\eta\wedge m)}}\Big )   + O\big (\frac {\log n} {n^{2H(\frac k n)+\eta\wedge 1}} \big )\\
&=&\Var \big ( V_n^a B_{H(\frac k n)}(\frac k n)\big) \Big ( 1 + O\big (\frac {\log n} {n^{\eta\wedge 1}} \big ) +O\Big ( \frac {1} {n^{2(\eta \wedge m) -2H(\frac k n)}}\Big )  \Big )  .
\end{eqnarray*}
Hence the bound of Proposition \ref{propMBM1} can be deduced. \hfill $\Box$\\
~\\

\noindent {\it Proof  of Corollary \ref{corolMBM}.} Fact (i) is immediate from Proposition \ref{propMBM1}. By Property  \ref{popy1},
$\E (V_n^{a^{*}}X_{(a^+,a^-)}(k/n))^2 \sim (4- 4^{H(k/n)}) n^{-2H(k/n)}    $ and hence $\E (V_n^{a^{*}} X_{(a^+,a^-)}(k/n))^2 \ge C n^{-2H(k/n)} $ for some $C>0$ independent of  $1\le k \le n-2$ and $n$.  Use the last fact,
Proposition \ref{propMBM2} (\ref{lexp}), (\ref{lexp12}), boundedness of $A(t,t') $ and $B(t,t'),$ and $|V_n^{a^*} H(\frac k n)| \le
C n^{- (\eta \wedge 2)} $ to conclude inequality (\ref{corMBM}) for $|k-k'| > k_0$ and $k_0$ large enough; for $|k-k'| \le k_0$
(\ref{corMBM}) is trivial since the l.h.s. of  (\ref{corMBM}) is bounded by $1$.

It remains to prove fact (ii). 
Let $k'- k =: \tilde k \ge 1, \,  H :=  H(k/n), \, H' := H(k'/n)$.  Let first $ \tilde k  \le n^{1-\alpha}$. Then in view of (\ref{corMBM}) it suffices to check that
\begin{equation} \label{iibdd}
\frac 1 {\tilde k^{4-H - H'}} \, \le \, \frac{C n^{\theta (1-\alpha)}}{\tilde k^\gamma}, \qquad
\frac {1}{n^{2(\eta \wedge 2)-H - H'}}\, \le \,  \frac{C n^{\theta (1-\alpha)}}{\tilde k^\gamma}.
\end{equation}
The second inequality in (\ref{iibdd}) holds for $\tilde k  \le n^{1-\alpha}$ because of (\ref{iineq}) and
$\min (2 (\eta \wedge 2) - 2H, 2 (\eta \wedge 2) - 2H') \le 2 (\eta \wedge 2) - H - H'$. The first inequality in (\ref{iibdd})
holds trivially since $\gamma - (4 - H - H') < 0.$
Next, let $\tilde k > n^{1-\alpha}$. Then using (\ref{corMBM}) it suffices to check the inequalities
$$
\frac 1 {\tilde k^{4-H - H'}}\, \le \, \frac{C} {n^{(\gamma - \theta) (1-\alpha)}}, \qquad
\frac {1}{n^{2(\eta \wedge 2)-H - H'}} \, \le \, \frac{C} {n^{(\gamma - \theta) (1-\alpha)}}.
$$
which easily follow from (\ref{iineq}), similarly as above. \hfill $\Box$\\
~\\
\noindent {\it Proof of Proposition \ref{limitMBM}..} Part (i) is immediate from Theorem \ref{Limitgeneral}  (i) and Corollary \ref{corolMBM} (i).
Part (ii) follows from Theorem \ref{Limitgeneral}  (ii) and Corollary \ref{corolMBM}, by taking
$\theta = 0$ and $\gamma > 1/2 $ sufficiently close to $1/2$ (in such case the inequality $\gamma <
4 - 2H(t)$ of (\ref{iineq}) is satisfied since $H(t) < 1 $). To prove (iii),  use Corollary \ref{corolMBM} and
Theorem \ref{Limitgeneral}  (iii) with $\theta = 0$ and $\gamma > 1/2 $ which can be chosen under \eqref{MBM1}, with  $\mu =  (2-3\alpha)/4 $ and condition (\ref{MBM2}). Note that the lower bound  (\ref{MBM2}) implies $(\eta \wedge 1)
> (2-3\alpha)/4 $ and \eqref{MBM1} implies $2\big ((\eta \wedge 2) - \sup_{t\in (0,1)}H(t) \big )> (2-3\alpha)/4$, or condition $\kappa \ge \mu $ of Theorem \ref{Limitgeneral}  (iii). Parts (iv) and (v) follows similarly from Corollary \ref{corolMBM},
Theorem \ref{Limitgeneral}  (iv) (with $\theta = 0$, $\kappa\geq \mu_1=(1-3\alpha)/4$ satisfied from \eqref{MBM1}) and Theorem \ref{Limitgeneral2} (with $\alpha\geq 1 -2\kappa$ and $(1+2(\eta\wedge 2))^{-1}\geq 1-2(\eta\wedge 1)$ for all $\eta>0$). \hfill $\Box$\\
~\\

{\bf Acknowledgments.} The authors are extremely grateful to an
anonymous referee for a very careful reading and numerous relevant
suggestions and criticisms that strongly improved
the first version of the paper.

\end{document}